\documentclass{ucthesis}

\usepackage{amsfonts}
\usepackage{amsmath}
\usepackage{amsthm} 
\usepackage{epsfig}
\usepackage{graphics}
\usepackage{natbib}
\usepackage{latexsym}
\usepackage{amssymb,exscale}
\usepackage[all]{xy}

\usepackage{bookman}

\theoremstyle{plain} 
    \newtheorem{theorem}{Theorem}
    \newtheorem{lemma}[theorem]{Lemma}
    \newtheorem{proposition}[theorem]{Proposition}

    \newtheorem{conjecture}[theorem]{Conjecture}
\theoremstyle{definition} 
    \newtheorem{definition}[theorem]{Definition}
    
    \newtheorem{remark}[theorem]{Remark}
    \newtheorem{example}[theorem]{Example}

\def\Gam{\Gamma}
\def\Del{\Delta}

\def\Ome{\Omega}
\def\alp{\alpha}
\def\bet{\beta}
\def\gam{\gamma}
\def\del{\delta}
\def\eps{\epsilon}
\def\lam{\lambda}

\def\Ome{\Omega}

\def\ups{\upsilon}
\def\zet{\zeta}

\def\C{\mathbb{C}}
\def\F{{\bf F}}
\def\G{{\bf{G}}} 
\def\H{{\bf{H}}} 

\def\Z{\mathbb{Z}}
\def\M{\Omega}
\def\N{\mathbb{N}}
\def\Q{{\mathcal Q}}
\def\R{\mathbb{R}}

\def\X{{\mathcal X}}

\def\Z{\mathbb{Z}}
\def\ZZ{{\mathcal Z}}

\def\f{{\bf f}}
\def\g{{\bf g}}

\def\m{{\bf m}}
\def\n{{\bf n}}

\def\summ{\sum\limits}
\def\intt{\int\limits}
\def\prodd{\prod\limits}
\def\limm{\lim\limits}
\def\tends{\rightarrow}

\def\fr{\frac}
\def\Mid{\left\vert \right.}
\def\l{\left}
\def\r{\right}
\def\<{\langle}
\def\>{\rangle}
\def\hsp{\hspace}
\def\mb{\mbox}

\newcommand{\E}{\mbox{\bf E}}

\def\bar{\overline}
\def\P{{\bf P}}
\def\Var{\mbox{Var}}

\def\eqd{\stackrel{d}{=}}
\def\given{\left.\vphantom{\hbox{\Large (}}\right|}
\def\ind{{\mathbf 1}}
\def\d{\partial}
\def\dz{\frac{\partial }{ \partial z}}
\def\dzbar{\frac{\partial }{ \partial {\bar z}}}
\def\dw{\frac{\partial }{ \partial w}}
\def\dwbar{\frac{\partial }{ \partial {\bar zw}}}

\def\lap{\Delta}

\newcommand\mnote[1]{} 
\newcommand\be{\begin{equation}}

\newcommand\ee{\end{equation}}

\newcommand{\comment}[1]{}

\newcommand{\sm}{{\raise0.3ex\hbox{$\scriptstyle \setminus$}}}

\newcommand{\D}{{\mathbb D}}

\def\mb{\mbox}
\def\l{\left}
\def\r{\right}

\def\lam{\lambda}
\def\alp{\alpha}
\def\eps{\epsilon}
\def\tends{\rightarrow}

\newcommand{\T}{\mathbb T}

\numberwithin{equation}{section}
\numberwithin{theorem}{section}

\swapnumbers

\theoremstyle{plain}

\renewcommand{\phi}{\varphi}

\renewcommand{\S}{\mathbb S}

\renewcommand{\phi}{\varphi}

\def\Kdet{{\mathbb K}}

\begin{document} 
 \bibliographystyle{habbrv}

\title{Zeros of Random Analytic Functions}
\author{Manjunath Krishnapur}

\degreeyear{2006}
\degreesemester{Spring}
\degree{Doctor of Philosophy}
\chair{Professor Yuval Peres}
\othermembers{
Professor Michael Christ\\
Professor Steven N. Evans\\
Professor Jim Pitman}

\numberofmembers{4}
\prevdegrees{M.\ Stat.\ (Indian Statistical Institute, Calcutta, India) 2001}
\field{Statistics}
\campus{Berkeley}
 
\maketitle

\approvalpage

\copyrightpage
 
\begin{abstract}
The dominant theme of this thesis is that {\it random matrix valued analytic functions}, generalizing both random matrices and
random analytic functions, for many purposes can 
 (and perhaps should) be effectively studied in that level of generality. 

We study zeros of random analytic functions in one complex
variable. It is known that there is a one parameter family of Gaussian
analytic functions  with zero sets that are stationary  in each of the
three symmetric spaces, namely the plane, the
sphere and the unit disk, under the corresponding group of isometries. 

We show a way to generate non Gaussian random analytic functions whose
zero sets are also stationary in the same domains. There are particular cases
where the exact distribution of the zero set turns out to belong to an
important class of point processes known as determinantal point
processes.

Apart from questions regarding the exact distribution of 
zero sets, we also study certain asymptotic properties. We
 show asymptotic normality for smooth statistics applied to zeros of these
 random analytic functions. Lastly, we present some results on 
 certain large deviation problems  for the zeros of the planar and
 hyperbolic  Gaussian analytic functions.

\abstractsignature
\end{abstract}
 
\begin{frontmatter}
 
\begin{dedication}
\null\vfil
{\large
\begin{center}
To all my teachers. \\
\vspace{12pt}
\end{center}}
\vfil\null
\end{dedication}
 
\tableofcontents

\begin{acknowledgements}
It is a truth universally acknowledged, that a single author 
 presenting a piece of research must owe a lot to the inputs of many
 people. This definitely applies to this thesis, and I am happy
 to acknowledge the help I have received from various quarters.

 As regards my Ph.D., above all, I would like to thank my advisor,
 Yuval Peres, for advising me all these years and especially for
 keeping up an endless supply of problems to work on, when in a quest
 for research problems I was still
 executing a random walk on many topics in probability theory. I am just as 
 grateful to B\'{a}lint Vir\'{a}g and Yuval Peres for drawing me into
 the study of random analytic functions and for generously sharing 
 their insights and problems with me. I met random matrices through
 Steve Evans, whose course on this subject in my first semester at
 Berkeley remains one  of the most useful courses I have ever
 taken. Mikhail Sodin, through  his papers as well as through direct
 conversations, deeply  
  influenced my perspective of the subject. I thank my thesis committee
 members for reading the draft and making several useful suggestions.

 I have also immeasurably benefited from  the many courses I took,
 most of them in the statistics and mathematics departments. I learnt a
 great deal more of mathematics in these courses than I could have on
 my own, for which I thank all the instructors. Also
 I thank my fellow graduate students with whom I had many wonderful  
 discussions. Particularly, in the beginning years I learnt 
 significantly from Antar Bandyopadhyay, Noam Berger and G\'{a}bor
 Pete, and on the main topic of this thesis I had many great
 discussions with Ron Peled and my collaborator Ben Hough. I thank
 Sourav Chatterjee and Mikhail Sodin for getting me interested in
 Normal approximation problems. I am very
 grateful also to the department staff who have been immensely helpful
 throughout. 

While I learnt a  lot of mathematics after coming to Berkeley, none of 
 this would have  been possible without 
 the training in mathematics and statistics that I received from my
 professors at the Indian Statistical Institute. Without their dedicated
 teaching I would not  have become a probabilist. I particularly
 accuse Alladi Sitaram, Sundaram Thangavelu and S.M. Srivastava of 
 inspiring me into thinking that I should become a researcher. No less
 important were the courses of S.C.Bagchi, Arup Bose,  V. Pati,
 S. Ramasubramaniam, T.S.S.R.K. Rao and many others.  
  
Impersonal teachers have their own significance, and include for
example, the authors of many books, whose ideas I may have absorbed
and today assume to be my own. Among them I would particularly like to
mention the great 
probabilist and expositor, Mark Kac. When I applied for my PhD I
mentioned in my ``Statement of Purpose''  that what got me interested
in probability theory for the first time, was the theorem by Kac that
the average number of real roots of a random polynomial with
independent standard normal coefficients is asymptotic to $\frac{2}{\pi}\log
n$. It is a pleasure to me that this thesis can be seen 
as continuing the same theme of zeros of random polynomials, but
hopefully does more than fill some much needed gaps! 

 Even more influential on my development as a
 person and on my   
 attitude towards learning, were my parents, my teachers
 at school and college, and my relatives and friends (the public
 opinion is that I could stand much more development, but I claim for myself
 all credit for that). It would be silly to  
 even try to adequately acknowledge in words, their roles in my
 making. My father's huge r\'{e}pertoire of stories fired my imagination in my
 early childhood and made me think beyond everyday concerns. Apart
 from many other things, my  mother saw to it that I paid  due 
  attention to my studies, till the time came when I
 realized that it was a pleasure. My brother and sister, my aunts and 
  uncles and my friends were quite as important  
  in shaping me. Almost none of them is a mathematician, but they have
 that high respect (even without full comprehension) for knowledge
 that is so widespread in India. I greatly value 
  my friends for their great company and for never quite
 giving up on me, even though I have always been most irregular in
 returning their e-mails or phone calls. I am sure that all of them
 will feel happy on seeing my thesis.  


\end{acknowledgements}
 
\end{frontmatter}


\chapter{Introduction}\label{chap:intro}
\section{Introductory remarks}
Random analytic functions on the one hand and random matrices on the
other are two well studied topics in probability theory and
mathematical physics. One of the chief interests to a probabilist in
these objects is the kind of point processes one gets, by taking the
set of zeros or the set of eigenvalues, as the case may be. Both these 
  kinds of point processes typically have the property of 
  ``repulsion'', meaning that the points distribute themselves more evenly
  than they would if they were thrown down independently. That is an appealing
  feature because, while there are ways to construct point processes that
  are more clumped than independent points, there are not many natural
  ways in which to generate point processes with less clumping. 

This fact and several others (more empirical than mathematical) have
  led to a folk wisdom that random analytic functions  and random
  matrices  share many similarities. Differing responses to
  this statement have been heard, including one that points out the
  obvious tautology 
  here (after all the characteristic polynomial of a random matrix is 
  a random polynomial), and another that says that there is  not much similarity but instead evokes an ``anthropic'' reasoning (the same set of
  people work on both these fields). Without denying the validity of
  these explanations, in this thesis we take a more positive approach
  attempting to provide a unifying framework that includes both random
  matrices and random analytic functions (see caveats below).

\noindent {\em This is the simple but seemingly useful idea of
  considering random matrix-valued analytic functions, and the set of
  points 
  where it becomes singular (i.e., the zeros of the determinant). The
  linear polynomials reduce to random matrices and the $1\times 1$ matrices
  correspond to random analytic functions.}

A little explanation is in order. When we talk of random analytic
functions, we tacitly mean that we are somehow specifying the
distribution of coefficients or some closely related quantities (otherwise
any random set of points would be the zeros of a random analytic
function). Furthermore  it is usually difficult to analyse a random
analytic function (especially to get exact properties) except in the
case of Gaussian 
coefficients. So the essence of the above paragraph is that the
determinant of a Gaussian matrix-valued analytic function is a
non-Gaussian analytic function in itself, but nevertheless amenable to
analysis because it is built out of Gaussian analytic functions.

Secondly, the earlier claim about random matrices falling within our
framework should be toned down. The chief, although not
the whole, emphasis in random matrix theory is on the study of  
Hermitian random matrices and their (real) eigenvalues, for physical
as well as mathematical reasons. When we go to
higher polynomials there is perhaps no natural way to get the zeros 
 to lie on the real line. This may explain why these objects have not
 been studied before.  What we study here are zeros in the
 complex plane, for which of course there is no such
 problem. Nevertheless 
 we believe that it is also interesting mathematically to study polynomials
 with random Hermitian or random unitary coefficients (we do not do
 this here) even though the zeros are spread out in the complex plane.

We now outline the contents of the thesis briefly.
\begin{itemize}
\item In the remaining sections of this chapter, we give a quick
introduction to the basic notions of a point process, correlation
functions and Gaussian analytic functions. Most importantly, we recall
the three canonical families of Gaussian analytic functions on the plane, the
sphere and the unit disk (hyperbolic plane). 

\item In Chapter~\ref{chap:invzeros} we give a recipe for generating
 a slew of (non-Gaussian) random analytic functions whose zeros are
stationary in the  plane, the sphere and the unit disk. We make some
basic computations on the distribution of zeros that will be used later.

\item In Chapter~\ref{chap:invdets} we recall the notion of a determinantal
point process, and characterize the stationary determinantal point
processes in the three fundamental domains. Of these the planar ones
are known to be (limits of) the distribution of eigenvalues of certain
random matrices (the Ginibre ensemble) while the processes on the
sphere and disk are new (these processes themselves have been considered
before in \cite{caillol}, but an independent probabilistic
meaning was not known).

\item In Chapter~\ref{chap:matgaf}, we present an evocative analogy which 
suggests that the determinantal point processes on the sphere and the
disk, introduced in Chapter~\ref{chap:invdets},  are in fact the
singular points of certain random matrix-valued analytic
functions that were introduced in Chapter~\ref{chap:invzeros}. 

\item In Chapter~\ref{chap:sphere} we prove that the stationary
determinantal processes on the sphere introduced
in Chapter~\ref{chap:invdets} are the singular points of the
random matrix analytic function $zA-B$ (in this case, more simply, the
eigenvalues of $A^{-1}B$).

\item In Chapter~\ref{chap:disk} we give partial proof that the
determinantal processes on the disk introduced in
Chapter~\ref{chap:invdets} are the singular points of the random matrix
analytic function $A_0+zA_1+z^2A_2+\ldots $.  

\item In Chapter~\ref{chap:normality} we show asymptotic normality for
smooth statistics applied to the zeros of random analytic functions
introduced in Chapter~\ref{chap:invzeros} following a method of Sodin
and Tsirelson who showed the same for the canonical models of Gaussian
analytic functions. 

\item In Chapter~\ref{chap:overcrowd} and Chapter~\ref{chap:moderate} we move
away from the line of presentation so far, and return to canonical {\it
  Gaussian} analytic functions. We deal with two large deviation type
problems for zeros of the planar Gaussian analytic function, one posed by
Yuval Peres, which we solve fully and another due to Mikhail Sodin, which we solve partially.
\end{itemize}

\section{Basic notions and definitions}
\subsection{Point processes, Correlation functions}\label{sec:corrfun}
A {\bf point process} in a locally compact Polish space
$\M$ is a random integer-valued positive
Radon measure $\X$ on $\M$. (Recall that a
 Radon measure is a Borel measure which is finite
on compact sets.) If $\X$ almost surely assigns at most measure $1$
to singletons, it is a {\bf simple} point process;
in this case $\X$ can be identified with a random discrete subset
of $\M$, and $\X(D)$ represents the  number of points of
this set that fall in $D$.

The distribution of a point process can, in most
cases, be described by its {\bf correlation functions}
(also known as joint intensities) w.r.t a fixed Radon measure $\mu$ on
$\M$. 
\begin{definition}\label{def:jtint} The correlation functions of a point
  process $\X$  w.r.t.~$\mu$ are functions (if any exist) $\rho_n:\M^n\tends
  [0,\infty)$ for $n\ge 1$, such that for any family
  of  mutually disjoint Borel subsets $D_1,\ldots,D_k$ of $\M$, and for
  any non-negative integers $n_1,\ldots ,n_k,$
\begin{equation}\label{eq:corrfun1}
  \E\l[\prod_{i=1}^k \genfrac{(}{)}{0pt}{}{\X(D_i)}{n_i} n_i!\r]=
  \intt_{\prodd_i D_i^{n_i}} \rho_n(x_1,\ldots ,x_n) d\mu(x_1)\ldots
  d\mu(x_n),
\end{equation}
where $n=\summ_{i=1}^k n_i$.
\end{definition}
\begin{remark} It is a natural question to ask for conditions that
  guarantee the existence of correlation functions and conditions
  under which they determine the distribution of the point
  process. Such conditions do exist, see Lenard's 
  ~\cite{lenard1},\cite{lenard2},\cite{lenard3} or the survey by
  Soshnikov~\cite{sos1}. But the conditions are too complicated and
  not relevant for our purposes. In any case, when the joint distribution of
 $\X(D_1),\ldots,\X(D_k)$ is determined by its moments, the
  correlation functions determine the distribution of $\X$.
\end{remark}

\begin{remark}\label{rem:jtint}
For overlapping sets, the situation is more
complicated. Restricting attention to simple point processes,
 $\rho_n$ is not the intensity measure of
$\X^n$, but that of $\X^{\wedge n}$, the set of
ordered $n$-tuples of distinct points of $\X$. Indeed,
(\ref{eq:corrfun1}) implies (see \cite{lenard1,
lenard2, pervir}) that for any Borel set $B\subset
\M^n$ we have
 \begin{equation}\label{soshcor1} \E\,
\#(B\cap \X^{\wedge n}) =\intt_B \rho_n(x_1,\dots
,x_n)\,d\mu(x_1)\dots d\mu(x_n) \,. 
\end{equation}
\end{remark}

Assuming that $\X$ is simple, the correlation functions may
be interpreted as follows: 
\begin{itemize}
\item If $\M$ is finite and $\mu = \textrm{counting measure}$ then
  $\rho_k(x_1, \dots, x_k)$ is the probability that $x_1, \dots, x_k
  \in \X$. 
\item If $\M$ is open in $\R^d$ and $\mu=\textrm{Lebesgue measure}$,
  if $\rho_n$ exist and are continuous, then
\begin{equation}\label{eq:equivalentdefinitionofintensities}
\rho_k(x_1, \dots, x_k) = \lim_{\epsilon \rightarrow 0} \frac{\P[\X
  \textrm{ has a point in each of
  }B_\eps(x_j)]}{(\textrm{Vol}(B_\eps))^k}. 
\end{equation}
Conversely, if for every $k\ge 1$, the right hand side of
(\ref{eq:equivalentdefinitionofintensities}) exists and is continuous
in $x_i$, $1\le i \le k$, then it is the $k$-point correlation functions of $\X$.
\end{itemize}

For us $\M$ will always be an open subset of the plane (or the sphere
$S^2$) and $\X$ will be a simple point process. $\mu$ may always be
taken to be the Lebesgue measure on $\M$, but we often find it
convenient to use some other measure that is mutually absolutely
continuous with the Lebesgue measure. 

\subsection{Complex Gaussian distribution}
A {\bf standard complex Gaussian} is a complex-valued random variable
with probability density $\frac{1}{\pi}e^{-|z|^2}$ w.r.t the Lebesgue
measure on the complex plane. Equivalently, one may define it as
$X+iY$, where $X$ and $Y$ are i.i.d.~N($0,\frac{1}{2}$) random
variables. 

Let $a_k$, $1\le k\le n$ be i.i.d. standard complex Gaussians. Let
${\bf a}$ denote the column vector $(a_1,\ldots ,a_n)^t$. Then if $B$ is an $m\times n$ matrix, 
$B{\bf a}+{\bf \mu}$ is said to be an $m$-dimensional complex Gaussian
vector with mean ${\bf \mu}$ (an $m\times 1$ vector) and covariance
$\Sigma=BB^*$ (an $m\times m$ matrix). We denote its distribution by
$\C N_m \l({\bf \mu},\Sigma \r)$. 

Here are some basic properties of complex Gaussian random variables.
\begin{itemize}
\item If ${\bf a}$ is a complex Gaussian, its distribution is
  determined by ${\bf \mu}=\E\l[{\bf a} \r]$ and $\Sigma=\E\l[\l({\bf
  a}-{\bf \mu}\r)\l({\bf a}-{\bf \mu}\r)^* \r]$. All moments of the
  form 
  \begin{equation*}
    \E\l[(a_k-\mu_k)(a_j-\mu_j) \r], \hsp{1cm} 1\le k,j \le n, 
  \end{equation*}
vanish. This is the case even for $j=k$.

\item If $a$ is a standard complex Gaussian, then $|a|^2$ and
  $\frac{a}{|a|}$ are independent, and have exponential distribution 
 with mean $1$ and uniform distribution on the circle $\{z:|z|=1\}$,
  respectively. 

\item If $a_n$, $n\ge 1$ are i.i.d.~$\C N(0,1)$, then by an easy
  application of Borel-Cantelli,
  \begin{equation}
    \label{eq:limsupofiid}
    \lim\sup_{n\tends \infty} |a_n|^{\frac{1}{n}} = 1, \hsp{1cm}
    \mb{almost surely.}
  \end{equation}
In fact, equation (\ref{eq:limsupofiid}) is valid for any
i.i.d.~sequence of complex-valued random variables $a_n$, such that
$\E \l[ \max \l\{\log|a_1|,0 \r\} \r] <\infty$. Equation 
(\ref{eq:limsupofiid}) is useful to compute the radii of convergence
of random power series with independent coefficients.

\item {\bf Wick Expansions:}  The Wick or the Feynman diagram
  expansion  is an expansion of $L^2$ functions of a Gaussian measure in an
orthonormal basis consisting of polynomials of the underlying
Gaussians. Following the presentation in the book by
Janson~\cite{janson}, we state the essential facts in the limited
  context that we shall need later.
 More details and complete proofs of the assertions can be found in
  \cite{janson}.

Let $a_1,\ldots ,a_p$ be i.i.d. standard complex Gaussians.
 Consider the collection of all  monomials $\prodd_k
a_k^{m_k}\bar{a}_k^{n_k}$ in these $p$ variables, and orthonormalise 
 them by projecting the polynomials of degree $N$ on the orthogonal
 complement of the polynomials of degree $N-1$. 

 This procedure is the same as applying
Gram-Schmidt to the monomials after arranging them in increasing order
of the degree (how we order monomials of the same degree is
immaterial because distinct monomials of the same degree are  clearly
orthogonal).  
Thus we get an orthonormal basis of all square integrable functions of the $a_k$s, and
 the basis elements, termed {\bf Wick powers},  are denoted by 
\begin{equation*}
:\prodd_k \fr{a_k^{m_k}\bar{a}_k^{n_k}}{\sqrt{m_k ! n_k !}}:= \prodd_k
:\fr{a_k^{m_k}\bar{a}_k^{n_k}}{\sqrt{m_k ! n_k !}}:, 
\end{equation*}
the equality a consequence of the independence of $a_k$s.

These Wick polynomials are known explicitly (see ~\cite{janson})- 
\begin{equation*}
  :a^m\bar{a}^n: = \summ_{r=0}^{m\wedge n} (-1)^r r !
  \genfrac{(}{)}{0pt}{}{m}{r} \genfrac{(}{)}{0pt}{}{n}{r} a^{m-r}\bar{a}^{n-r},
\end{equation*}
although this is not particularly important to us.
It is quite well known that products of random variables that are
jointly Gaussian can be described by summing over the weights of
certain combinatorial entities. There is a similar formula (known as
Wick formula or Feynman diagram formula) for expectation of product of
Wick powers. We shall only need the following special
case.

{\bf Wick/Feynman diagram formula:} Let $(b_1,\ldots ,b_s)$ have a complex
Gaussian distribution with mean zero. Then 
\begin{equation}\label{eq:feynmanformula}
\E\l[\prodd_{j=1}^s :b_j^{m_j}{\bar b_j}^{n_j}: \r] = \summ_{\gam}
\ups(\gam),  
\end{equation}
where the sum is over all {\it complete Feynman diagrams} $\gam$ {\it 
 without self  interaction} (henceforth we shall just say {\bf Feynman
 diagram}). To define this, consider a collection of
 $\summ_j(m_j+n_j)$ vertices
 with $m_j$  of the vertices labeled $j$ and $n_j$ of the vertices
 labeled ${\bar 
 j}$, for $1\le j\le s$. All the vertices labeled $j$ are also
 supposed to be distinguishable although we shall not introduce any
 more notation to distinguish them. Now, each $\gam$ is a matching of
 these vertices (a subgraph in which each vertex has degree
 $1$), such that each edge in $\gam$ connects a vertex labeled $i$ to
 a vertex labeled ${\bar j}$ for some $i\not=j$.   
   
The value $\ups(\gam)$ of the diagram is the product of the weights of
all the edges in $\gam$, and the weight of an edge joining a vertex
labeled $i$ to a vertex labeled ${\bar j}$ ($i\not= j$) is $\E\l[
b_i{\bar b_j}\r]$.   

\begin{example}\label{eg:wickfortwo}
 Let $s=2$. Then we must consider Feynman diagrams on
  the labels $\{ 1, \bar{1},2,\bar{2} \}$ with $m_1$ vertices labeled
  $1$, $m_2$ vertices labeled $2$, $n_1$ vertices labeled $\bar{1}$
  and $n_2$ vertices labeled $\bar{2}$. Since a Feynman diagram (in
  our terminology as explained above) must connect $1$s to $\bar{2}$s
  and vice-versa, and must give every vertex degree one, there are no
  Feynman diagrams unless $m_1=n_2$ and $m_2=n_1$, in which case there
  are $m_1!n_1!$ such diagrams. Thus
\begin{equation*}
 \E\l[\prodd_{j=1}^2 :b_j^{m_j}{\bar b_j}^{n_j}: \r] = \l\{
\begin{array}{ll}
 m_1!n_1!\E[b_1\bar{b}_2]^{m_1}\E[b_2\bar{b}_1]^{n_1} & \mb{ if
 }m_1=n_2,m_2=n_1. \\
 0 & \mb{ otherwise.} \end{array} \r.
\end{equation*}
\end{example}

\end{itemize}

\subsection{Gaussian analytic functions} 
Endow  the space of analytic functions on a region $\M$ with
 the topology of uniform convergence on compact sets.
  This makes it a complete separable metric space which is the
  standard setting for doing probability theory (To see completeness, 
  if $\{f_n\}$ is a Cauchy sequence, then $f_n$ converges uniformly on
 compact sets to some continuous function $f$. Then it is easy to see
 that $f$ must be analytic because its integral on any closed contour
 is zero since $\intt_{\gam} f = \limm_{n\tends \infty}\intt_{\gam}
 f_n $ and the 
 latter vanishes for every $n$, by analyticity of $f_n$).  

\begin{definition} Let $\f$ be a random variable taking values in the
 space of analytic functions on a region 
  $\M\subset \C$. We say $\f$ is a Gaussian analytic function
  (GAF) on $\M$ if  $(\f(z_1),\ldots,\f(z_n))$ has a mean zero
  complex Gaussian distribution for every $z_1,\ldots,z_n\in \M$.     
\end{definition}
It is easy to see the following properties of GAFs.
\begin{itemize} 
\item $\{\f^{(k)}\}$ are jointly Gaussian, i.e., the joint
  distribution of $\f$ and finitely many derivatives of $\f$ at
  finitely many points,
  \begin{equation*}
    \l\{\f^{(k)}(z_j):0\le k\le n, 1\le j\le m\r\},
  \end{equation*}
has a (mean zero) complex Gaussian distribution. 
\item The distribution of a Gaussian analytic function is determined
  by its covariance kernel $\l(\E\l[\f(z)\bar{\f(w)}\r] \r)_{z,w\in
  \M}$ denoted by $K_\f(z,w)$ or just $K(z,w)$ if there is no
  ambiguity as to which $\f$ is under consideration.
\end{itemize}

\subsection{Stationary zero sets of Gaussian analytic
  functions}\label{subsec:invgafs} 
 Our interest is in the zero set of a random analytic function. Unless
 one's intention is to model a particular physical phenomenon by a 
point process, there is one criterion that makes some point processes
more interesting than others, namely, {\it stationarity} under a large group
of transformations (stationarity of a random process means invariance of its 
distribution under a group action. It is also called {\it invariance},
 especially when the stationarity is in ``space'' rather than
 ``time'', but we use both terms interchangeably). There are three 
particular two dimensional domains on which the group of conformal
automorphisms act transitively (There are two others that we do not
consider here, the cylinder or the punctured plane, and the two
dimensional torus). We introduce these domains now.

\begin{itemize}
\item {\bf The Complex Plane $\C$}: The group of transformations
  \begin{equation}\label{eq:isometriesofplane}
    \phi_{\lam,\beta} (z) = \lam z +\bet, \hsp{1cm} z\in \C
  \end{equation}
where $|\lam|=1$ and $\bet\in \C$, is nothing but the Euclidean motion group.
These transformations preserve the Euclidean metric
$ds^2=dx^2+dy^2$  and the Lebesgue measure
$dm(z)=dx dy$ on the plane. 

\item {\bf The Sphere $\S^2=\C\cup \{\infty\}$}:  The group of transformations
  \begin{equation}\label{eq:isometriesofsphere}
    \phi_{\alp,\beta} (z) = \frac{\alp z +\bet}{-{\bar \bet}z + {\bar
    \alp}}, \hsp{2cm} z\in \C\cup \{\infty\}
  \end{equation}
where $\alp,\bet\in \C$ and $|\alp|^2+|\bet|^2=1$, is the group of
linear fractional transformations mapping $\C\cup \{\infty\}$ to itself
bijectively. These transformations preserve the spherical metric
$ds^2=\frac{dx^2+dy^2}{(1+|z|^2)^2}$ and the spherical area measure
$\frac{dm(z)}{(1+|z|^2)^2}$.  We call it the spherical metric because
it is the push forward of the usual metric on the sphere inherited
from $\R^3$, onto $\C\cup \{\infty\}$ under the stereographic
projection, and the measure is the push forward of the spherical area
measure. The transformations (\ref{eq:isometriesofsphere}) are just
the rotations of the sphere under this identification with $\C \cup
\{\infty\}$. 

\item {\bf The Hyperbolic Plane $\D=\{z: |z|<1 \}$}: 
 The group of transformations
  \begin{equation}\label{eq:isometriesofdisk}
    \phi_{\alp,\beta} (z) = \frac{\alp z +\bet}{{\bar \bet}z + {\bar
    \alp}}, \hsp{2cm} z\in \D
  \end{equation}
where $\alp,\bet\in \C$ and $|\alp|^2-|\bet|^2=1$, is the group of
linear fractional transformations mapping the unit disk
$\D=\{z:|z|<1\}$ to itself  bijectively. These transformations
preserve the hyperbolic metric $ds^2=\frac{dx^2+dy^2}{(1-|z|^2)^2}$
and the hyperbolic area measure $\frac{dm(z)}{(1-|z|^2)^2}$ (this
normalization differs from the usual one, with curvature $-1$, by a
factor of $4$, but it makes the analogy with the other two cases more
formally similar). This is one of the many models discovered by
Poincar\'{e} for the hyperbolic geometry of Bolyai, Gauss and
Lobachevsky (see~\cite{cfkp} for an introduction).
\end{itemize}

Note that in each case, the group of transformations acts transitively
on the corresponding space, i.e., for every $z,w$ in the domain, there
is a transformation $\phi$ such that $\phi(z)=w$. This means that in
these spaces every point is just like every other point. Now we
introduce three families of GAFs whose relation to these symmetric
spaces will be made clear in Proposition~\ref{prop:stationarygafs}.

In each case, the domain of the random analytic function can be found
 from (\ref{eq:limsupofiid}). Indeed, (\ref{eq:limsupofiid}) implies
 that when $a_n$ are i.i.d. standard complex Gaussians, $\summ_n a_n
 c_n z^n$ has the same radius of convergence as $\summ_n c_n z^n$.
\begin{itemize}
\item {\bf The Complex Plane $\C$}: Define for $L>0$,
  \begin{equation}\label{eq:planargaf}
    \f(z) = \summ_{n=0}^{\infty} a_n \frac{\sqrt{L^n}}{\sqrt{n!}}z^n.
  \end{equation}
For every  $L>0$, this is a random analytic function in the entire plane.

\item {\bf The Sphere $\S^2$}: Define for $L\in \N=\l\{1,2,3,\ldots  \r\}$,
  \begin{equation}\label{eq:sphericalgaf}
    \f(z) = \summ_{n=0}^{L} a_n \frac{\sqrt{L(L-1)\ldots
    (L-n+1)}}{\sqrt{n!}}z^n. 
  \end{equation}
For every  $L\in \N$, this is a random analytic function on
$\S^2=\C\cup {\infty}$ with a pole at $\infty$ (i.e., it is a polynomial). 

\item {\bf The Hyperbolic Plane $\D$}: Define for $L>0$,
  \begin{equation}\label{eq:hyperbolicgaf}
 \f(z) = \summ_{n=0}^{\infty} a_n \frac{\sqrt{L(L+1)\ldots
    (L+n-1)}}{\sqrt{n!}}z^n. 
  \end{equation}
For every  $L>0$, this is a random analytic function in the unit disk 
$\D=\{z:|z|<1\}$. 
\end{itemize}

\begin{remark}
  Although we wrote (\ref{eq:planargaf}) for every $L>0$, they are identical up to a scaling
  of the complex plane. However, the functions in
  (\ref{eq:sphericalgaf}) and (\ref{eq:hyperbolicgaf})  are truly
  different for different $L$, i.e., there is no transformation of the
  $S^2$ and $\D$, that makes $\f_L$ and $\f_{L'}$ the same, for
  $L\not= L'$. This is particularly obvious for the sphere, because
  $L$ then denotes the number of zeros of $\f$.
\end{remark}

We just quote the following proposition from \cite{ST1}. (The proof is
 contained in  the proof of
 Proposition~\ref{prop:stationaryrafs}). These random analytic
 functions were discovered in several stages and (partially) by
 several authors. The main contributions are due to  Bogomolny, Bohigas and
  Leboeuf \cite{bbl92} and \cite{bbl96}, Kostlan~\cite{kostlan93}, 
  Shub and Smale~\cite{shubsmale}. Some of them are natural
  generalizations (to complex coefficients) of random polynomials
  studied by Mark Kac in his founding papers starting with
 \cite{kac}. The special case $L=2$, in the unit disk was 
  derived also by  Diaconis and Evans~\cite{diaeva} as the limit of the
  logarithmic derivative of characteristic  polynomials of random
  unitary matrices. The uniqueness in
 Proposition~\ref{prop:stationarygafs} also was perhaps known, but a
 much stronger form of uniqueness (that the first intensity of zeros of any
 Gaussian analytic function determines the distribution of the
 Gaussian analytic function itself, up to multiplication by arbitrary deterministic non-vanishing analytic functions) was found by Sodin~\cite{sodin}.
 \begin{proposition}\label{prop:stationarygafs}
  The zero sets of the GAF $\f$ in equations (\ref{eq:planargaf}),
  (\ref{eq:sphericalgaf}) and (\ref{eq:hyperbolicgaf}) are invariant
  (in distribution) under the transformations defined in equations
  (\ref{eq:isometriesofplane}), (\ref{eq:isometriesofsphere}) and
  (\ref{eq:isometriesofdisk}) respectively. This holds for every
  allowed value of the parameter $L$, namely $L>0$ for the plane and
  the disk and  $L\in \N$ for the sphere. 

  Moreover, these are the only Gaussian analytic functions (up to
  multiplication by deterministic non vanishing analytic functions)
  with stationary zero sets in these domains.
\end{proposition}


\chapter{Stationary zero sets of random analytic
functions}\label{chap:invzeros}   

As we saw in Proposition~\ref{prop:stationarygafs}, on each of the three domains
$\C/\S^2/\D$, there is a one parameter family of Gaussian analytic
functions whose zero sets are stationary under the corresponding group
of isometries. Moreover, these are the only Gaussian analytic
functions on these domains with these properties. Indeed
Hannay~\cite{hannay} likens the uniqueness of the Gaussian analytic
function in (\ref{eq:planargaf}) to that of the Poisson process or the
thermal blackbody radiation. 


 Here we stick to the three domains  $\C/\S^2/\D$ and ask for random
 analytic functions whose zero sets are stationary. By
 Proposition~\ref{prop:stationarygafs}, we must necessarily seek among
 non-Gaussian analytic 
 functions. A natural idea might by to replace i.i.d. Gaussians in the
 coefficients by i.i.d. complex-valued random variables from some
 other distribution. However, these seem difficult to analyse. Gaussian
 analytic functions have the nice property that the evaluations of the
 function and its derivatives are all Gaussian with distributions that
 we can explicitly work with and this fails in other cases. In fact we
 do not know of another  
 example of a power series with i.i.d. coefficients whose zero set is
 stationary (on any of these three domains).  We resolve this deadlock by  
 constructing non-Gaussian analytic functions using Gaussian analytic
 functions as building blocks.

\section{A recipe for stationary zero sets of random
  analytic functions}
Let $\Q$ be a (non-random) homogeneous polynomial in $k$ variables
  with complex coefficients and let $\f$ be any Gaussian analytic
  function (not necessarily one of the canonical models defined in
  Section~\ref{subsec:invgafs}). Then if 
  $\f_i,i\le k$ are i.i.d. copies of $\f$, then $\Q(\f_1,\ldots
  ,\f_k)$ is a random analytic function on the same  domain as $\f$. 

\begin{proposition}\label{prop:stationaryrafs} Let $\Q$ be a
  homogeneous polynomial of degree $d$ 
  in $k$ variables 
  with complex coefficients, and let $\f$ be one of the canonical
  models of  Gaussian functions in (\ref{eq:planargaf}),
  (\ref{eq:sphericalgaf}) or (\ref{eq:hyperbolicgaf}). If $\f_i$, $1\le i \le
  k$ are i.i.d. copies of $\f$, then the zero set of the  random
  analytic function  
  \begin{equation*}
    \F(z) :=  \Q\l( \f_1(z),\ldots ,\f_k(z) \r)
  \end{equation*}
  is stationary under the same group of isometries as the zero set of $\f$.  
\end{proposition}

\begin{proof}
  First we recall the proof of invariance of the zero set of the
  Gaussian analytic functions in (\ref{eq:planargaf}),
  (\ref{eq:sphericalgaf}) and (\ref{eq:hyperbolicgaf}).. Fix
  an isometry $\phi$ of $\M$ (given in (\ref{eq:isometriesofplane}),
  (\ref{eq:isometriesofsphere}) and (\ref{eq:isometriesofdisk})). In
  each of the three cases, there is a 
  deterministic non-vanishing function $\Delta_{\phi,L}$ such that
  \begin{equation}\label{eq:shiftinggaf}
    \f(z) \eqd \Delta_{\phi,L}(z) \f\l(\phi(z)\r),
  \end{equation}
where in fact
\begin{equation*}
  \Delta_{\phi,L}(z) = \l\{ \begin{array}{cc}
        e^{Lz\lam \bar{\bet}+ \frac{1}{2} L|\bet|^2}  &
        \mb{domain}=\C. \\
        \phi'(z)^{\frac{L}{2}} & \mb{domain}=\S^2. \\
        \phi'(z)^{-\frac{L}{2}} &\mb{domain}=\D. \end{array} \r.
\end{equation*}
Note that the equality in (\ref{eq:shiftinggaf}) is for the entire
process, not just for a fixed $z$. Therefore, the zero set of $\f$ is
invariant in distribution under the action of $\phi$. (To prove equation
(\ref{eq:shiftinggaf}), just compute the covariance kernels of the
Gaussian processes on the left and right hand sides). 

Coming back to $\F$, we see that
\begin{eqnarray*}
  \Delta(z)^d \F(\phi(z)) &=& \Q\l( \Delta(z)\f_1\l(\phi(z)\r),\ldots
  ,\Delta(z)\f_k\l(\phi(z) \r) 
  \r) \\
 &\eqd& \Q\l(\f_1(z),\ldots ,\f_k(z) \r) \hsp{7mm} \mb{ from
   }(\ref{eq:shiftinggaf}) \\ 
  &=&  \F(z).
\end{eqnarray*}
This implies that the zero set of $\F$ is invariant in distribution
under the action of $\phi$.
\end{proof}

This is a very simple observation, but note that while $\F$ is built
 in a simple manner out of copies of $\f$, the zero set of $\F$ is by no means a
simple transformation of the zero sets of $\f_1,\ldots \f_k$ (except
in trivial cases such as when $\Q(\zet_1,\zet_2)=\zet_1 \zet_2$). Thus 
 the sets of zeros that we get are genuinely new point processes, but have
 the advantage of being based on {\it Gaussian} analytic functions,
 and therefore amenable to analysis. We illustrate this next, by
 computing the first and second correlations (joint intensities) for
 the zeros of $\F$. The tool that we use to study functions such as
 $\F$ is the Wick expansion, suggested to us
 by Mikhail Sodin (see the paper by Sodin and Tsirelson~\cite{ST1} for
 a use of Wick expansions in the context of Gaussian analytic functions).
We call random analytic functions of the kind described in
 Proposition~\ref{prop:stationaryrafs} as {\bf polygafs}.

\section{How to study the zeros of a polygaf?} 
If $F$ is any analytic function (not random) on $\M$, let $dn_F$ denote the
counting measure, with appropriate multiplicities, on the zeros of
$F$. Then, 
\begin{equation}
  \label{eq:lplogf}
  \frac{1}{2\pi}\lap \log
|F(z)| = dn_F(z)
\end{equation}
in the sense of distributions. This just means that for any $\phi\in
C_c^{\infty}(\M)$,  
 \begin{equation}\label{eq:distributionallaplacian}
   \intt_{\M} \phi(z) dn_{F}(z) =  \intt_{\M} \lap \phi(z)
   \frac{1}{2\pi} \log |F(z)| dm(z),
 \end{equation}
where $m$ is the Lebesgue measure. Therefore when $\F$ is any random
analytic function, understanding the distribution of the zero set
depends on being able to do computations with $\log |\F|$ (When $\F$ is Gaussian, there are
other approaches to studying the zero set of $\F$, but it appears that
the approach outlined here is the only one that is equally convenient
for our more general setting. The other methods make use of the
probability density of $\F$ evaluated at several points in the domain
etc, which are not available to us here). 

Now from (\ref{eq:distributionallaplacian}), if $\phi_i$, $1\le i \le
k$ are smooth functions with disjoint supports in $\M$, we get that
\begin{eqnarray*}
  \E\l[ \prodd_{i=1}^k \intt_{\M} \phi_i(z) dn_{\F}(z) \r] &=&
  (2\pi)^{-k}\intt_{\M^k}  \prodd_{i=1}^k \lap_{z_i} \phi_i(z_i) \E\l[
  \prodd_{i=1}^k \log |\F(z_i) | \r] \prodd_{i=1}^k dm(z_i) \\
 &=& (2\pi)^{-k}\intt_{\M^k} \l( \prodd_{i=1}^k \phi_i(z_i) 
   \lap_{z_i} \r) \E\l[ 
  \prodd_{i=1}^k \log |\F(z_i) | \r] \prodd_{i=1}^k dm(z_i).
\end{eqnarray*}
In the last line we integrated by parts.

In \ref{sec:corrfun} we defined the correlation functions in terms of
the moments of the joint counts of the number of points falling in
several regions. Fixing $k$ {\it distinct} points $w_1,\ldots ,w_k$ in
$\M$  and letting $\phi_i$ be a bump function in a small neighbourhood of
$z_i$, by elementary measure theoretical arguments one can deduce that the
 $k$-point correlation function $\rho_k$ of the zero set of $\F$,  
 w.r.t Lebesgue measure is given by 
\begin{equation}
  \label{eq:corrfunofraf}
  \rho_k(w_1,\ldots ,w_k) = \frac{1}{(2\pi)^k} \l( \prodd_{i=1}^k
  \lap_{w_i} \r) \E \l[  
  \prodd_{i=1}^k \log |\F(w_i) | \r],
\end{equation}
for distinct $w_1,\ldots ,w_k$. 
\begin{remark}
  Hammersley~\cite{hamm} gave a formula for the correlation functions
  of zeros of random polynomials in terms of the distribution of the
  coefficients. (\ref{eq:corrfunofraf}) is an alternative way of
  expressing the same. In this form, for $k=1$ it is sometimes called
  Edelman-Kostlan formula. 
\end{remark}



The way to analyse $\log |\F|$ is via Wick expansions that were
outlined in Chapter~\ref{chap:intro}.
\begin{example} The particular example of Wick expansions that is of
  interest to us is the 
  following: Let $\Q$ be a homogeneous polynomial in $k$ variables with complex 
coefficients. If $a_i,i\le k$ are i.i.d. $\C N(0,1)$ random
variables, then 
\begin{equation*}
\E\l[\given \log |\Q(a_1,\ldots,a_k)| \given^p\r]<\infty
\end{equation*}  
for every finite $p$. Hence we can expand $\log |\Q(a_1,\ldots,a_k)|$ in Wick powers as   
\begin{equation}
  \label{eq:wickforP}
  \log |\Q(a_1,\ldots ,a_k)| = \summ_{\m,\n \in \Z_+^k}
  \frac{C_{\m,\n}}{\sqrt{\m ! \n !}} \prodd_{j=1}^k :a_j^{m_j}
  \bar{a_j}^{n_j}:, 
\end{equation}
where $\m=(m_1,\ldots,m_k),\n=(n_1,\ldots,n_k)$, $\m!=\prodd_{j=1}^k
m_j!$ and  
\begin{equation*}
  C_{\m,\n} = \frac{1}{\sqrt{\m! \n!}}\E\l[\log |\Q(a_1,\ldots ,a_k)| \prodd_{j=1}^k :\bar{a_j}^{m_j} a_j^{n_j}:\r],
\end{equation*}
and the equality in (\ref{eq:wickforP}) is in the $L^2$ sense 
(it could be better, of course).

 We record two observations for later use.
\begin{itemize}
\item $C_{\m,\n}=\bar{C}_{\n,\m}$ for all $\m,\n \in \Z_+^k$, because
  $\bar{a}_k$ are also i.i.d. $\C N(0,1)$.
\item $C_{\m,\n}=0$ unless $\m_{\bullet}=\n_{\bullet}$, where
  $\m_{\bullet}:=\summ_j 
  m_j$. To see this, note for any $\lam$ with $|\lam|=1$, $\lam a_j$
  are also i.i.d. $\C N(0,1)$ and hence, by the homogeneity of $\Q$, it is
  also true that $\log |\Q(\lam a_1,\ldots ,\lam a_k)|=\log
  |\Q(a_1,\ldots ,a_k)|$. Therefore, from the equation above for
  $C_{\m,\n}$, we see that
  $C_{\m,\n}=\lam^{\m_{\bullet}-\n_{\bullet}}C_{\m,\n}$, which cannot
  be true unless $\m_{\bullet}=\n_{\bullet}$ or 
  $C_{\m,\n}=0$. 
\end{itemize}
\end{example}

\section{Distribution of the zero set of a polygaf}\label{sec:distribzeros}
Now let $\f_i$ be i.i.d. copies of
$\f$, a Gaussian analytic function on a domain $\Ome$ (not
necessarily one of the canonical GAFs on the plane, sphere or disk). As before $\Q$ is a homogeneous polynomial. 

Define $\F(z)=\Q(\f_1(z),\ldots ,\f_k(z))$. If $K(z,w)$ is the
covariance kernel of $\f$, then set
\begin{equation*}
  \hat{\F}(z) = \frac{\F(z)}{K(z,z)^{d/2}} =
  \Q\l(\hat{\f_1}(z),\ldots,\hat{\f_k}(z)\r), 
\end{equation*}
where $\hat{\f_j}(z)=\frac{\f_j(z)}{\sqrt{K(z,z)}}$ and $d$ is the degree of
$\Q$. Then from (\ref{eq:wickforP}) we can write,
\begin{equation*}
  \log |\hat{\F}(z)| = \summ_{\m,\n \in \Z_+^k}
  \frac{C_{\m,\n}}{\sqrt{\m! \n!}} \prodd_{j=1}^k :\hat{\f_j}(z)^{m_j}
  \bar{\hat{\f_j}(z)^{n_j}}:
\end{equation*}
We work with $\hat{\F}$ rather than $\F$, because $\hat{\f_j}(z)$ are
i.i.d standard Gaussians for any $z$, and so we can directly use the
 Wick formulas that we stated in the previous chapter.

\noindent {\bf First Intensity: } The first intensity (or $1$-point
  correlation function) as given by (\ref{eq:corrfunofraf}) is
\begin{equation}
  \label{eq:edelmankostlan}
  \rho_1(z)=\frac{1}{2\pi} \lap_z \E\l[\log |\F(z)| \r].
\end{equation}
Therefore,
\begin{eqnarray*}
  \rho_1(z) &=& \frac{1}{2\pi} \lap_z \E\l[\log |\F(z)| \r] \\
       &=& \frac{1}{2\pi} \lap_z \E\l[\log |\hat{\F}(z)|
       \r]+\frac{1}{2\pi} \lap_z \log K(z,z)^{d/2} \\
       &=& \frac{1}{2\pi} \lap_z C_{\underline{0},\underline{0}} +
       \frac{1}{2\pi} \lap_z \log K(z,z)^{d/2}
\end{eqnarray*}
from (\ref{eq:wickforP}). Therefore we obtain
\begin{equation}\label{eq:firstintensityforpolygafs}
     \rho_1(z) = \frac{d}{4\pi} \lap_z \log K(z,z).
\end{equation}
As a special case, set $k=1$, $\Q(\zeta)=\zeta$ in formula
(\ref{eq:firstintensityforpolygafs}) to deduce that the intensity of zeros of $\f$ is $\frac{1}{4\pi} \lap_z \log K(z,z)$ (This is known as the Edelman-Kostlan formula). Thus we see 
that the intensity of zeros of $\F$ is $d$ times the intensity of
zeroes of $\f$. This simple relationship between the intensities is
surprisingly not quite obvious from the definition. 

\noindent {\bf Two point Correlations: } Again from 
(\ref{eq:corrfunofraf}) we get the $2$-point correlation. It is easy
to see that
\begin{eqnarray*}
  \rho_2(z,w)-\rho_1(z)\rho_1(w) &=& \frac{1}{(2\pi)^2} \lap_z \lap_w \E\l[\log
  |\hat{\F}(z)| \log |\hat{\F}(w)|\r]
\end{eqnarray*}
From (\ref{eq:wickforP}) the right hand side can be written as
\begin{equation*}
     \frac{1}{(2\pi)^2} \lap_z \lap_w \summ_{\m,\n,\m',\n'}
  \frac{C_{\m,\n}C_{\m',\n'}}{\sqrt{\m! \n! \m'! \n'!}} \prodd_{j=1}^k
  \E\l[:\hat{\f_j}(z)^{m_j}
  \bar{\hat{\f_j}(z)^{n_j}}::\hat{\f_j}(z)^{m'_j}
  \bar{\hat{\f_j}(z)^{n'_j}}: \r].
\end{equation*}
This is precisely the situation elucidated in
Example~\ref{eg:wickfortwo}. Thus, only terms with $\m=\n',\n=\m'$
survive. 

Now make use of the observations made earlier-
  (1) $C_{\m,\n}=\bar{C}_{\n,\m}$, and (2) $C_{\m,\n}=0$ unless
  $\m_{\bullet}=\n_{\bullet}$. Grouping together terms by
  $p=\m_{\bullet}$, we get, 
\begin{equation}\label{eq:2ptcorrelationformula}
\rho_2(z,w)-\rho_1(z)\rho_1(w) = \frac{1}{(2\pi)^2}
  \summ_{p=0}^{\infty} 
  |\tilde{C}_p|^2 \lap_z \lap_w \frac{|K(z,w)|^{2p}}{K(z,z)^pK(w,w)^p},
\end{equation}
where $|\tilde{C}_p|^2= \summ_{\m_{\bullet}=\n_{\bullet}=p}
|C_{\m,\n}|^2$.

\begin{remark} Equation (\ref{eq:2ptcorrelationformula}) has the
  appealing feature that the effects of the two ingredients of $\F$,
  namely the polynomial $\Q$ and the Gaussian analytic function $\f$,
  are clearly separated. $|\tilde{C}_p|^2$ depends only on $\Q$ while 
  $\frac{|K(z,w)|^{2p}}{K(z,z)^pK(w,w)^p}$ depends only on $\f$.
  This observation is  crucially used in the next chapter, when we compute
  the correlations  for specific polynomials $\Q$. One can write
  analogous but more complicated  expressions for higher correlations, with $\lap_{z_k}$,
  $1\le k\le m$ applied to a sum over Feynman diagrams. 
\end{remark}


%


\chapter{Stationary determinantal point processes}\label{chap:invdets}
  
In this chapter we move away from random analytic functions
and talk about a different class of point processes. In the next
chapter, we return to zeros of random analytic functions and show that
there are point processes in the intersection of the two classes.

\begin{figure}\label{fig:threeprocesses}
\centering
\includegraphics[height=1.75in]{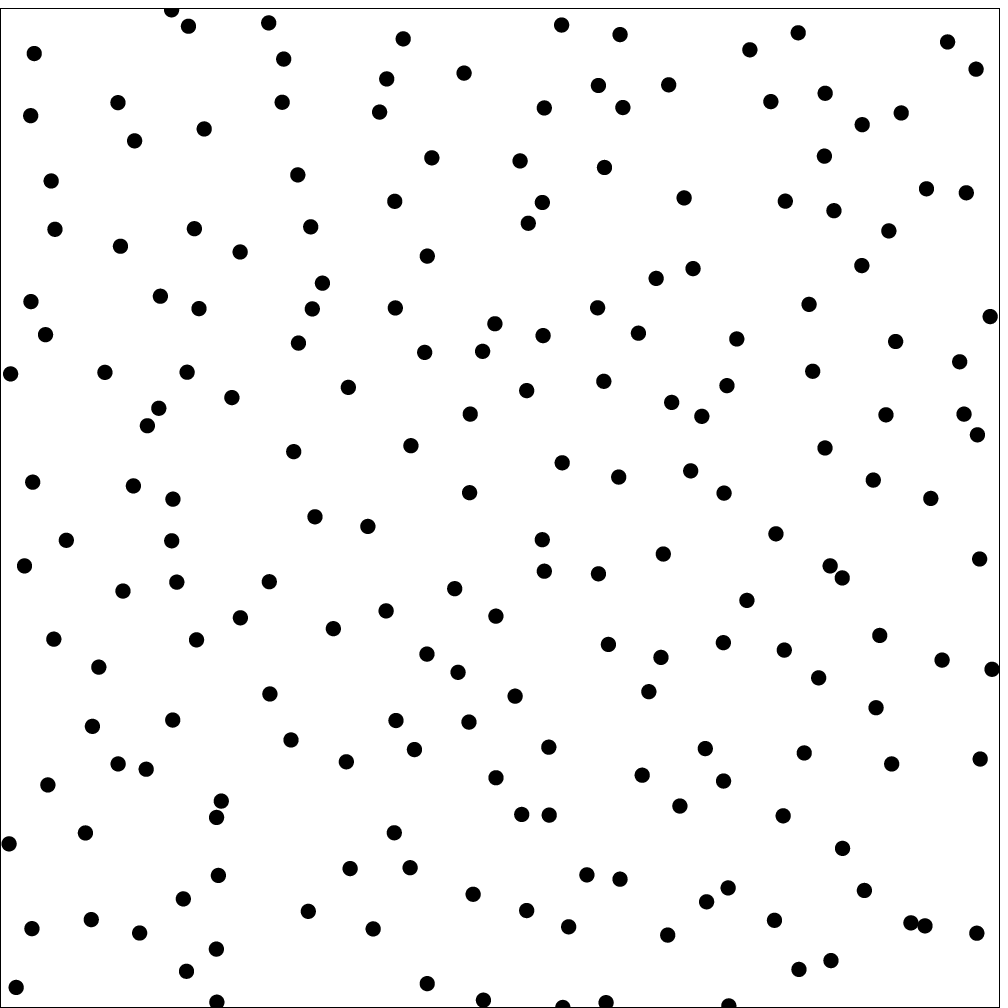}\hspace{.25in}
\includegraphics[height=1.75in]{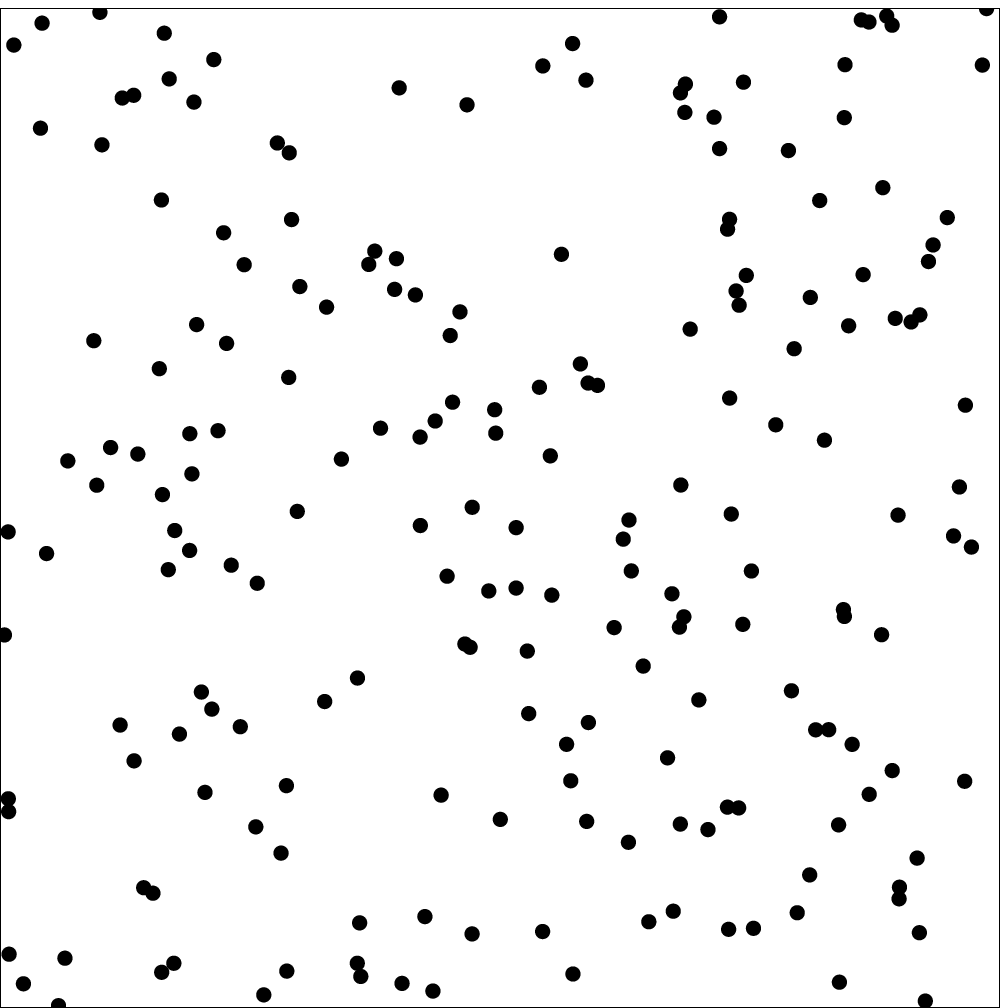}\hspace{.25in}
\includegraphics[height=1.75in]{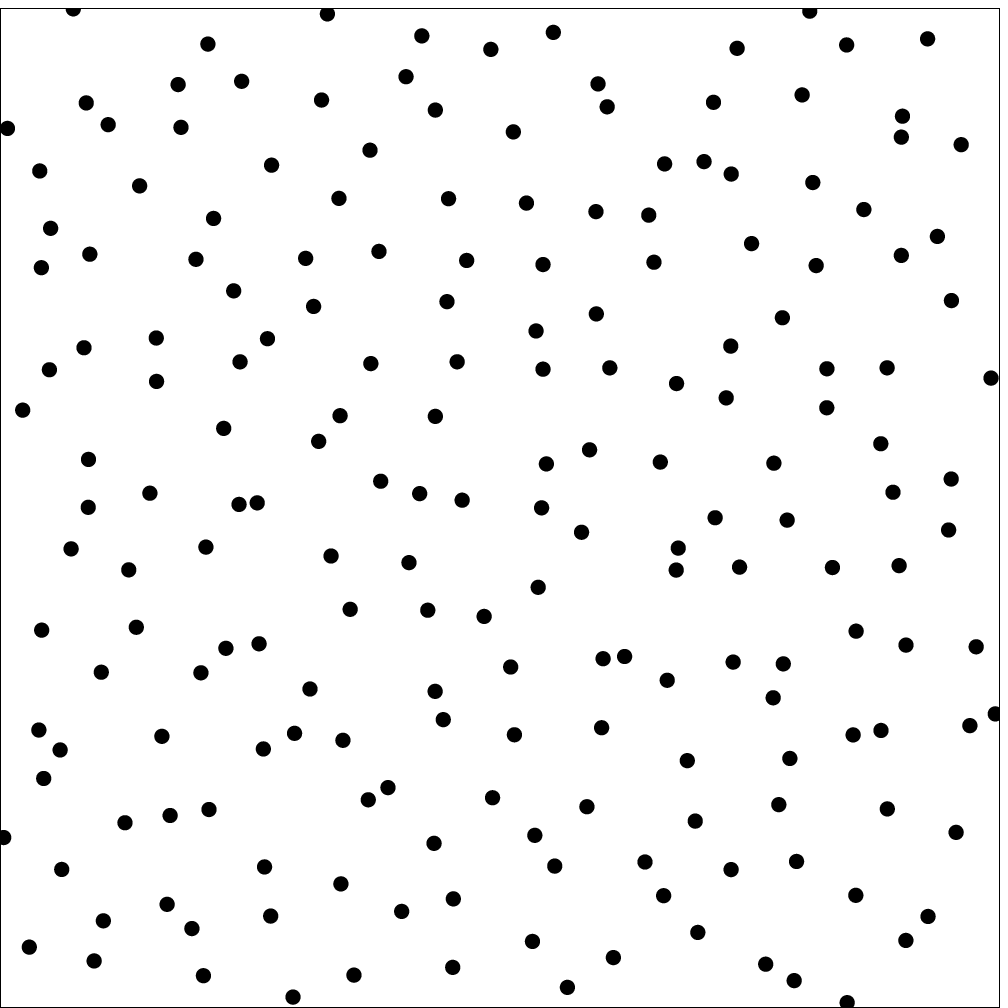}
\caption{\label{fig:gin&gaf} Samples of a translation invariant
  determinantal process(left) and zeros of a Gaussian analytic
   function(right). Determinantal processes exhibit repulsion at all
   distances, and the zeros repel at short distances only. The picture
   in the center shows a Poisson process (picture due to B\'{a}lint Vir\'{a}g)}
\end{figure}

One of the main qualitative properties of zero sets of random analytic
functions is that they have the property of ``repulsion'', also called
``negative correlation'', at short ranges. This terminology is a little
misleading because correlations are never negative! The precise
meaning of negative correlations is that $\rho_2(x,y)<\rho_1(x) \rho_1(y)$ 
for $x,y$ that are sufficiently close. There is another class of point
processes that has this repulsion  property at all distances in a very
strong sense. These point processes were introduced by
Macchi~\cite{mac} and are known as {\bf Determinantal
  (Fermionic) point processes}.  See Figure~\ref{fig:threeprocesses}
for a visual comparison of zeros and eigenvalues with a Poisson process. 

\begin{definition}\label{def:det} A point process $\X$ on $\M$ is said to be a
  {\bf determinantal process} with kernel $\Kdet:\M^2 \tends \C$, if
  it is simple (i.e., there are no coincident points almost surely) and
  its correlation functions w.r.t a measure $\mu$ satisfy: 
    \begin{equation}\label{eq:detdefn}
    \rho_k(x_1,\ldots,x_k)=\det\l( \Kdet(x_i,x_j)\r)_{1\le i,j \le k},
    \end{equation}
for every $k\ge 1$ and $x_1,\ldots ,x_k \in \M$. We shall always assume
that $\Kdet$ is the projection kernel of a closed subspace of $L^2$,
i.e., that $\Kdet(x,y)=\summ_j \phi_j(z)\bar{\phi_j(y)}$ for a sequence of
functions (there may be infinitely many of them) $\{\phi_j \}$
orthonormal in $L^2(\mu)$. The distribution is determined by giving
the kernel $\Kdet$ or the Hilbert space on which it is a
projection. 
\end{definition}

\begin{remark}
  This definition may look artificial at first sight, but the
  motivation comes from quantum mechanics, where the
 probability  densities are given by the absolute square of a complex-valued function called the amplitude. If $\phi_1,\ldots ,\phi_n$
  (assumed to be orthogonal) are 
  single particle wave functions of electrons, the most natural
  $n$-particle wave  function is not $\prodd_k \phi_k(x_k)$ (as would
  have been the case for ``independence'') because firstly it has no
  symmetry in $x_1,\ldots ,x_n$, and secondly it shows no properties such as
  repulsion. Hence the idea is to anti-symmetrize this, to get $\det
  \l(\phi_i(x_j) \r)$ (This is for Fermions. For Bosons, one symmetrizes and gets the {\it permanent} of $\l(\phi_i(x_j) \r)$, but that is another story).  The probability density is given by the 
  absolute square of this, which can be written as $\det\l(
  \Kdet(x_i,x_j)\r)_{1\le i,j \le n}$, with $\Kdet(x,y)=\summ_{j=1}^n
  \phi_j(x) {\bar \phi_j(y)}$. To generalize this notion to infinite
  particle systems, it is necessary to formulate the definition in
  terms of correlation functions as done in Definition~\ref{def:det}. For
  correlation functions to be non-negative in (\ref{eq:detdefn}) a natural
  assumption is to take $\Kdet$ to be a Hermitian non-negative
  definite kernel. Then it turns out that $0\le \Kdet \le 1$ is a
  necessary condition for there to exist a determinantal process with
  kernel $\Kdet$ (see ~\cite{sos1} or ~\cite{hkpv}). And then, any
  such process can be 
  expressed as a mixture of determinantal processes with projection
  kernels (see ~\cite{takshi1},\cite{takshi2},\cite{hkpv}). The rank
  of the projection is the number of points in the process. Observe
  that $\rho_2(x,y)-\rho_1(x)\rho_1(y)=-|\Kdet(x,y)|^2$ is
  negative. More general inequalities like this are a clear
  consequence of the determinantal form of the correlation functions.
\end{remark}
\begin{remark} Another justification for this definition is the huge
  number of instances 
  of determinantal processes that have arisen so far in random matrix
  theory and combinatorics, many of them predating the definition. To
  name a few,
\begin{itemize}
\item Non-intersecting random 
  walks on the line (Karlin and McGregor~\cite{KarMg},\cite{johan}).
\item Eigenvalues of a random 
  unitary matrix chosen according to the Haar
  measure (Dyson~\cite{dys1}). (There are many more random matrix examples).
\item Subset of edges of a finite graph present in a uniformly chosen
  spanning tree (Burton and  Pemantle~\cite{burpe}).
\item  (An encoding of) Young diagrams sampled from the 
  poissonised Plancherel measure of the 
  symmetric group (Borodin, Okounkov and Olshanski~\cite{borokools}).
\end{itemize} 
 The extensive survey of  Soshnikov~\cite{sos1} gives many more
 examples and details.
\end{remark}

While determinantal processes on $\R$ or $\Z$ have been studied extensively
because they arise in random matrix theory and combinatorics, in two
dimensions they seem to have been largely untouched.

To get a determinantal point process is trivial. Just take a
reproducing kernel Hilbert space $\H$ of functions on $\M$, and let
$\Kdet$ be the reproducing kernel of $\H$. Then there exists a
determinantal point process with kernel $\Kdet$ (trivial when
dim($\H$)$<\infty$, otherwise it can be constructed by taking limits of
finite dimensional ones. See ~\cite{sos1} or \cite{hkpv} for
details). However to deserve our attention, the process must have some
attractive features in addition  to its determinantal nature. We
adhere to two guiding principles. 
\begin{itemize}
\item That the process be invariant in distribution (i.e., stationary)
  under a rich group
  of transformations.
\item That the process arise in a natural way probabilistically.
\end{itemize}
Although imprecise, the first principle suggests that we
consider the same three domains $\C,\S^2,\D$. And since we would
ultimately want to relate these to zeros of random analytic functions,
we consider Hilbert spaces of {\it analytic} functions on these spaces.
  We show that the Hilbert spaces that give rise to {\it stationary}
  determinantal processes are
 precisely the well known {\bf Bargmann-Fock} spaces of analytic
 functions. 
 But before going into this, we should point out that these
 determinantal processes were already studied by
 Caillol~\cite{caillol} under the name ``One component plasma on the
 sphere'' (the 
 two-component version was studied by Forrester, Jancovici
 and Madore~\cite{forjanmad}), and Jancovici and
 T{\'e}llez~\cite{jantell} from the 
 point of view of  constructing  Coulomb gases on these spaces. We
 arrived at these processes prompted by a question of B\'{a}lint
 Vir\'{a}g as to what natural determinantal processes can be defined on the
 two dimensional sphere. As we remarked earlier, one can {\it define} a
 determinantal process 
  (which can be regarded as Coulomb gas at a particular temperature
 $\beta=2$) by choosing one's favourite Hilbert space with a
 kernel. What is new 
 here is that we show that these are unique in a certain sense,
 and most importantly, we show in the next chapter how to get these
 determinantal processes as zeros of random analytic functions.

\begin{theorem}\label{thm:invdets} Let $\M$ be one of $\C, \S^2,\D$
  and let $\mu$ be an arbitrary 
radially symmetric Radon measure on $\M$. Let  $\phi_k$ be complex
analytic functions on $\M$ and belong to $L^2(\mu)$. Then the determinantal
process with kernel 
\begin{equation}
\Kdet(z,w) = \summ_{k=1}^N \phi_k(z) {\bar \phi_k(w)}
\end{equation}
is invariant in distribution under the corresponding group of
isometries if and only if it is one of the following.
\begin{itemize}
\item $\M=\C$, $d\mu(z)=e^{-\alp |z|^2}dm(z)$, $N=\infty$,
  $\phi_k(z)=\sqrt{\fr{\alp}{\pi}}z^k$, where $\alp\in (0,\infty)$,
  and the kernel is
\begin{equation}\label{eq:detplanekernel}
\Kdet_{\alp}^{\C}(z,w) = \fr{\alp}{\pi} e^{\alp z{\bar w}}.
\end{equation}
The Hilbert space is the space of analytic functions in
$L^2(\C,\frac{1}{\pi}e^{-\alp |z|^2})$. 
We call this process \underline{\mb{Det-$\C$-$\alp$}}.
\item $\M=\S^2$, $d\mu(z)=\fr{dm(z)}{(1+|z|^2)^{\fr{\alp+1}{2}}}$,
    $N=\alp$, $\phi_k(z)=\sqrt{\fr{\alp}{\pi}}z^k$, where $\alp\in
    \{1,2,3,\ldots \}$, and the kernel is 
\begin{equation}\label{eq:detspherekernel}
\Kdet_{\alp}^{\S^2}(z,w) = \fr{\alp}{\pi} (1+z{\bar
    w})^{\alp-1}.
\end{equation}
The Hilbert space is the space of analytic functions in
$L^2(\S^2,\frac{1}{\pi (1+|z|^2)^{\alp+1}})$. 
We call this process \underline{\mb{Det-$\S^2$-$\alp$}}.
\item $\M=\D$, $d\mu(z)=(1-|z|^2)^{\fr{1}{2}(\alp-1)}dm(z)$,
    $N=\infty$, $\phi_k(z)=\sqrt{\fr{\alp}{\pi}}z^k$, where $\alp\in
    (0,\infty)$, and the kernel is 
\begin{equation}\label{eq:detdiskkernel}
\Kdet_{\alp}^{\D}(z,w) = \fr{\alp}{\pi} \fr{1}{(1-z{\bar w})^{\alp+1}}.
\end{equation}
The Hilbert space is the space of analytic functions in
$L^2(\D,\frac{1}{\pi}(1-|z|^2)^{\alp-1})$.
We call this process \underline{\mb{Det-$\D$-$\alp$}}.
\end{itemize}
\end{theorem}

\begin{remark}
  Note the similarity to the classification of Gaussian analytic
  function with stationary zeros in
  Proposition~\ref{prop:stationarygafs}. Just as there, here too,
  Det-$\C-\alp$ processes are all identical up to scale, whereas
  the Det-$\S^2-\alp$ and Det-$\D-\alp$  are genuine one
  parameter families of point processes.

 There are (at least) two ways of using a positive
  definite kernel in probability theory. One is to use it as the covariance
  kernel of a Gaussian process and another is to use it as the
  kernel of a determinantal process (it has to be a projection kernel
  for the latter). We are not aware of any probabilistic connection
  between the two.

 When we use the reproducing kernels of the Bargmann-Fock spaces as
 covariance kernels of Gaussian processes, we get the Gaussian
 analytic functions introduced 
 in (\ref{eq:planargaf}), (\ref{eq:sphericalgaf}) and
 (\ref{eq:hyperbolicgaf}). When we use
 them as kernels for determinantal processes we  get stationary point
 processes in these domains.  
\end{remark}

\begin{proof}[Proof of Theorem~\ref{thm:invdets}]
Let
\begin{equation*}
\rho^2(z) = \l\{ \begin{array}{cc}
              \frac{1}{\pi} & \mb{ if }\M=\C. \\
              \frac{1}{\pi(1+|z|^2)^2} & \mb{ if }\M=\S^2. \\
              \frac{1}{\pi(1-|z|^2)^2} & \mb{ if }\M=\D. \\
              \end{array} \r.
\end{equation*}
Then $\rho^2(z)dm(z)$ is the unique invariant measure (up to
multiplication by a constant) on $\M$. 

If $\X$ is a determinantal process with kernel $\Kdet$ on $\M$, with
distribution invariant under the corresponding isometry group, then the first
intensity of the process,  $\Kdet(z,z)d\mu(z)$ must be equal to $\alp
\rho^2(z)dm(z)$ for some $\alp>0$.

Express the correlation functions of $\X$ w.r.t the measure
$\rho^2(z)dm(z)$ instead of the Lebesgue measure. Then the kernel becomes 
\begin{equation*}
\frac{\alp \Kdet(z,w)}{\sqrt{\Kdet(z,z)}\sqrt{\Kdet(w,w)}}.
\end{equation*} 
Invariance of the second correlation function implies that
\begin{equation}
  \label{eq:2ptcorr}
  \frac{\alp
  |\Kdet(\phi(z),\phi(w))|^2}{\Kdet(\phi(z),\phi(z))\Kdet(\phi(w),\phi(w))}=\mb{Const}(z,w), 
\end{equation}
for every isometry $\phi$ of $\M$, where by $\mb{Const}(z,w)$ we mean
that it does not depend on $\phi$.

The idea is this. We differentiate equation (\ref{eq:2ptcorr})
w.r.t $\phi$ and equate to zero. The derivatives w.r.t $\phi$ can be
written as derivatives w.r.t $z,w$. That gives us differential
equations for $\Kdet$ that are easy to solve.

Firstly, fix any $(z_0,w_0)$ such that $\Kdet(z_0,w_0)\not=0$. Without
loss of generality take $z_0=0=w_0$. Then there
is a neighbourhood $N$ of $0$ in the Complex plane such that if
$z,w \in N$, and $\phi$ is close to identity, then
$|\Kdet(z,w)-\Kdet(0,0)|<\frac{1}{2}|\Kdet(0,0)|$. Let $S$ be the disk of
radius $\frac{1}{2}|\Kdet(0,0)|$ centered at $\Kdet(0,0)$. Then $S$
cannot intersect both the positive and negative parts of real axis
because $S$ is convex and does not contain $0$. Moreover ${\bar S}$ and
  $S$ intersect the real line at the same points. Therefore by
  removing the positive {\it or} the negative half line, we can define a
 continuous branch of logarithm on $S\cup {\bar S}$. Henceforth ``$\log$''
 will denote this function.

Taking logarithms in Equation(\ref{eq:2ptcorr}) we get
\begin{equation}\label{eq:2ptcorr2}
  \log \Kdet(\phi(z),\phi(w)) +\log \Kdet(\phi(w),\phi(z)) - \log
  \Kdet(\phi(z),\phi(z)) -\log \Kdet(\phi(w),\phi(w)) 
\end{equation}
is equal to const$(z,w)$, not depending on $\phi$.
To differentiate w.r.t $\phi$ we parameterize it with complex numbers
as follows.

\begin{itemize}
\item {Complex plane} Write $\phi(z)=\lam z+\alp$, where $\alp\in C$, $|\lam|=1$. 
\item {Sphere} Write $\phi(z)=\frac{\alp z +\bet}{-{\bar \bet}z+{\bar
      \alp}}$, where $\alp,\bet\in \C$, $|\alp|^2+|\bet|^2=1$.
\item {Disk} Write $\phi(z)=\frac{\alp z +\bet}{{\bar \bet}z+{\bar
      \alp}}$, where $\alp,\bet\in \C$, $|\alp|^2-|\bet|^2=1$. 
\end{itemize}

Let us deal with the planar case first.

{\bf Complex plane} Write $\phi_{\alp}(z)=z+\alp$, where $\alp\in C$. Then
for $\eps$ is small enough, if $|\alp|<\eps$ and $z,w \in N$,
then $\Kdet(\phi_{\alp}(z),\phi_{\alp}(w))\in S$. Apply $\frac{\d^2}{\d
  \alp \d {\bar \alp}}$ to equation (\ref{eq:2ptcorr2}) and evaluate at
$\alp=0$. We get
\begin{eqnarray*}
  0 &=& \frac{\d^2}{\d z \d {\bar w}}\log \Kdet(z,w) + \frac{\d^2}{\d w \d
  {\bar z}}\log \Kdet(w,z) -\frac{\d^2}{\d z \d {\bar z}}\log \Kdet(z,z) -
  \frac{\d^2}{\d w \d {\bar w}}\log \Kdet(w,w) \\
    &=& Q(z,w)+Q(w,z)-Q(z,z)-Q(w,w).
\end{eqnarray*}
Where $Q(z,w)=\frac{\d^2}{\d z \d {\bar w}}\log \Kdet(z,w)$. This is well
defined for $z,w \in N$ and is analytic in $z,{\bar
  w}$. Applying $\frac{\d^2}{\d z \d {\bar w}}$ to
$Q(z,w)+Q(w,z)-Q(z,z)-Q(w,w)$ we deduce that $\frac{\d^2}{\d z \d
  {\bar w}}Q(z,w)=0$. By expanding Q locally as a power series in
$z,{\bar w}$, and observing the symmetry $Q(z,w)={\bar Q(w,z)}$, we conclude 
that $Q(z,w)=h(z)+{\bar h(w)}$ for some $h$ holomorphic on $N$. Then
\begin{equation*}
  \log \Kdet(z,w) = {\bar w}h_1(z)+z{\bar h_1(w)} +C, \hsp{2cm} \forall
  z,w\in N,
\end{equation*}
where $h_1'=h$ and $C$ is a constant. 

Now again consider (\ref{eq:2ptcorr2}) and apply $\frac{\d}{\d
  \alp}$ to it. We get
\begin{eqnarray*}
  0 &=& \frac{\d}{\d z} \log \Kdet(z,w) + \frac{\d}{\d w} \log \Kdet(w,z) -
  \frac{\d}{\d z} \log \Kdet(z,z) - \frac{\d}{\d w} \log \Kdet(w,w) \\
    &=&  {\bar w}h(z)+{\bar h_1(w)}+{\bar z}h(w)+{\bar h_1(z)}-{\bar
  z}h(z)-{\bar h_1(z)}-{\bar w}h(w)-{\bar h_1(w)}\\
    &=& ({\bar z}-{\bar w})(h(w)-h(z)).
\end{eqnarray*}
Since this has to hold for every $z,w\in N$, we must have
$h\equiv 0$. That means that $h_1(z)=az+b$ for some $a,b$. Therefore, 
\begin{eqnarray*}
  \log \Kdet(z,w) &=& \bar{w} (az+b) + z\bar{(aw+b)} +C \\
              &=& (a+\bar{a}) (z+\frac{b}{a+\bar{a}})\bar{(w+\frac{b}{a+{\bar a}})} +C'
\end{eqnarray*}
Making a change of variables we get $\Kdet(z,w)=e^{\alp z{\bar w}}$. Then
$\alp$ has to be positive. This concludes the planar case.

{\bf Sphere} In this case, we again differentiate
Equation(\ref{eq:2ptcorr}) w.r.t $\alp,\bet$ and their conjugates. However
$\phi_{\alp,\bet}$ depends on the parameters and their conjugates as
well, and that makes the equation longer. A simplification is obtained
by noting the following:

Let $g(z,w)$ be analytic in $z,{\bar w}$. Then with $\phi=\phi_{\alp,\bet}$, 
\begin{equation}
  \label{eq:spherediff}
(\frac{\d}{\d \bet} + {\bar w}\frac{\d}{\d \alp})(\frac{\d}{\d {\bar
    \bet}}+z\frac{\d}{\d {\bar \alp}})
    g(\phi(z),\phi(w))
= \frac{(1+z{\bar w})^2}{(-{\bar \bet} z+\bar{\alp})(-\bet {\bar
  w}+\alp)}(\d_1{\bar \d_2}g)(\phi(z),\phi(w)).
\end{equation}
Here $\d_1,{\bar \d_2}$ denote the derivatives w.r.t the first and
second arguments.
 
Apply $(\frac{\d}{\d \bet} + {\bar w}\frac{\d}{\d \alp})(\frac{\d}{\d {\bar
    \bet}}+z\frac{\d}{\d {\bar \alp}})$ to Equation (\ref{eq:2ptcorr2})
  and evaluate at $\alp=1,\bet=0$. Then (\ref{eq:spherediff})
  yields, 
  \begin{equation*}
    0 = Q(z,w)+Q(w,z)-Q(z,z)-Q(w,w),
  \end{equation*}
where, this time $Q(z,w)=(1+z{\bar w})^2  \frac{\d^2}{\d z \d {\bar
    w}}\log \Kdet(z,w)$. But all the considerations that applied to $Q$ in the
    Planar case also apply here and we get $\Kdet(z,w)=(1+z{\bar
    w})^{\alp}$ (perhaps after a change of variables). Now $\alp$ will have to be a positive integer
    (because integrating $\Kdet(z,z)$ over the whole space should give
    $N!$, where $N$ is the total number o points in the process).

{\bf Disk} Analogous to the spherical case, here we observe that for
any $g$ analytic in $z$ and anti-analytic in $w$,
\begin{equation}
  \label{eq:diskdiff} 
(\frac{\d}{\d \bet} - {\bar w}\frac{\d}{\d \alp})(\frac{\d}{\d {\bar
    \bet}} - z\frac{\d}{\d {\bar \alp}})
    g(\phi_{\alp,\bet}(z),\phi_{\alp,\bet}(w))
=\frac{(1-z{\bar w})^2}{({\bar \bet} z+\bar{\alp})(\bet {\bar
  w}+\alp)}(\d_1{\bar \d_2}g)(\phi_{\alp,\bet}(z),\phi_{\alp,\bet}(w)).
\end{equation}


Applying $(\frac{\d}{\d \bet} - {\bar w}\frac{\d}{\d \alp})(\frac{\d}{\d {\bar
    \bet}} - z\frac{\d}{\d {\bar \alp}})$ to (\ref{eq:2ptcorr2})
  and using equation (\ref{eq:diskdiff}) gives us 
\begin{equation*}
    0 = Q(z,w)+Q(w,z)-Q(z,z)-Q(w,w),
  \end{equation*}
where, this time $Q(z,w)=(1-z{\bar w})^2  \frac{\d^2}{\d z \d {\bar
    w}}\log \Kdet(z,w)$. As before this leads us to $\Kdet(z,w)=(1-z{\bar
    w})^{\alp}$ and $\alp$ will have to be positive.

This completes the proof of the theorem.
\end{proof}


\chapter{Random matrix-valued analytic functions}\label{chap:matgaf}
In Chapter~\ref{chap:invzeros}, we saw that by choosing a Gaussian
analytic function $\f$ with a stationary zero set and a (non-random)
homogeneous 
polynomial $\Q$, we could construct a random analytic function with stationary
zeros. There are two complementary questions that arise naturally.
\begin{enumerate}
\item Can one study these random analytic functions in this
  generality without having to appeal to special $\Q$ and $\f$?
\item Are there particular examples of $\Q$ and $\f$ that are somehow special?
\end{enumerate}

The answer to both these questions is yes. Regarding the first 
question, we already saw in Chapter~\ref{chap:invzeros} that the
correlation functions of the zero set can be computed in a general
fashion. We shall use these computations in Chapter~\ref{chap:normality} to
prove asymptotic normality for the zero sets in general.  In the
current chapter and the next two, we answer the second question and
show that zeros of 
random analytic functions sometimes (but far from frequently, let
alone always) turn out to be determinantal point processes.

\section{Determinantal processes that are zeros of RAFs:  Known results}

We remarked earlier that the focus in random matrix theory has been on
Hermitian random matrices. In the preface to his book ``Random matrices'', 
Mehta~\cite{mehta} says ``The theory of non-Hermitian random
matrices, though not applicable to any physical problems, is a
fascinating subject and must be studied for its own sake. In this
direction an impressive step [has been taken by] Ginibre ...'' Ginibre
found the exact distribution of eigenvalues of three (two, strictly
speaking) ensembles of non-Hermitian random matrices. We quote the one
that is relevant to us. 

\begin{theorem}[Ginibre(1965)~\cite{gin} ]\label{thm:ginibre}
Let $A$ be an $n\times n$ matrix with i.i.d.~standard complex Gaussian
entries. Then the eigenvalues of $A$ have density
\begin{equation}
  \label{eq:ginibredensity}
  \rho_n(z_1,\ldots,z_n) = \frac{1}{\pi^n\prodd_{k=1}^{n-1}k! }
  e^{-\summ_{k=1}^n |z_k|^2} \prod_{i<j} |z_i-z_j|^2.  
\end{equation}
Equivalently, one may say that the eigenvalues of $A$ form a
determinantal point process with kernel
\begin{equation}\label{eq:gindetform}
  \Kdet_n(z,w) = \summ_{k=0}^{n-1}
  \frac{(z\bar{w})^k}{k!},
\end{equation}
w.r.t the reference measure $d\mu(z)=\frac{1}{\pi}e^{-|z|^2}$. 
The corresponding Hilbert space $\H=\mb{span}\{1,z,\ldots ,z^{n-1}
\}\subset L^2(\C,\frac{e^{-|z|^2}}{\pi}dm(z))$.
\end{theorem}
Despite the enthusiastic response, as shown by Mehta's quote above,
there do not seem to be any significant {\it exact} results beyond
Ginibre's. The following beautiful result of Peres and
Vir\'{a}g~\cite{pervir} seems to be the next such.
\begin{theorem}[Peres and Vir\'ag(2003)~\cite{pervir} ]\label{thm:pervir}
 Let $\f$ be the random power series whose coefficients are
 i.i.d. standard complex Gaussians (this is the case $\M=\D, L=1$ in
  (\ref{eq:hyperbolicgaf}). Then the zeros of $\f$ form a determinantal point process
 on the unit disk $\D$ with the kernel (the Bergman kernel of the unit disk)
 \begin{equation*}
   \Kdet(z,w) = \frac{1}{\pi(1-z\bar{w})^2},
 \end{equation*}
w.r.t the reference measure $d\mu(z)=\frac{1}{\pi}dm(z)$ on $\D$.
 The corresponding Hilbert space $\H=\mb{span}\{1,z, z^2 \ldots , 
\}\subset L^2(\D,\frac{dm(z)}{\pi})$ is the space of all analytic
functions in $L^2(\D)$.
\end{theorem}

\begin{remark} Observing that Theorem~\ref{thm:pervir} identifies the
 distribution of zeros of the Gaussian analytic function (with $L=1$) defined in
 (\ref{eq:hyperbolicgaf}) as being Det-$\D-1$ (recall the definition of Det-$\D-1$ from (\ref{eq:detdiskkernel})), one is tempted to guess
 that the Gaussian analytic functions defined in (\ref{eq:planargaf}) 
  (\ref{eq:sphericalgaf}) and (\ref{eq:hyperbolicgaf}) might have zeros
 distributed like the Det-$\C-L$, Det-$\S^2-L$ and Det-$\D-L$ (defined
 in
 (\ref{eq:detplanekernel}), (\ref{eq:detspherekernel}) and
 \ref{eq:detdiskkernel}),  respectively). However these canonical 
 Gaussian analytic functions do not have determinantal zero
 sets. Indeed, it was 
 observed by Peres and Vir\'{a}g in their paper that these zero sets
 do {\it not} have negative correlations at large distances and hence,
 cannot be determinantal). Therefore this 
 beautiful might-have-been story is completely false! Nevertheless, to
 quote Einstein, ``The Lord is subtle, but not malicious''. In the next
 section we shall see how a completely different but equally
 compelling picture might well be true.
\end{remark}

\section{Determinantal processes that are zeros of RAFs: An analogy}
First let us list all the {\it Gaussian} analytic functions whose zero
sets we know to be determinantal. This includes
Theorem~\ref{thm:pervir} and two trivial cases (a one-point point
process is always determinantal!). 

$\bullet$ $z-a=0$ : One zero, with standard complex Gaussian
distribution on $\C$. $\H=\mb{span}\{1\}$ in
$L^2(\C,\frac{e^{-|z|^2}}{\pi}dm(z))$. 

$\bullet$ $za-b=0$ : One zero, distributed uniformly on $\S^2$ upon
stereographic projection from the plane. $\H=\mb{span}\{1\}$ in
$L^2(\S^2=\C \cup\{\infty\},\frac{1}{\pi(1+|z|^2)^2}dm(z))$ 

$\bullet$ $a_0+za_1+z^2a_2+\dots =0$ : Peres and Vir\'{a}g~\cite{pervir}:
  Infinitely many zeroes in the disk. A determinantal point
  process with kernel $\pi^{-1}
  (1-|z|^2)^{-2}$ w.r.t Lebesgue measure on the unit disk. Equivalently, $\H=\mb{span}\{1,z,z^2,\ldots \}$ in  $L^2(\D,\frac{1}{\pi}dm(z))$.  

Note that Ginibre's result (Theorem~\ref{thm:ginibre}) can be seen as
regarding the zeros of the random analytic function $zI-A$, which can
be thought of as a matrix version of the first of the above
examples. This suggests that we consider the matrix versions of the
other two, i.e., we look at 
\begin{itemize}
\item $\det (zA-B)=0$, where $A,B$ are $n\times n$ independent matrices with
  i.i.d. standard complex Gaussian entries.

\item $\det\l(A_0+zA_1+z^2A_2+\ldots \r)=0$, where $A_k$ are
  independent $n\times n$ matrices with each one having i.i.d. complex
  Gaussian entries.
\end{itemize}

The analogy strongly suggests that the solutions to these equations
should give us the determinantal point processes corresponding to the
Bargmann-Fock spaces on the sphere and the unit disk (but only for
integer values of the parameter, since the size of the matrix, namely
$n$, is discrete). Before going into the statements and proofs, we
 make some big-picture remarks and connect these objects to the random
 analytic functions studied in Chapter~\ref{chap:invzeros}.

\begin{remark} Note that here we are looking at the set of $z$ for which a random
matrix-valued  analytic function ($zA-B$ or $A_0+zA_1+z^2A_2+\ldots $)
becomes singular. This concept is an obvious generalization of both
random matrices (which correspond to the case when the analytic
function is linear) and Gaussian analytic functions (which correspond
to the case when the matrices have size $1\times 1$). In spite of this
natural appeal, the concept of a random matrix-valued analytic
function does not seem to have been considered in the literature. 
 
One possible reason could be that the focus in random matrix theory
has been almost entirely on eigenvalues in one dimension (real line or
the circle) for physical reasons as well as the strong mathematical
connections with orthogonal polynomials, representation theory
etc. Moreover the eigenvalues have a physical meaning in quantum
mechanics. Note the difficulty of forcing the zeros to lie on the real line, 
except by considering eigenvalues of a Hermitian matrix. Nevertheless,
the idea of matrixifying seems to be useful, not only as suggested
above with Gaussian matrix coefficients, but also polynomials with
coefficients that are Haar-distributed unitary matrices.
\end{remark}

\section{Matrix-valued GAFs and polygafs}
Now we want to point out the connection with homogeneous polynomials
applied to i.i.d. copies of Gaussian analytic functions (polygafs, that is). 

Consider $\det(zA-B)$. This is the same as applying the homogeneous
polynomial $\Q=$''$\det$'' in $n^2$ variables, to $n^2$ i.i.d. copies of
the Gaussian analytic function $\f(z)=az-b$ (which is the case
$\M=\S^2,L=1$ in (\ref{eq:sphericalgaf})).

Similarly $\det\l(A_0+zA_1+z^2A_2+\ldots \r)$ is the  homogeneous
polynomial $\Q=''\det''$ in $n^2$ variables, applied to $n^2$ i.i.d. copies
of the  Gaussian analytic function $\f(z)=a_0+a_1z+a_2 z^2+\ldots $
(which is the case  $\M=\D,L=1$ in (\ref{eq:hyperbolicgaf})). 

In other words, we have already shown in 
Proposition~\ref{prop:stationaryrafs} that the zero sets of these 
RAFs are stationary in $\M$. In the next two chapters we investigate the
distributions in greater depth. We shall show that in the first case
($\M=\S^2$) we do get determinantal processes, whereas in the second,
we show partial results in this direction. The precise statements of
the conjectures  are as follows:

\begin{conjecture}\label{conj:sphere} Let $A,B$ be i.i.d $n\times n$ matrices with
i.i.d. standard complex Gaussian entries. The zeros of $\det(zA-B)$ form a
determinantal point process with kernel
\begin{equation*}
  \Kdet(z,w) = \frac{n}{\pi}
 \frac{(1+z\bar{w})^{n-1}}{(1+|z|^2)^{\frac{n+1}{2}}(1+|w|^2)^{\frac{n+1}{2}}},    
\end{equation*}
w.r.t the Lebesgue measure on $\C$. Equivalently, $\Kdet$ is the
projection kernel on the subspace of {\it analytic} functions in
$L^2 (\S^2, \frac{n}{\pi (1+|z|^2)^{n+1}}dm(z))$.
\end{conjecture}

\begin{conjecture}\label{conj:disk} Let $A_k$ be i.i.d. $n\times n$ matrices with
i.i.d. standard complex Gaussian entries. The zeros of $\det\l(
A_0+zA_1+z^2A_2+\ldots \r)$ form a determinantal point process on $\D$
with kernel
\begin{equation*}
  \Kdet(z,w) = \frac{n}{\pi}
 \frac{(1-|z|^2)^{\frac{n-1}{2}}(1-|w|^2)^{\frac{n-1}{2}}}{(1-z\bar{w})^{n+1}},    
\end{equation*}
w.r.t the Lebesgue measure on $\D$. Equivalently, $\Kdet$ is the
projection kernel on the subspace of {\it analytic} functions in
$L^2 (\D, \frac{n}{\pi}(1-|z|^2)^{n-1}dm(z))$.
\end{conjecture}

\begin{remark} The Det-$\C-1$ process are obtained 
  from Ginibre's theorem~\ref{thm:ginibre} by letting $n\tends
  \infty$. So from the point of view of determinantal processes, our
  problems can be stated as finding a probabilistic meaning to
  Det-$\S^2-\alp$ and Det-$\D-\alp$. 
\end{remark}


\chapter{Matrix analytic functions on the sphere}\label{chap:sphere}
In this section we prove Conjecture~\ref{conj:sphere} stated at the end of
Chapter~\ref{chap:matgaf}, i.e., we show that the processes
\mb{Det-$\S^2$-$\alp$}  arise as the singular points of the matrix GAF
$zA-B$ or equivalently, zeros of the polygaf $\det(zA-B)$. Recall  
that for $\M=\S^2$, $\alp$ is a positive integer (the number of points in the
process). We shall denote it by $n$ in this section. 

\begin{theorem}\label{thm:dets2} Let $A,B$ be independent $n\times n$ random
  matrices with i.i.d. standard complex Gaussian entries. Then the set
  $\X$ of zeros of $\det(zA-B)=0$, i.e., the eigenvalues of $A^{-1}B$,
  has the distribution \mb{Det-$\S^2$-$n$}.  
\end{theorem}
We need the following lemma.
\begin{lemma}\label{lem:invdist} Let $\X$ be a point process on $\C$
  with $n$ points almost surely. Assume that the $n$-point correlation
  function (equivalently the density) of $\X$ has the form
\[ p(z_1,\ldots ,z_n) = \Mid \Del(z_1,\ldots ,z_n) \Mid^2 V(|z_1|^2,\ldots ,|z_n|^2). \]
Here $\Del(z_1,\ldots ,z_n)$ denotes the Vandermonde factor $\prodd_{i<j}(z_j-z_i)$. 

Suppose also that $\X$ has a distribution invariant under automorphisms of the sphere $\S^2$, i.e.,  under the transformations $\phi_{\alp,\bet}(z)=\frac{\alp z +\bet}{-\bar{\bet}z+\bar{\alp}}$, for any $\alp,\bet$ satisfying $|\alp|^2+|\bet|^2=1$.  
Then 
\begin{equation}\label{eq:specialformofdensity}
V(|z_1|^2,\ldots ,|z_n|^2)=\mb{Const.} \prodd_{k=1}^n
\frac{1}{(1+|z_k|^2)^{n+1}}.  
\end{equation} 
\end{lemma}

\begin{proof}[Proof of Lemma~\ref{lem:invdist}] The claim is that the probability density of the $n$
  points of $\X$ (in exchangeable random order) is 
  \begin{equation*}
    q(z_1,\ldots,z_n) :=\mb{Const.} \Mid \Del(z_1,\ldots ,z_n) \Mid^2
\prodd_{k=1}^n \frac{1}{(1+|z_k|^2)^{n+1}}
  \end{equation*}
First let us check that the density $q$ is invariant under the
isometries of $\S^2$. For this let ${\alp,\bet}(z)=\frac{\alp z
  +\bet}{-\bar{\bet}z+\bar{\alp}}$, with $\alp,\bet$ satisfying
$|\alp|^2+|\bet|^2=1$.  Then, 
\begin{equation}\label{eq:derphi}
\phi^{\prime}(z)= \frac{1}{(-{\bar \bet}z +{\bar \alp})^2}.
\end{equation}
\begin{equation}\label{eq:metphi}
1+\Mid \phi(z) \Mid^2 = \frac{1+\Mid z \Mid^2}{\Mid -{\bar \bet}z +{\bar \alp}\Mid^2}.
\end{equation}
\begin{equation}\label{eq:diffphi}
\phi(z) - \phi(w) = \frac{z-w}{(-{\bar \bet}z +{\bar \alp})(-{\bar \bet}w +{\bar \alp})}.
\end{equation}
From (\ref{eq:derphi}),(\ref{eq:metphi}) and (\ref{eq:diffphi}), it
follows that
\begin{equation}\label{eq:invforq}
 q\l(\phi(z_1),\ldots ,\phi(z_n)\r)\prodd_{k=1}^n |\phi^{\prime}(z_k)|^2 =
 q\l(z_1,\ldots, z_n\r), 
\end{equation}
which shows the invariance of $q$.  

Invariance of $\X$ means that $\forall \alp,\bet$ with
$|\alp|^2+|\bet|^2=1$, and for every $z_1,\ldots ,z_n$, we have 
\begin{equation}\label{eq:invforp}
 p\l(\phi(z_1),\ldots ,\phi(z_n)\r)\prodd_{k=1}^n |\phi^{\prime}(z_k)|^2 =
 p\l(z_1,\ldots, z_n\r). 
\end{equation}
Set $W(z_1,\ldots ,z_n)=\frac{p(z_1,\ldots
  ,z_n)}{q(z_1,\ldots,z_n)}$. Then, from (\ref{eq:invforp}) and
  (\ref{eq:invforq}), we get
\begin{itemize}
\item $W\l(z_1,\ldots ,z_n\r)$ is a function of $|z_k|^2$, $1\le k \le n$, only.
\item $W\l(\phi(z_1),\ldots ,\phi(z_n)\r)=W\l(z_1,\ldots ,z_n \r)$ for
  every $z_1,\ldots ,z_n$. 
\end{itemize}   
We claim that these two statements imply that $W$ is a constant. To
see this fix $z_k=r_k e^{i\theta_k}$, $1\le k\le n$, such that $r_1
<r_k$ for $k\ge 2$. Let $\alp=\frac{1}{\sqrt{1+r_1^2}},
\bet=-\frac{z_1}{\sqrt{1+r_1^2}}$. Then $|\alp|^2+|\bet|^2=1$ and so
$\phi_{\alp,\bet}$ is an isometry of $\S^2$. From the above stated
properties of $W$, we deduce,
\begin{eqnarray*}
  W(z_1,\ldots,z_n) &=& W\l(\phi(z_1),\ldots,\phi(z_n)\r) \\
                    &=& W\l(0,\frac{z_2-z_1}{1+z_2{\bar z_1}},\ldots
                    ,\frac{z_n-z_1}{1+z_n{\bar z_1}}\r) \\
                    &=& W\l(0, \given
                    \frac{r_2e^{i\theta_2}-z_1}{1+r_2e^{i\theta_2}{\bar
                    z_1}} \given, \ldots
                    , \given \frac{r_ne^{i\theta_n}-z_1}{1+r_ne^{i\theta_n}{\bar
                    z_1}} \given \r). 
\end{eqnarray*}
Take $z_1=1$ and $1<r_k<1+\eps$. Then as $\theta_k$, $2\le k\le n$ vary
independently over $[0,2\pi]$, the quantities $\Mid
                    \frac{r_ke^{i\theta_k}-z_1}{1+r_ke^{i\theta_k}{\bar
                    z_1}} \Mid$ vary over the intervals
                $\l[\frac{r_k-1}{r_k+1},\frac{r_k+1}{r_k-1}\r]$. By
                our choice of $r_k$s, this means that
                \begin{equation*}
                  W(0,t_2,\ldots, t_n) = \mb{Constant} \hsp{2cm}
                  \forall t_k\in \l[\eps,\frac{1}{\eps}\r].
                \end{equation*}
$\eps$ is arbitrary, hence $W(0,t_2,\ldots,t_n)$ is
constant. Therefore $W(z_1,\ldots ,z_n)$ is constant.

This shows that $p(z_1,\ldots, z_n)=\mb{Const.} q(z_1,\ldots ,z_n)$. 
\end{proof}

\begin{proof}[Proof of Theorem~\ref{thm:dets2}] Firstly we observe that the
  distribution of $\X$ is invariant under conformal  automorphisms of
  $\S^2$. This is a direct consequence of
  Proposition~\ref{prop:stationaryrafs}, with $\Q=$''$\det$'' and
  $\f(z)=az-b$. Still we give another simple direct proof. Let $\alp,\bet$ be such that
  $|\alp|^2+|\bet|^2=1$. Set 
\[ C=\alp A + {\bar \bet}B \hsp{1cm} \mb{and} \hsp{1cm} D= -\bet A + {\bar \alp} B. \]
Then $C$ and $D$ are i.i.d. matrices with  i.i.d. standard complex
Gaussian entries. Therefore the eigenvalue set of $C^{-1}D$ is also
distributed as $\X$.  

Now $\X$ is the set of solutions to the equation
\[ \det \l( zA-B\r) = 0,\] say $\X=\{z_1,\ldots ,z_n\}$.
 By our observation $\X$ also has the same distribution as the set of
 solutions to \[ 0=\det \l( zC-D\r) = \det\l((z\alp+\bet)A - (-{\bar
 \bet}z +{\bar \alp})B \r), \] which is precisely
 $\l\{\phi_{\alp,\bet}(z_k)\r\}$. This proves the invariance. 

Now by Lemma~\ref{lem:invdist}, it suffices to show that the density
of points in $\X$ is of the form given in (\ref{eq:specialformofdensity}). 

We use the following well known matrix decomposition.

\noindent{\bf Schur decomposition:} Any diagonalizable matrix $M\in GL(n,\C)$
can be written as 
\begin{equation}\label{eq:schurdecomp}
  M = U(Z+T)U^*,
\end{equation}
where $U$ is unitary, $T$ is strictly upper triangular and $Z$ is
diagonal. Moreover the decomposition is almost unique, in the
following sense:

$M=V(W+S)V^*$ in addition to (\ref{eq:schurdecomp}), with $V,S,W$ being respectively  unitary, strictly upper triangular, and diagonal, if and only if
 the element of $W$ are a permutation of the elements of $Z$, and if
 this permutation is identity, then 
 $V=U\Theta$ and $\Theta S \Theta^* =T$ for some $\Theta$ that 
 is both diagonal and unitary, i.e.,
 $\Theta=$ Diag$(e^{i\theta_1},\ldots ,e^{i\theta_n})$ .

Corresponding to this matrix decomposition (\ref{eq:schurdecomp}),
 Ginibre~\cite{gin} proved the measure decomposition

\noindent{\bf Ginibre's measure decomposition:} If $M$ is decomposed as in
(\ref{eq:schurdecomp}), with the elements of 
$Z$ in a uniformly randomly chosen order, then 
\begin{equation}
  \label{eq:ginibredecomp}
  \prod_{i,j} dm(M_{ij}) = \l(\prod_{i<j} |z_i-z_j|^2 \prod_k
  dm(z_k) \r)\l(\prod_{i<j} dm(T_{ij}) \r) d\nu(U)
\end{equation}
where $\nu$ is a finite measure 
on the unitary group $U(n)$ such that $d\nu(U\Theta)=d\nu(U)$ for every
diagonal unitary $\Theta$. 

Conditional on $A$, the matrix $M:=A^{-1}B$ has the density
\begin{equation*}
e^{-\mb{tr}(M^*A^*AM)}|A|^{2n}
\end{equation*} 
w.r.t. the Lebesgue measure on $GL(n,\C)\subset \C^{n^2}$. 
From the measure decomposition  (\ref{eq:ginibredecomp}) we get the
density of $Z$, $T$, $U$, $A$ to be  
\begin{equation*} 
\l(\prod_{i<j} |z_i-z_j|^2 \prod_{k=1}^n
dm(z_k)\r) e^{-\mb{tr}(A^*A(I+MM^*))} \Mid A 
\Mid^{2n}  d \nu(U) \prod_{i<j} dm(T_{ij}) \prod_{i,j}dm(A_{ij}) . 
\end{equation*}
(We have omitted constants entirely)
Thus the density of $Z$ is obtained by integrating over $T,U,A$. Now write 
$Z=\Theta R$ where $\Theta$ and $R$ are diagonal matrices with  the
polar and radial parts of $z_k$, respectively. Then 
\[ MM^*= U\Theta(R+\Theta^*T)(R+\Theta^*T)^*\Theta^*U^*. \] 
As stated earlier, $d \nu(U\Theta)=d \nu(U)$. The elements of
$\Theta^*T$ are the same as elements of $T$, but multiplied by complex
numbers of absolute value $1$. Hence, $\Theta^*T$ has the same
``distribution'' as $T$. Thus replacing $U$ by $\Theta^*U$ and $T$ by
$\Theta^*T$ we see that the density of $Z$  
is of the form $ \prod_{i<j} |z_i-z_j|^2 V(R)$. This is the
form of the density required to apply Lemma~\ref{lem:invdist}. Thus we
conclude that the eigenvalue density is 
\begin{equation}\label{eq:densityofX}
 \mb{Const.}\prodd_{i<j} |z_i-z_j|^2  \prod_{k=1}^n \frac{1}{
 (1+|z_k|^2)^{n+1}}.
\end{equation} 

To compute the constant, note that 
\begin{equation*}
\l\{ \sqrt{\frac{n}{\pi} \binom{n-1}{k}}\frac{z^k}{(1+|z|^2)^{\frac{n+1}{2}}} \r\}_{0\le k \le n-1}
\end{equation*}
is an orthonormal set. Projection on the Hilbert space generated by
 these functions  gives a determinantal process whose kernel
is as given in  the definition of $Det-\S^2-n$. Writing out the density shows
that this is the same as the eigenvalue density that we have
 determined. Hence the constants must match.
\end{proof}


\chapter{Matrix analytic functions on the disk}\label{chap:disk}

Conjecture~\ref{conj:disk} at the end of Chapter~\ref{chap:matgaf} asserted that the
singular points of the matrix GAF $A_0+zA_1+z^2 A_2+\ldots$ or
equivalently, the zeros of the polygaf $\det(A_0+zA_1+z^2
A_2+\ldots)$ are distributed as
$Det-\D-n$. In other words they form a determinantal point process with kernel
\begin{equation*}
  \Kdet_n(z,w) = \fr{n}{\pi} \fr{1}{(1-z{\bar w})^{n+1}}
\end{equation*}
w.r.t the reference measure
$d\mu(z)=(1-|z|^2)^{\fr{1}{2}(n-1)}dm(z)$ on $\D$.
Here $A_k$ are i.i.d. $n\times n$ matrices with
i.i.d. standard complex Gaussian entries.

As already emphasized, Proposition~\ref{prop:stationaryrafs} applies
to show that the zeros of the polygaf $A_0+zA_1+z^2 A_2+\ldots$  are
stationary on the unit disk. In this chapter we shall prove that the
first and second correlation functions agree with those  of the
determinantal process with kernel $\Kdet_n$.

\noindent {\bf First intensity:} In the notation of
Chapter~\ref{chap:invzeros} we have $\Q=${\it 
  determinant}, a  homogeneous polynomial in $n^2$ variables
  $\zet_{ij},i,j\le n$ and $\f_{ij}(z)$ i.i.d. copies of the power
  series $\summ_{n=0}^{\infty} a_n z^n$.

  From (\ref{eq:firstintensityforpolygafs}), the intensity of zeros is
$\frac{n}{4\pi} \lap \log (1-|z|^2)^{-1}$.   By an elementary
computation this comes out to be $\frac{n}{\pi(1-|z|^2)^2}$. Since
this is the same as $\Kdet_n(z,z)$, it follows that the polygaf under
consideration has the same intensity of zeros as the determinantal
process with kernel $\Kdet$ (we omit $n$ in the subscript often).  

\noindent {\bf 2-point correlations:} We shall prove that
\begin{equation}\label{eq:equalityof2ptcorr}
  \rho_2(z,w)-\rho_1(z)\rho_1(w) = - |\Kdet(z,w)|^2,
\end{equation}
which shows that the 2-point correlations for the zeros of the polygaf
agree with those of the determinantal process with kernel
$\Kdet$. (Henceforth correlations are expressed w.r.t. the Lebesgue measure).

\noindent {\it We prove (\ref{eq:equalityof2ptcorr}) by using
Theorem~\ref{thm:dets2} and a change of variables!}

Let  $\g_{ij}(z)=a_{ij}z+b_{ij}$, where $a_{ij},b_{ij}$ are all
i.i.d. $\C N(0,1)$. Then $K(z,w)=1+z\bar{w}$. Consider the polygaf
$\G(z)=\det(\g_{ij}(z))$. This is precisely the polygaf considered in
Chapter~\ref{chap:sphere}.  Thus we know from Theorem~\ref{thm:dets2}
(which we proved by certain matrix decompositions, not at all by using the formulas for correlation functions of polygafs) that the zeros of $\G$  are determinantal with kernel
\begin{equation*}
\frac{n}{\pi}\frac{(1+z\bar{w})^{n-1}}{(1+|z|^2)^{\frac{n+1}{2}}(1+|w|^2)^{\frac{n+1}{2}}}.  
\end{equation*}      

  But the the general formula (\ref{eq:2ptcorrelationformula})
  for two-point correlations of zeros of polygafs applies to $\G$ also
  and hence it must be the case that
\begin{equation}
  \label{eq:2ptspherical}
  \frac{16}{(2\pi)^2} \dz \dzbar \dw \dwbar \summ_{p=0}^{\infty}
  |\tilde{C}_p|^2
  \frac{(1+z\bar{w})^p(1+w\bar{z})^p}{(1+z\bar{z})^p(1+w\bar{w})^p} =
  - \frac{n^2}{\pi^2}\frac{(1+z\bar{w})^{n-1}(1+w\bar{z})^{n-1}}{(1+z\bar{z})^{n+1}(1+w\bar{w})^{n+1}},
\end{equation}
where $\tilde{C_p}$ depend only on $\Q$, not on the GAFs that we feed
in.

From (\ref{eq:2ptspherical}) we get
\def\dz{{\partial \over \partial z}}
\def\dzbar{{\partial \over \partial {\bar z}}}
\begin{equation}\label{eq:4variables}
  \frac{16}{(2\pi)^2} \frac{\partial}{\partial x_1}
  \frac{\partial}{\partial {\bar y_1}}\frac{\partial}{\partial
  x_2}\frac{\partial}{\partial {\bar y_2}} \summ_{p=0}^{\infty} 
  |\tilde{C}_p|^2
  \frac{(1+x_1\bar{y_2})^p(1+x_2\bar{y_1})^p}{(1+x_1\bar{y_1})^p(1+x_2\bar{y_2})^p} =
  - \frac{n^2}{\pi^2}\frac{(1+x_1\bar{y_2})^{n-1}(1+x_2\bar{y_1})^{n-1}}{(1+x_1\bar{x_2})^{n+1}(1+y_1\bar{y_2})^{n+1}}.
\end{equation}
This is because, both sides of (\ref{eq:4variables}) are analytic in
$x_1,x_2$ and anti-analytic in $y_1,y_2$ and moreover,
(\ref{eq:2ptspherical}) says that the two are equal on the diagonal
$\{x_1=y_1,x_2=y_2\}$. Thus, by a standard (and elementary) fact that can be found in any introductory book on several variable complex analysis, see for example Rudin~\cite{rud}, the two sides must be equal for all
$x_1,x_2,y_1,y_2$. 

Now make the substitution $x_1=iz,x_2=iw,y_1=i\bar{w},y_2=i\bar{z}$ to
get
\begin{equation}
  \label{eq:2pthyperbolic}
  \frac{16}{(2\pi)^2} \dz \dzbar \dw \dwbar \summ_{p=0}^{\infty}
  |\tilde{C}_p|^2
  \frac{(1-z\bar{z})^p(1-w\bar{w})^p}{(1-z\bar{w})^p(1-w\bar{z})^p} =
  - \frac{n^2}{\pi^2}\frac{(1-z\bar{z})^{n-1}(1-w\bar{w})^{n-1}}{(1-z\bar{w})^{n+1}(1-w\bar{z})^{n+1}}.
\end{equation}
But again using (\ref{eq:2ptcorrelationformula}) the left hand side is
precisely what we get for $\rho_2(z,w)-\rho_1(z)\rho_1(w)$, when $\Q$ is the
determinant of $n^2$ variables and  $\f_{ij}$ are i.i.d. copies of
$\f(z)=\summ_{n=0}^{\infty} a_n z^n$, with $a_n$ being i.i.d. standard  
complex Gaussians. And the right side of (\ref{eq:2pthyperbolic}) is
$\Kdet(z,w)\Kdet(w,z)$. 

We already know that $\rho_1(z)=\Kdet(z,z)$. Therefore, the
 two point correlation $\rho_2(z,w)$ is 
\begin{equation*}
  \det\l[\begin{array}{cc} \Kdet(z,z) & \Kdet(z,w) \\
                              \Kdet(w,z) & \Kdet(w,w) \end{array} \r].
\end{equation*}


\chapter{Asymptotic Normality}\label{chap:normality}

\section{Background:  Results for Gaussian analytic functions}
 Sodin and Tsirelson~\cite{ST1} proved
asymptotic normality for  smooth ($C_c^2$) statistics applied to the
zeros of the three canonical models of Gaussian analytic functions in 
(\ref{eq:planargaf}), (\ref{eq:sphericalgaf}) and
(\ref{eq:hyperbolicgaf}), as the density parameter $L\tends
\infty$. More precisely, they showed that for any real valued $\phi
\in C_c^2(\M)$, if   
\begin{equation*}
  \ZZ_L(\phi) = \summ_{z\in \f_L^{-1}\{0\}} \phi(z),
\end{equation*}
then,
\begin{equation*}
  \frac{\ZZ_L(\phi)-\E\l[\ZZ_L(\phi)\r]}{\sqrt{\Var(\ZZ_L(\phi))}} \tends N(0,1),
\end{equation*}
and also that
\begin{equation*}
  \Var(\ZZ_L(\phi)) = \frac{\kappa}{L} \|\lap^* \phi \|_{L^2(m^*)}^2 +
  o(L^{-1}), 
\end{equation*}
for a constant $\kappa$ that is described explicitly and the same for
all the three geometries. (Note that as 
$L\tends \infty$, the variance goes to zero!)  Here $m^*$ is the
invariant measure on $\M$ and 
$\lap^*$ is the invariant Laplacian.  In other words 
\begin{equation}
  \label{eq:invmeasure}
  dm^*(z) = \l\{ \begin{array}{ll}
                 dm(z) & \mb{ if } \M=\C \\
                 \frac{dm(z)}{(1+|z|^2)^2} & \mb{ if }\M=\S^2 \\
                 \frac{dm(z)}{(1-|z|^2)^2} & \mb{ if }\M=\D
                 \end{array} \r.
\end{equation}
and
\begin{equation}
  \label{eq:invlaplacian}
  \lap^* = \l\{ \begin{array}{ll}
                 \lap & \mb{ if } \M=\C. \\
                 (1+|z|^2)^2 \lap & \mb{ if }\M=\S^2. \\
                 (1-|z|^2)^2 \lap & \mb{ if }\M=\D. \end{array} \r.
\end{equation}

\section{Our results: for polygafs}
In this article we modify the method of Sodin and Tsirelson to obtain
central limit theorems for smooth statistics of zeros of polygafs.
 One point of this exercise is to demonstrate that the random analytic
functions can be studied at the level of generality introduced in
Chapter~\ref{chap:invzeros} (Recall that polygafs include matrix GAFs as 
very special cases).

\begin{theorem}\label{thm:normality}
Let $\F_L(z)=\Q(\f_1^{(L)}(z),\ldots ,\f_p^{(L)}(z))$ where,
\begin{itemize}
\item  $\Q$ is a fixed non-random homogeneous polynomial in $p$ complex
  variables,   

\item $\f_i^{(L)}$ are i.i.d. GAFs in (\ref{eq:planargaf})
or (\ref{eq:sphericalgaf}) or (\ref{eq:hyperbolicgaf}),

\item $\phi:\M\tends \R$ is a $C_c^2$ function on $\M$.
\end{itemize}
Set
\begin{equation*}
  \ZZ_L(\phi) = \summ_{z\in \F_L^{-1}\{0\}} \phi(z).
\end{equation*}
Then 
\begin{equation*}
  \frac{\ZZ_L(\phi)-\E\l[\ZZ_L(\phi)\r]}{\sqrt{\Var(\ZZ_L(\phi))}}
  \tends N(0,1). 
\end{equation*}
Moreover
\begin{equation}\label{eq:varasymptotics}
  \Var(\ZZ_L(\phi)) = \frac{\kappa(\Q)}{L} \|\lap^* \phi \|_{L^2(m^*)}^2 +
  o(L^{-1}), 
\end{equation}
for a constant $\kappa(\Q)$ that is described explicitly and is the same for
all the three geometries.

\end{theorem}

Let $\hat{\ZZ}_L(\phi)=\ZZ_L(\phi)-\E\l[\ZZ_L(\phi)\r]$. Let $\chi$
denote a standard (real) normal random variable. The idea of the proof is to 
 show that 
\begin{equation}
  \label{eq:convofmoments}
  \E\l[\hat{\ZZ}_L(\phi)^s\r] = \E\l[ \hat{\ZZ}_L(\phi)^2\r]^{\fr{s}{2}} \E\l[\chi^s 
  \r] + \E\l[ \hat{\ZZ}_L(\phi)^2\r]^{\fr{s}{2}} o(1)  \hsp{1.5cm} \mb{ as } L\tends \infty, 
\end{equation}
for $s=1,2,3,...$.
Then the moments of
 $\frac{\ZZ_L(\phi)-\E\l[\ZZ_L(\phi)\r]}{\sqrt{\Var(\ZZ_L(\phi))}}$
 converge to those of $\chi$ and convergence in distribution
 follows. To show  (\ref{eq:convofmoments}), we need to compute the moments of 
$\hat{\ZZ}_L(\phi)$.  


Recall the formula (\ref{eq:distributionallaplacian}) 
 \begin{equation}
   \intt_{\M} \phi(z) dn_{\F}(z) =  \intt_{\M} \lap \phi(z)
   \frac{1}{2\pi} \log |\F(z)| dm(z).
 \end{equation}
From this we can also write
\begin{equation*}
  \hat{\ZZ}_L(\phi)= \intt_{\M} \lap \phi(z)
   \frac{1}{2\pi} \log |\hat{\F}(z)| dm(z),
\end{equation*}
where $\hat{\F}(z)=\frac{1}{K(z,z)^{d/2}} \F(z)$ as defined in
  Section~\ref{sec:distribzeros}. Then  one can write the moments of
  $\hat{\ZZ}_L(\phi)$ as  
\begin{equation}\label{eq:centralmomentsofZL1}
\E\l[\hat{\ZZ}_L(\phi)^s \r] = (2\pi)^{-s} \intt_{\M^s}
 \l(\prodd_{j=1}^s \lap_{z_j}\phi(z_j) \r) 
 \E\l[\prodd_{j=1}^s \log  |\hat{\F}_L(z_j)|\r]
 \prodd_{j=1}^s dm(z_j)
\end{equation}
for $s=1,2,\ldots $.

\subsection{Central moments of $\ZZ_L(\phi)$}
From the homogeneity of $\Q$, we can write $\hat{\F}_L(z) =
\Q(\g_1^{(L)}(z),\ldots 
,\g_p^{(L)}(z))$, where $\g_i(z)=\frac{\f_i(z)}{\sqrt{K(z,z)}}$. These $\g_i^{(L)}$ are no longer analytic functions, but
 they are independent complex Gaussian processes on $\M$
 with constant variance $1$. Thus for any fixed $z$, we have that 
$\g_i^{(L)}(z)$, $1\le i\le p$ are i.i.d. standard complex Gaussians,
and from (\ref{eq:wickforP}) it follows that 
\begin{equation*}
  \log |\hat{\F}_L(z)| = \summ_{\m,\n\in \Z_+^p} \frac{C_{\m,\n}}{\sqrt{\m! \n!}} \prodd_{k=1}^p :\g_k^{m_k}(z)\bar{\g_k}^{n_k}(z):
\end{equation*}
with coefficients $C_{\m,\n}$ that are the same as in
(\ref{eq:wickforP}) and depend only on $\Q$ but not on $z$ or
even the GAF $\f$.  We get 
\begin{equation*}
\E\l[\prodd_{j=1}^s \log  |\hat{\F}_L(z_j)|\r] = \summ_{\l\{\m_j,\n_j\r\}_{1\le j\le s}}  \l( \prodd_{j=1}^s \frac{C_{\m_j,\n_j}}{\sqrt{\m_j! \n_j!}} \r) 
\prodd_{k=1}^p \E\l[\prodd_{j=1}^s :\g_k^{m_{kj}}(z_j)\bar{\g_k}^{n_{kj}}(z): \r].
\end{equation*}
Here $\m_j=(m_{1j},\ldots ,m_{pj})$ and likewise for $\n_j$. 

Each of the expectations on the right hand side of the above equation
can be ``evaluated'' by the Feynman diagram formula
(\ref{eq:feynmanformula}). We denote by 
$\hat{K}_L(z,w)$, the quantity $\E\l[ \g(z){\bar \g(w)}
\r]=\frac{K(z,w)} {\sqrt{K(z,z) K(w,w)}}$. Also we write $\ups(\gam
;z_1,\ldots ,z_s)$ for the value of a Feynman diagram $\gam$, with
edge weights given by the covariance matrix of the the Gaussian vector
$(\g(z_1),\ldots ,\g(z_s))$. Then, we get

\begin{eqnarray*}
\E\l[\prodd_{j=1}^s \log  |\hat{\F}_L(z_j)|\r] &=& \summ_{\l\{\m_j,\n_j\r\}_{1\le j\le s}}  \l( \prodd_{j=1}^s \frac{C_{\m_j,\n_j}}{\sqrt{\m_j! \n_j!}} \r) 
\summ_{\gam_1,\ldots ,\gam_p}  \prodd_{k=1}^p \ups(\gam_k; z_1,\ldots ,z_s) \\
 &=& \summ_{\gam_1,\ldots ,\gam_p} \l( \prodd_{j=1}^s \frac{C_{\m_j,\n_j}}{\sqrt{\m_j! \n_j!}} \r)  \prodd_{k=1}^p \ups(\gam_k; z_1,\ldots ,z_s) \\
&=& \summ_{\gam_1,\ldots ,\gam_p} \l( \prodd_{j=1}^s \frac{C_{\m_j,\n_j}}{\sqrt{\m_j! \n_j!}} \r)  \ups(\cup_{k=1}^p \gam_k; z_1,\ldots ,z_s). 
\end{eqnarray*}
Here in the last sum, $\gam_k$, $1\le k\le p$ vary over all possible
Feynman diagrams on labels $\{1,{\bar 1},\ldots s,{\bar s}\}$ and
$\m_j,\n_j \in \Z_+^p$ are  such that the number of vertices labeled
$j$ in $\gam_k$ is $m_{kj}$ and the number of vertices labeled ${\bar
  j}$ in $\gam_k$ is $n_{kj}$. 

Put this together with (\ref{eq:centralmomentsofZL1}) to  deduce that 
$\E\l[\hat{\ZZ}_L(\phi)^s \r]$ is equal to
\begin{equation}\label{eq:ZLmoment}
 (2\pi)^{-s}\summ_{\gam_1,\ldots ,\gam_p}
\l( \prodd_{j=1}^s \frac{C_{\m_j,\n_j}}{\sqrt{\m_j! \n_j!}} \r)
 \intt_{\M^s}  \l(\prodd_{j=1}^s \lap_{z_j} \phi(z_j)\r) 
\ups(\cup_{k=1}^p \gam_k; z_1,\ldots ,z_s)\prodd_{j=1}^s dm(z_j). 
\end{equation}
Here again, the sum is over all legal diagrams $\gam_k$.

Our goal is to prove (\ref{eq:convofmoments}). To that end, we now
consider the second moment, which is obtained by setting $s=2$. We get 
\begin{equation}\label{eq:secondmomentformula}
\E\l[\prodd_{j=1}^2 \log  |\hat{\F}_L(z_j)|\r] = \summ_{\gam_1,\ldots
  ,\gam_p} \l( \prodd_{j=1}^2 \frac{C_{\m_j,\n_j}}{\sqrt{\m_j! \n_j!}}
\r) \ups(\cup_{k=1}^p \gam_k; z_1,z_2), 
\end{equation} 
where, each $\gam_k$ is now a Feynman diagram on vertices labeled
$\{1,{\bar 1}, 2,{\bar 2}\}$.  

Now suppose $s$ is even.  Write (\ref{eq:secondmomentformula}) with
$z_1,z_2$ replaced by $z_{2k-1},z_{2k}$ for $k\le \frac{s}{2}$ and
multiply them together. On the right hand side, we get (the product of
the values  of Feynman diagrams is the value of the union of the Feynman
diagrams)  
\begin{equation*}
\summ_{\gam_1,\ldots ,\gam_p} \l( \prodd_{j=1}^s
\frac{C_{\m_j,\n_j}}{\sqrt{\m_j! \n_j!}}  \r) \ups(\cup_{k=1}^p
\gam_k; z_1,\ldots ,z_s), 
\end{equation*}
where the sum is over all Feynman diagrams $\gam_k$ on vertices
labeled by $\{1,{\bar 1},\ldots s,{\bar s}\}$, {\it such that}, when
all the vertices labeled by $j,{\bar j}$ are identified for each $j$,
then $\cup_{k=1}^p \gam_k$  has connected components 
$\{1,2\}$, $\{3,4\}, \ldots , \{s-1,s\}$.

To get $\E\l[ \ZZ_L(\phi)^2\r]^{\fr{s}{2}}$, 
integrate against $\fr{1}{(2\pi)^s}
\prodd_{j=1}^s \lap_{z_j} \phi(z_j)$ w.r.t Lebesgue measure over $\M^s$. This
yields that $\E\l[ \ZZ_L(\phi)^2 \r]^{\fr{s}{2}}$ is equal to
\begin{equation}\label{eq:diagramspliton12etc}
(2\pi)^{-s}
 \summ_{\gam_1,\ldots ,\gam_p} \l( \prodd_{j=1}^s
 \frac{C_{\m_j,\n_j}}{\sqrt{\m_j! \n_j!}}\r) \intt_{\M^s}
 \l(\prodd_{j=1}^s \lap_{z_j}\phi(z_j) \r) \ups
 (\cup_{k=1}^p \gam_k; z_1,\ldots ,z_s) \prodd_{j=1}^s dm(z_j). 
\end{equation}
Here again the sum is over all diagrams that have connected components
$\{1,2\}$, $\{3,4\}, \ldots ,\{s-1,s\}$ (upon merging vertices labeled
$j,{\bar j}$).  

Instead of pairing $\{1,2,\ldots ,s\}$ as $\{1,2\}, \ldots
,\{s-1,s\}$, we could use any other matching. Write the expression analogous
to (\ref{eq:diagramspliton12etc}) for each matching, and add them all
up. Recall that the 
number of matchings of $\{1,2,\ldots ,s\}$ is $\E\l[\chi^s \r]$. Thus
we deduce that $\E\l[ \ZZ_L(\phi)^2 \r]^{\fr{s}{2}} \E\l[\chi^s \r]$
is equal to  
\begin{equation}\label{eq:diagramsplit}
  (2\pi)^{-s}
 \summ_{\gam_1,\ldots ,\gam_p} \l( \prodd_{j=1}^s
 \frac{C_{\m_j,\n_j}}{\sqrt{\m_j! \n_j!}} \r) \intt_{\M^s}
  \l(\prodd_{j=1}^s \lap_{z_j}\phi(z_j) \r)
 \ups(\cup_{k=1}^p \gam_k; z_1,\ldots ,z_s) \prodd_{j=1}^s dm(z_j), 
\end{equation}
where the sum is over all diagrams $\gam_k$, $1\le k\le p$ such that
$\cup \gam_k$ ``splits'', i.e., when vertices $j,{\bar j}$ are
merged, we get $s/2$ component each of size $2$. 

Compare (\ref{eq:diagramsplit}) with (\ref{eq:ZLmoment}) (the
expressions are  incomplete without the commentaries that follows after
the equations!). The terms on the right hand side of
(\ref{eq:ZLmoment}) that are 
absent in (\ref{eq:diagramsplit}) are precisely those, for which
$\cup_{k=1}^p \gam_k$ does {\it not} split into $s/2$ components. The proof
of (\ref{eq:convofmoments}) will be complete once we show that these
terms together contribute a negligible amount compared to $\E\l[
\ZZ_L(\phi)^2 \r]^{\fr{s}{2}}$.

\subsection{Estimating the second moment}
This is the case $s=2$, which we already dealt with in detail, in
Chapter~\ref{chap:invzeros}. Particularly, from
(\ref{eq:2ptcorrelationformula}) we can write
\begin{equation}\label{eq:secondmomentZL}
\E\l[ \hat{\ZZ}_L(\phi)^2 \r] = \frac{1}{(2\pi)^2}
  \summ_{p=0}^{\infty} 
  |\tilde{C}_p|^2 \intt_{\M^2} \lap \phi(z) \lap \phi(w)
  |\hat{K}_L(z,w)|^{2p} dm(z) dm(w), 
\end{equation}
where $|\tilde{C}_p|^2= \summ_{\m_{\bullet}=\n_{\bullet}=p}
|C_{\m,\n}|^2$. It is convenient to express everything in terms of the
invariant quantities of $\M$. From (\ref{eq:invmeasure}) and
(\ref{eq:invlaplacian}) we rewrite (\ref{eq:secondmomentZL}) as
\begin{eqnarray*}
  \E\l[ \hat{\ZZ}_L(\phi)^2 \r] &=& 
  \summ_{p=0}^{\infty} \frac{|\tilde{C}_p|^2}{(2\pi)^2}
   \intt_{\M^2} \lap^* \phi(z) \lap^* \phi(w)
  |\hat{K}_L(z,w)|^{2p} dm^*(z) dm^*(w).
\end{eqnarray*}
Now split $\lap^* \phi(z) \lap^* \phi(w)$ as $(\lap^* \phi(z))^2 +
  (\lap^* \phi(w))^2 - (\lap^*   \phi(z)-\lap^* \phi(w))^2$. 
We get a sum of three integrals. 
The first two integrals are equal by symmetry. We shall argue that the
last integral is negligible. Firstly note that
\begin{equation}
  \label{eq:KLhat}
  |\hat{K}_L(z,w)|^2 = \l\{ \begin{array}{ll}
                 e^{-L|z-w|^2} & \mb{ if } \M=\C \\
                 \frac{|1+z\bar{w}|^{2L}}{(1+|z|^2)^L(1+|w|^2)^L} &
                 \mb{ if }\M=\S^2 \\ 
                 \frac{(1-|z|^2)^L(1-|w|^2)^L}{|1-z\bar{w}|^{2L}} 
                 & \mb{ if }\M=\D. \end{array} \r. 
\end{equation}
Thus $ |\hat{K}_L(z,w)|^2$ is $1$ when $z=w$ and decays rapidly as
$(z,w)$ moves away from the diagonal. Since $(\lap^*\phi(z)-\lap^*
\phi(w))^2$ vanishes on the diagonal, it is easy to calculate that 
\begin{equation}
  \label{eq:boundonsecondintegral}
  \intt_{\M^2} (\lap^*
  \phi(z)-\lap^* \phi(w))^2 |\hat{K}_L(z,w)|^{2p} dm^*(z) dm^*(w) =
  o\l( \frac{1}{Lp} \r). 
\end{equation}
The first two integrals give us
\begin{equation}
  \label{eq:boundonfirstintegral}
  \intt_{\M^2} (\lap^* \phi(z))^2 |\hat{K}_L(z,w)|^{2p} dm^*(z)
  dm^*(w). 
\end{equation}
Since $\hat{K}_L$ is also invariant under isometries of $\M$,
fixing $z$ and integrating w.r.t $w$ we get 
\begin{equation*}
 \l(\intt_M  |\hat{K}_L(0,w)|^{2p} dm^*(w) \r) \l( \intt_{\M} (\lap^*
  \phi(z))^2 dm^*(z) \r).
\end{equation*}
Now from (\ref{eq:KLhat}) it can be checked by direct computation that
\begin{equation}\label{eq:integralofKLhat}
  \intt_M  |\hat{K}_L(0,w)|^{2p} = \l\{ \begin{array}{ll}
                 \frac{\pi}{Lp} & \mb{ if } \M=\C. \\
                 \frac{\pi}{Lp+1} & \mb{ if }\M=\S^2. \\
                 \frac{\pi}{Lp-1} & \mb{ if }\M=\D.
                 \end{array} \r.
\end{equation}
From (\ref{eq:boundonsecondintegral}), (\ref{eq:boundonfirstintegral})
 and  (\ref{eq:integralofKLhat}), we get
\begin{equation}
  \label{eq:secmomofZLest}
   \E\l[ \hat{\ZZ}_L(\phi)^2 \r]= \l( \summ_{p=0}^{\infty}
   \frac{|\tilde{C}_p|^2}{4\pi L p} \r) \| \lap^* \phi \|_{L^2(m^*)}^2
   +  o\l( \frac{1}{L} \r).
\end{equation}
This shows (\ref{eq:varasymptotics}) with 
\begin{equation*}
  \kappa(\Q) = \summ_{p=0}^{\infty}
   \frac{|\tilde{C}_p|^2}{4\pi p}.
\end{equation*}

\subsection{Estimating the values of non-split diagrams}
We want to show that the contribution of non-split diagrams to $\E[
\hat{\ZZ}_L(\phi)^s ]$ is negligible. Note that this includes all the
diagrams in the case when $s$ is odd.
Consider any $p$-tuple of diagrams $(\gam_1,\ldots ,\gam_p)$ on labels
$\{1,{\bar 1},\ldots s,{\bar s} \}$ such that $\cup_k \gam_k$ does not
split into pairs when $j,\bar{j}$ are merged. We bound
\begin{equation*}
  \intt_{\M^s}
  \l(\prodd_{j=1}^s \lap_{z_j}\phi(z_j) \r)
 \ups(\cup_{k=1}^p \gam_k; z_1,\ldots ,z_s) \prodd_{j=1}^s dm(z_j),
\end{equation*}
in absolute value by the obvious 
\begin{equation*}
   \|\lap \phi \|_{\infty}^s 
 \intt_{A^s} |\ups(\cup_{k=1}^p \gam_k; z_1,\ldots ,z_s)|
 \prodd_{j=1}^s dm(z_j), 
\end{equation*}
where $A$ is the support of $\phi$. Then split the integral of
$|\ups(\cup_{k=1}^p \gam_k; z_1,\ldots 
,z_s)|$ as a product of integrals over the connected components of
$\cup_k \gam_k$. 

Without loss of generality, let $\{1,2,\ldots ,r\}$ be a connected
component of $\cup_k \gam_k$ (when $j,\bar{j}$ are merged). Since the
 weights of edges are bounded by $1$, deleting some edges will only
 increase the integral that we want to bound. We delete enough edges 
 to get a spanning tree $T$ on $\{1,2,\ldots ,r\}$. So we are left with an
 integral of the form 
\begin{equation*}
\intt_{A^r} \prodd_{(i,j)\in T} |\hat{K}_L(z_i,z_j)|
 dm(z_1)\ldots dm(z_r). 
\end{equation*}
Integrate inwards starting with the leaves. From
(\ref{eq:integralofKLhat}) we get that 
\begin{equation*}
  \intt_{A^r} \prodd_{(i,j)\in T} |\hat{K}_L(z_i,z_j)|
 dm(z_1)\ldots dm(z_r) <  CL^{-r}.
\end{equation*}
Multiplying the contribution from each component, we get
\begin{equation*}
   \intt_{A^s} \ups(\cup_{k=1}^p \gam_k; z_1,\ldots ,z_s)
 \prodd_{j=1}^s dm(z_j) < CL^{-s+\mb{no. of components of }\cup \gam_k} .
\end{equation*}
Non-split diagrams are precisely those that have less than
$\frac{s}{2}$ components, whence the right hand side is
$O(L^{-\frac{s-1}{2}})$. 

Now we want to bound the total contribution
of all unsplit diagrams. This can be done in the following manner.

Approximate $\log |\F(z)|$ by polynomials
(by truncating the Wick expansion). For integration of $\lap \phi$
against a polynomial, our cruder bound on a single Feynman diagram
suffices, since there are only finitely many terms and each of them goes to zero. So we get asymptotic normality for $\lap \phi$ integrated
against polynomials. From this, one can deduce asymptotic normality
for $\ZZ_L(\phi)$. We skip the details (see \cite{ST1}).

To summarize, we have argued that 
\begin{equation*}
  \E\l[\hat{\ZZ}_L(\phi)^s\r] - \E\l[
  \hat{\ZZ}_L(\phi)^2\r]^{\fr{s}{2}} \E\l[\chi^s \r] = o(L^{-s/2})
\end{equation*}
and proved that $\E\l[
  \hat{\ZZ}_L(\phi)^2\r]^{\fr{s}{2}}$ is of order $L^{-s/2}$. Thus
  (\ref{eq:convofmoments}) follows.


\chapter{Overcrowding Problems}\label{chap:overcrowd}
\section{Statements of the problems}
 In this chapter we go back to the canonical models of Gaussian
 analytic functions defined in Chapter~\ref{chap:intro}. Namely,
 consider the following 
 Gaussian analytic functions (GAFs): 
\begin{itemize}
\item {\bf Planar GAF :} The function defined in  (\ref{eq:planargaf}),
\[ \g(z)=\summ_{n=0}^{\infty} \frac{a_n z^n}{\sqrt{n!}} \]
where $a_n$ are i.i.d. standard complex Gaussian random variables.

\item {\bf Hyperbolic GAFs :} For each $L>0$ the function defined in
  (\ref{eq:hyperbolicgaf}) 
\[ \f_{L}(z) = \summ_{n=0}^{\infty} \binom{-L}{n}^{1/2} a_n z^n \]
where as before $a_n$ are i.i.d.  standard complex Gaussians. Almost
surely, $\f_{L}$ is an analytic function in the unit disk (and no
more).   
\end{itemize}
 We denote the zero set by ${\mathcal Z}$. Let $n(r)$ denote the
 number of points of ${\mathcal Z}$ in the disk of radius $r$ around
 $0$ (The GAF will be clear from the context). By the invariance of the
 zero sets, the results carry over to disks centered elsewhere.
 
Yuval Peres asked the following question
 and conjectured that the probability decays as $e^{-cm^2\log(m)}$ in
 the planar case (personal communication). 
 
\noindent {\bf Question:} Fix $r>0$, ($r<1$ in the Hyperbolic case). Estimate
$\P\l[n(r)>m \r]$ as $m\tends \infty$. 
\begin{figure}
\includegraphics[height=2.8in]{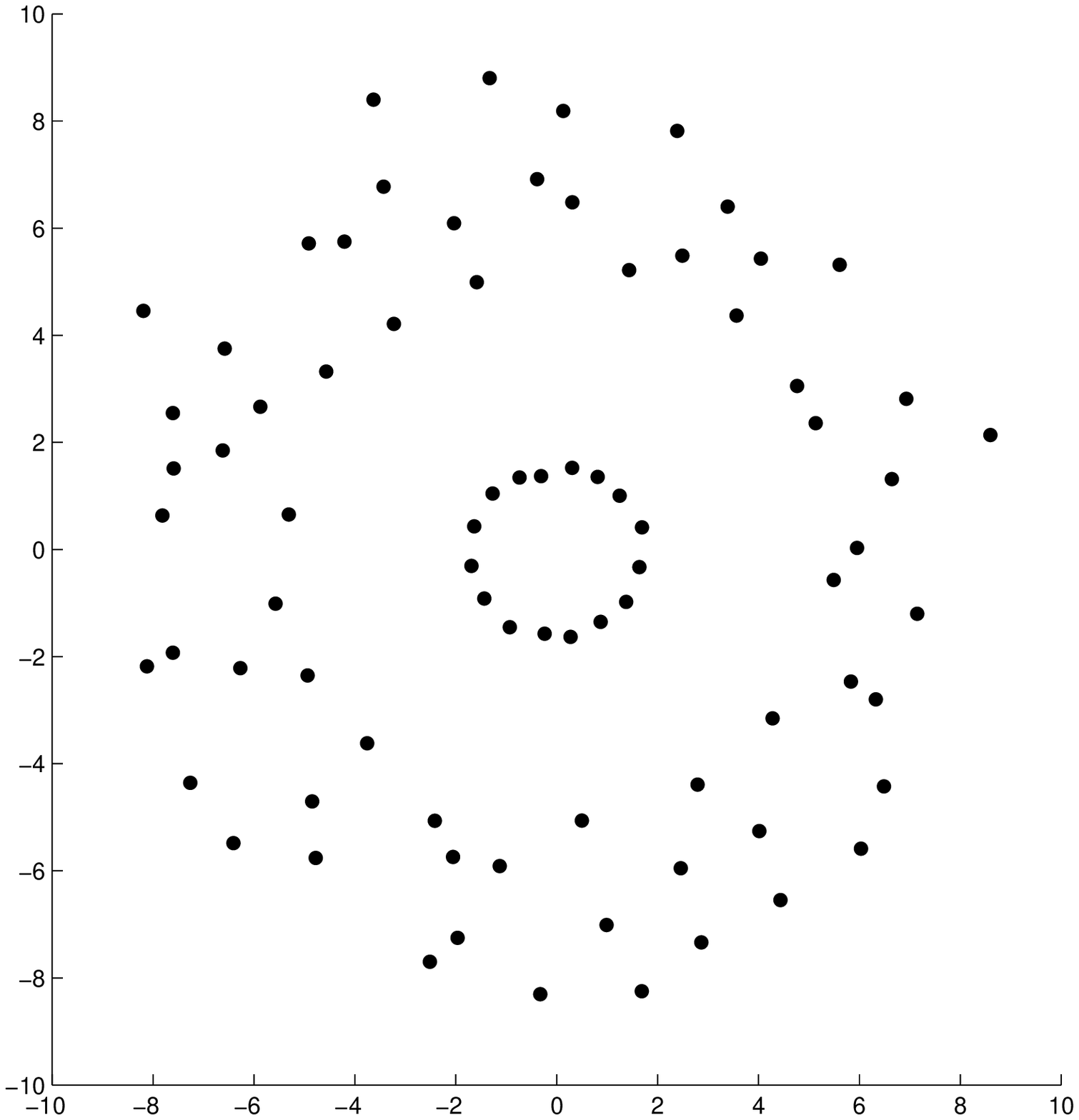}\hspace{.25in}
\includegraphics[height=2.8in]{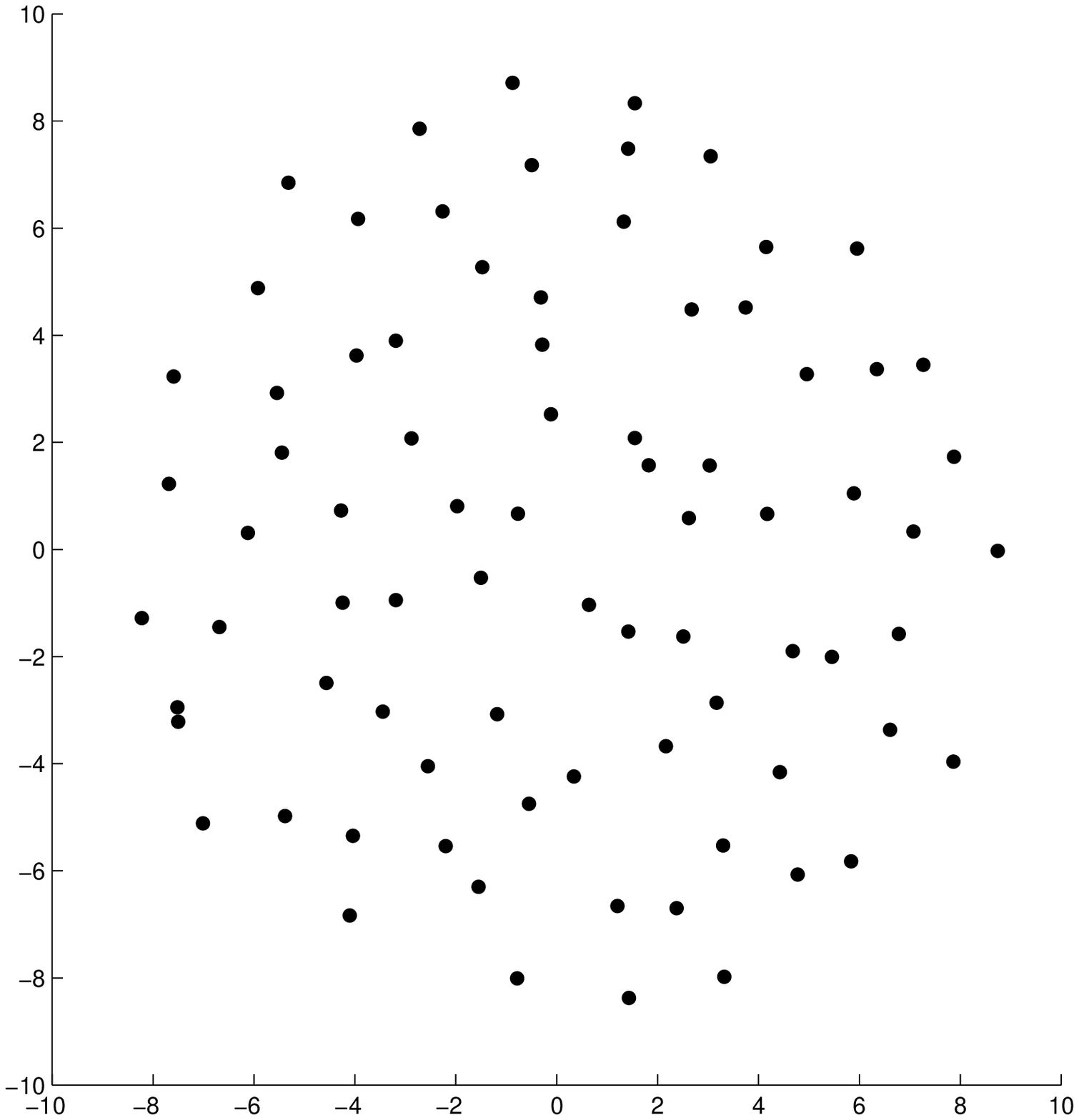}
  \centering
  \caption{Samples of the zero set of $\g$. Left: The zero process
    sampled under {\it certain sufficient
    conditions} (see the conditions in the lower bound of the proof of
    Theorem~\ref{thm:planar}. Take $\alp=0.5,r=2,m=16$) on the
    coefficients forcing $16$  zeros in the disk of radius
    $2$. Right: The unconditioned zero process.} 
  \label{fig:unconditioned}
\end{figure}

One motivation for such a question is in
Figure~\ref{fig:unconditioned}. There one can see the distribution of the 
zero set under certain conditions on the coefficients that force 
large number of zeros in the disk of radius $2$ (this is {\it not}
the zero set conditioned to have overcrowding - that seems harder to
simulate). The picture suggests that the distribution of the 
conditioned process may be worth studying on its own. A large
deviation estimate of the kind we derive will presumably be a
necessary step in such investigations.

The answer is different in the two settings. We prove-

\begin{theorem}\label{thm:planar} Consider the planar GAF $\g$. 
For any $\eps >0$, $\exists$ a constant $C_2$ (depending on $\eps,r$)
such that for every $m\ge 1$, 
\begin{equation*}
 e^{-\frac{1}{2} m^2\log(m)+O(m^2)} \le \P[n(r)\ge m] \le C_2
 e^{-(\frac{1}{2}-\eps)m^2\log(m)}. 
\end{equation*} 
In particular, $\P[n(r)\ge m]=e^{-\frac{1}{2} m^2\log(m)(1+o(1))}$.
\end{theorem}

\begin{theorem}\label{thm:hyperbolic} Fix $L>0$ and consider the
  GAF  $\f_{L}$. For any fixed $r<1$, there are constants
  $\bet,C_1,C_2$ (depending 
  on $L$ and $r$) such that for every $m\ge 1$, 
\begin{equation*} C_1(r)e^{- \frac{m^2}{|\log(r)|}} \le
  \P[n(r)\ge m] \le C_2(r) e^{-\bet(r) m^2}.
\end{equation*}
\end{theorem}

We prove Theorem~\ref{thm:planar} in Section~\ref{sec:planar},
Theorem~\ref{thm:hyperbolic} in Section~\ref{sec:hyperbolic}.

\section{Overcrowding - The planar case}\label{sec:planar}
In this section we prove Theorem~\ref{thm:planar}. Before that we
explain  why one expects the constant $\frac{1}{2}$ in the exponent in
Theorem~\ref{thm:planar}, by analogy with the Ginibre ensemble. 
 
\subsection{Ginibre ensemble} 
The Ginibre ensemble is the determinantal point process (earlier
we denoted this by Det-$\C-1$)  in the plane with kernel
\begin{equation}\label{eq:ginkernel}
 \Kdet(z,w)=\frac{1}{\pi}e^{-\frac{1}{2}|z|^2-\frac{1}{2}|w|^2+z{\bar w}}. 
\end{equation}
This process is of interest because it is the limit in distribution,
 as $n\tends\infty$, of the  
 point process of eigenvalues of an $n\times n$ matrix with
 i.i.d. standard complex Gaussian entries (Theorem~\ref{thm:ginibre}).

The Ginibre ensemble has many similarities to the zero set of $\g$. In
particular, the Ginibre ensemble is invariant in distribution under
Euclidean motions, 
has constant intensity $\frac{1}{\pi}$ in the plane and has the same
negative correlations as ${\mathcal Z}_{\g}$ at short
distances. Therefore there are other 
similarities too, for instance, see~\cite{denhan}. There are also
differences between the two point processes. For
instance, the Ginibre ensemble has all correlations negative, whereas
for the zero set of $\g$, long-range two-point correlations are
positive. However, in our problem, since we are considering a fixed
disk and looking at the event of having an excess of zeros in it, it
seems reasonable to expect the same behaviour for both these point
processes, since it is the short range interaction that is
relevant. In case of the Ginibre ensemble, the overcrowding problem is
easy to solve. 

\begin{theorem}\label{thm:ginibrecrowd}
Let $n_G(r)$ be the number of points of the Ginibre ensemble in the disk of radius $r$ around $0$ (by translation invariance, the same is true for any disk of radius $r$). Then for a fixed $r>0$,
\[ \P\l[n_G(r)\ge m\r] = e^{-\frac{1}{2}m^2\log(m) (1+o(1))}. \]
\end{theorem}
\begin{proof} By Kostlan~\cite{kostlan}, the set of absolute values of
  the points of the Ginibre ensemble has the same distribution as the
  set $\{R_1,R_2,\ldots \}$, where $R_n$ are independent, and $R_n^2$
  has Gamma($n,1$) distribution for every $n$. Hence $R_n^2\eqd
  \xi_1+\ldots +\xi_n$, where $\xi_k$ are i.i.d.  Exponential random
  variables with mean $1$, and it follows that 
\[ \P\l[R_n^2<r^2 \r]\ge \prodd_{k=1}^n \P\l[\xi_k<\frac{r^2}{n}\r]
  \ge \l(\frac{r^2}{2n}\r)^n, \] as long as $n\ge r^2$, because
  $\P\l[\xi_1<x\r]\ge \frac{x}{2}$ for $x<1$. Therefore we get
\begin{eqnarray} 
\P\l[n_G(r)\ge m\r] &\ge& \prodd_{n=1}^{m} \P\l[R_n^2<r^2 \r] \\
 &\ge& \prodd_{n=1}^{m} \l(\frac{r^2}{2n}\r)^n \\
 &=& \l(\frac{r^2}{2}\r)^{\frac{m(m+1)}{2}}e^{-{\summ_{n=1}^m
 n\log(n)}}.
\label{eq:lbfornr}
\end{eqnarray}
Here and elsewhere we shall encounter the term $\summ_{n=1}^m
n\log(n)$. We compute its asymptotics now. 
\[ n\log(n) \le x\log(x) \le (n+1)\log(n+1) \hsp{2cm}\mb{ for }n\le x\le n+1 \]
Integrate from $1$ to $m+1$ and note that 
\[ \intt_1^a x\log(x) dx = \frac{1}{2}a^2 \log(a) -
\frac{a^2}{4}+\frac{1}{4}, \] 
to get
\begin{equation}\label{eq:sumnlogn}
\summ_{n=1}^m n\log(n)\le  \frac{1}{2}(m+1)^2 \log(m+1) - \frac{(m+1)^2}{4}+ \frac{1}{4} \le \summ_{n=1}^{m+1} n\log(n). 
\end{equation}
Thus (\ref{eq:lbfornr}) gives
\begin{eqnarray*}
\P\l[n_G(r)\ge m\r] &\ge& e^{-\frac{1}{2}(m+1)^2\log(m+1)+\frac{(m+1)^2}{4}-\frac{1}{4}+\frac{m(m+1)}{2}\log(r^2/2)}\\
 &=& e^{-\frac{1}{2}m^2\log(m)+O(m^2)}.
\end{eqnarray*}
To prove the inequality in the other direction, note that
\begin{eqnarray*}
\P\l[n_G(r)\ge m\r] &\le& \P\l[\summ_{n=1}^{m^2}\ind(R_n^2<r^2)\ge m \r]+\summ_{n=m^2+1}^{\infty} \P\l[R_n^2<r^2\r] \\
&\le& {m^2 \choose m}\prodd_{n=1}^{m} \P\l[R_n^2<r^2 \r] +
\summ_{n>m^2}e^{-n\log(n)(1+o(1))}.
\end{eqnarray*}
In the second line, for the first summand we used the fact that $R_n^2$ are
stochastically increasing and for the second term we used the well known
fact $\P\l[R_n^2<r^2\r]=\P\l[\mb{Pois}(r^2)\ge n\r]$ and then the
usual bound on the tail of a Poisson random variable, namely $\P\l[\mb{Poisson}(\theta) \ge a \r]\le e^{-a\log(a/\theta)+a-\theta}$. 

Using the same idea to bound $\P\l[R_n^2<r^2\r]$ in the first summand,
we obtain
\begin{eqnarray*}
\P\l[n_G(r)\ge m\r] &\le& {m^2 \choose m}\prodd_{n=1}^{m} e^{-{n\log(n/r^2)-r^2+n}}+ e^{-m^2\log(m^2)(1+o(1))} \\
&\le& {m^2 \choose m}e^{\frac{m(m+1)}{2}(1+\log(r^2))-mr^2-\summ_{n=1}^m n\log(n)}+e^{-m^2\log(m^2)(1+o(1))} \\
&=& e^{-\frac{1}{2}m^2\log(m) (1+o(1))} \hsp{2cm} \l(\mb{using (\ref{eq:sumnlogn}) again} \r).
\end{eqnarray*}
In the last line we used ${m^2\choose m}<m^{2m}$. This completes the proof. 
\end{proof}

\subsection{Proof of Theorem~\ref{thm:planar}}
 Our method of proof is largely based on  that of Sodin and
 Tsirelson~\cite{ST3}. 
(They estimate the ``hole probability'', $\P\l[ n(r)=0\r]$ as $r\tends
 \infty$.) 

\begin{proof}[Proof of Theorem~\ref{thm:planar}]

\noindent {\bf Lower Bound } Suppose the $m^{\mb{th}}$ term dominates the sum of all the  other terms on $\d D(0;r)$, i.e., suppose
\begin{equation}\label{eq:dominate} 
\Mid\frac{a_m z^m}{\sqrt{m!}}\Mid \ge \Mid \summ_{n\neq m} \frac{a_n z^n}{\sqrt{n!}}\Mid \hspace{2cm} \mb{ whenever }|z|=r.
\end{equation}
Then, by Rouche's theorem $\g(z)$ and $\frac{a_m z^m}{\sqrt{m!}}$ have the same number of zeros in $D(0;r)$. Hence $n(r)=m$. 
Now we want to find a lower bound for the probability of the event in (\ref{eq:dominate}). Note that the left side of (\ref{eq:dominate}) is identically equal to $\frac{|a_m|r^m}{\sqrt{m!}}$. 

Now suppose the following happen-
\begin{enumerate}\label{enum:conditions}
\item $|a_n| \le n$ $\forall n\ge m+1$.
\item $|a_m| \ge (\alp+1)m$ where $\alp$ will be chosen shortly.
\item $|a_n|\frac{r^{n}}{\sqrt{n!}} < \frac{r^m}{\sqrt{m!}}$ for every $0\le n\le m-1$.
\end{enumerate}

Then the right hand side of (\ref{eq:dominate}) is bounded by
\begin{eqnarray*}
\mb{RHS of } (\ref{eq:dominate}) &\le& \summ_{n=0}^{m-1}|a_n|\frac{r^{n}}{\sqrt{n!}}  + \summ_{n=m+1}^{\infty} \frac{|a_n|r^n}{\sqrt{n!}} \\
&\le& \summ_{n=0}^{m-1}\frac{r^{m}}{\sqrt{m!}} + \summ_{n=m+1}^{\infty} \frac{nr^n}{ \sqrt{n!}} \\
&\le& m\frac{r^{m}}{\sqrt{m!}} + C \frac{m r^m}{\sqrt{m!}} \\
&=& (C+1)m\frac{r^{m}}{\sqrt{m!}} \\
&\le& |a_m|\frac{r^m}{\sqrt{m!}}
\end{eqnarray*}
if $\alp=C$. Thus if the above three events occur with $\alp=C$, then the $m^{\mb{th}}$ term dominates the sum of all the other terms on $\d D(0;r)$. Also these events have probabilities as follows.
\begin{enumerate}
\item $\P[|a_n| \le n$ $\forall n\ge m+1]\ge 1- \summ_{n=m+1}^{\infty} e^{-n^2}\ge 1-C'e^{-m^2}$. 
\item $\P[|a_m|\ge (C+1)m] = e^{-(C+1)^2m^2}$.
\item The third event has probability as follows. Recall again that
  $\P\l[\xi<x\r]\ge \frac{x}{2}$ if $x<1$ and $\xi$ has exponential with
  mean $1$. We apply this below with $x=\l(\frac{r^{m-n}\sqrt{n!}}{
    \sqrt{m!}}\r)^2$. This is clearly less than $1$ if $n\ge
  r^2$. Therefore if $m$ is sufficiently large it is easy to see that
  for all $0\le n\le m-1$, the same is valid. Thus
\begin{eqnarray*}
\P\l[|a_n|\le \frac{r^{m-n}\sqrt{n!}}{\sqrt{m!}} \hsp{2mm}\forall n\le
m-1 \r] &=& \prodd_{n=0}^{m-1} \P\l[|a_n|\le
\frac{r^{m-n}\sqrt{n!}}{\sqrt{m!}} \r] \\ 
&\ge& \prodd_{n=0}^{m-1} \frac{r^{2m-2n}n!}{2 m!} \\
&=& r^{m(m+1)}e^{\frac{1}{2}m^2\log(m)} 2^{-m} e^{-m^2\log(m)+O(m^2)} \\
&=& e^{-\frac{1}{2}m^2\log(m)+O(m^2)}.
\end{eqnarray*}
\end{enumerate}
Since these three events are independent, we get the lower bound in the theorem.

\noindent{\bf Upper Bound } By Jensen's formula, for any $R>r$ we have
\begin{equation}\label{eq:jensens}
n(r)\log\l(\frac{R}{r}\r) \le \intt_r^{R} \frac{n(u)}{u} du = \intt_0^{2\pi}\log|\g(Re^{i\theta })|\frac{d\theta}{2\pi} - \intt_0^{2\pi}\log|\g(re^{i\theta })|\frac{d\theta}{2\pi}.
\end{equation}

Let $R=R_m=\sqrt{m}$. Sodin and Tsirelson~\cite{ST3} show that 
\begin{equation}\label{eq:ubdformax}
\P\l[\log M(t)\ge\l(\frac{1}{2}+\eps\r)t^2\r]\le e^{-e^{\eps t^2}}
\end{equation}
where $M(t)=\max\{|\g(z)|: |z|\le t \}$. 

Now suppose $n(r)\ge m$ and $\log M(R_m)\le \l(\frac{1}{2}+\eps \r)m$ for some $\eps>0$. Then by (\ref{eq:jensens}) we have
\begin{eqnarray*}
-\intt_0^{2\pi}\log|\g(re^{i\theta})|\frac{d\theta}{2\pi} &\ge& m\log\l(\frac{\sqrt{m}}{r}\r) - \l(\frac{1}{2}+\eps\r)m \\
&=& \frac{1}{2} m\log(m) - m\log(r) -\l(\frac{1}{2}+\eps\r)m \\
&=& \frac{1}{2} m\log(m) - O(m).
\end{eqnarray*}
 Thus $\P[n(r)\ge m]$ is bounded by
\begin{eqnarray*}
 && \P\l[\log M(R_m)\ge \l(\frac{1}{2}+\eps\r)m\r] + \P\l[- \intt_0^{2\pi}\log|\g(re^{i\theta})|\frac{d\theta}{2\pi}\ge \frac{1}{2} m\log(m)-O(m) \r]\\
 &\le& e^{-e^{\eps m}} + \P\l[-
 \intt_0^{2\pi}\log|\g(re^{i\theta})|\frac{d\theta}{2\pi}\ge
 \frac{1}{2}m\log(m)(1+o(1)) \r] \hsp{1.2cm} \mb{by }\ref{eq:jensens}. 
\end{eqnarray*}
From Lemma~\ref{lem:ubdonintegral}, we deduce that for any $\del>0$, there is a constant $C_2$ such that  
\begin{eqnarray*} 
\P\l[- \intt_0^{2\pi}\log|\g(re^{i\theta})|\frac{d\theta}{2\pi} \ge \frac{1}{2} m\log(m) (1+o(1))\r] &\le& C_2 e^{-(2-\del)(\frac{m}{2} \log(m))^2/\log(\frac{m}{2} \log(m))}\\
&\le& C_2 e^{-(\frac{1}{2}-\frac{\del}{4})m^2\log(m) (1+o(1))}. 
\end{eqnarray*} From this, the upper bound follows. 
\end{proof}

\begin{lemma}\label{lem:ubdonintegral} For any given $\del>0$,
  $\exists C_2$ such that for every $m$,
\begin{equation*}
\P\l[- \intt_0^{2\pi}\log|\g(re^{i\theta})|\frac{d\theta}{2\pi}\ge m
\r]\le C_2e^{-\frac{(2-\del)m^2}{\log(m)}}.
\end{equation*}
 
\end{lemma}
\begin{proof} Let $P$ be the Poisson kernel on $D(0;r)$. Fix $\eps>0$
  and  let
 $A_{\eps}=\sup\{P(re^{i\theta},w):|w|=\eps,\theta\in[0,2\pi) \}$ and
 $B_{\eps}=\inf\{P(re^{i\theta},w):|w|=\eps,\theta\in[0,2\pi) \}$.
Since $\log |\g|$ is a subharmonic function, for any $w$ with $|w|=\eps$, we get
\begin{eqnarray*}
\log |\g(w)| &\le&
\intt_0^{2\pi}\log|\g(re^{i\theta})|P(re^{i\theta},w)\frac{d\theta}{2\pi}
\\ 
&\le& A_{\eps} \intt_0^{2\pi}\log_+|\g(re^{i\theta})|\frac{d\theta}{2\pi} - B_{\eps}\intt_0^{2\pi}\log_- |\g(re^{i\theta})|\frac{d\theta}{2\pi} \\
&=&  A_{\eps} \intt_0^{2\pi}\log_+|\g(re^{i\theta})|\frac{d\theta}{2\pi} + B_{\eps} \l( \intt_0^{2\pi}\log|\g(re^{i\theta})|\frac{d\theta}{2\pi}- \intt_0^{2\pi}\log_+|\g(re^{i\theta})|\frac{d\theta}{2\pi}\r) \\
&\le& B_{\eps}
\intt_0^{2\pi}\log|\g(re^{i\theta})|\frac{d\theta}{2\pi} + A_{\eps}
\log_+ M(r). 
\end{eqnarray*}
This implies $\log M(\eps) \le B_{\eps}
\intt_0^{2\pi}\log|\g(re^{i\theta})|\frac{d\theta}{2\pi} + A_{\eps}
\log_+ M(r)$.

Therefore if $\intt_0^{2\pi}\log|\g(re^{i\theta})|\frac{d\theta}{2\pi} \le -m$, then one of the following must happen. Either $\{\log M(\eps)\le -B_{\eps} m+\sqrt{m}$\} or  $\{A_{\eps}\log_+M(r)>\sqrt{m}\}$. 

Using (\ref{eq:ubdformax}), since $M(r)<M\l(Cm^{\frac{1}{4}}\r)$ for
any $C$, we see that 
\begin{equation*}
\P\l[\log_+M(r) >\frac{\sqrt{m}}{A_{\eps}}\r]\le e^{-e^{Cm}}
\end{equation*}
 for some  constant $C$ depending on $\eps$. Hence
\begin{eqnarray*}
\P\l[\intt_0^{2\pi}\log|\g(re^{i\theta})|\frac{d\theta}{2\pi} \le -m \r] &\le& e^{-e^{Cm}} + \P\l[ \log M(\eps) \le -B_{\eps} m+\sqrt{m} \r] \\
 &\le& e^{-e^{Cm}} + e^{-2B_{\eps}^2\frac{m^2}{\log(m)}(1+o(1))}
\end{eqnarray*}
where in the last line we have used Lemma~\ref{lem:maxmodulus}.

As $\eps\tends 0$, $B_{\eps}\tends 1$ and hence the proof is complete.

\end{proof}

Now we prove the upper bound on the maximum modulus in a disk of radius $r$ that was used in the last part of the proof of Lemma~\ref{lem:ubdonintegral}. For possible future use we prove a lower bound too.

\begin{lemma}\label{lem:maxmodulus} Fix $r>0$. There are constants $\alp,C_1,C_2$ such that 
\[ C_1 e^{-\frac{\alp m^2}{\log(m)}} \le \P[\log M(r)\le -m] \le C_2 e^{-\frac{2m^2}{\log(m)}(1+o(1))}.\]
\end{lemma}
\begin{proof}{\bf Lower bound } By Cauchy-Schwarz, 
$M(r)\le \l(\summ_{n=0}^{k-1} |a_n|^2 \r)^{1/2}e^{r^2/2} +
\summ_{n=k}^{\infty} \frac{|a_n|r^n}{\sqrt{n!}}$. We shall choose $k$
later. We will bound from below the probability that each of these
summands is less than $\frac{e^{-m}}{2}$.  

Let $\phi_k$ denote the density of $\Gam(k,1)$.
\begin{eqnarray*}
\P\l[\l(\summ_{n=0}^{k-1} |a_n|^2 \r)^{1/2}e^{r^2/2}\le
\frac{e^{-m}}{2}\r] &=&\P\l[\summ_{n=0}^{k-1} |a_n|^2 \le
\frac{e^{-2m}e^{r^2}}{4} \r]\\ 
&\ge& \phi_k\l(\frac{e^{-2m}e^{r^2}}{8}\r)\frac{e^{-2m}e^{r^2}}{8} \\
&=& e^{-2mk-k\log(k)+O(k)} 
\end{eqnarray*}

Also if $|a_n|\le n^2$ $\forall n\ge k$, then the second summand
\[ \summ_{n=k}^{\infty} |a_n|\frac{r^n}{\sqrt{n!}} \le C\frac{r^kk^2}{\sqrt{k!}} \le C e^{-k\log(k)/3} \]
Also the event $\{|a_n|\le n^2$ $\forall n\ge k \}$ has probability at least $1-\summ_{n=k+1}^{\infty} e^{-n^4} \ge 1-Ce^{-k^4}$.

Thus if we set $k=\frac{\gam m}{\log(m)}$ for a sufficiently large $\gam$, then both the terms are less than $e^{-\frac{m}{2}}$ with probability at least $e^{-2\gam m^2/\log(m)}$. 

\noindent{\bf Upper bound } By Cauchy's theorem, \[ a_n=\frac{\sqrt{n!}}{2\pi i}\intt_{C_r} \frac{\g(\zet)}{\zet^{n+1}} d\zet, \] where $C_r$ is the curve $C_r(t)=re^{it}$, $0\le t\le 2\pi$. Therefore, 
\[ |a_n| \le \frac{M(r)\sqrt{n!}}{r^n}. \]
Thus we get
\[ \P[M(r)\le e^{-m}] \le \prodd_{n=0}^{\infty} \P\l[|a_n|\le \frac{e^{-m}\sqrt{n!}}{r^n}\r]. \]
$|a_n|^2$ are i.i.d.  exponential random variables with mean
$1$. Therefore, 
\begin{equation*}
\P\l[|a_n|\le \frac{e^{-m}\sqrt{n!}}{r^n}\r]\le  \frac{e^{-2m}
 n!}{r^{2n}}.
\end{equation*}
 Using this bound for $n\le k :=\frac{\bet m}{\log(m)}$, we get
\begin{eqnarray*}
\P[M(r)\le e^{-m}] &\le& \prodd_{n=0}^k  \frac{e^{-2m} n!}{r^{2n}} \\
&\le& C e^{-2mk + \frac{k^2}{2}\log(k)+O(k^2)} \\
&\le& Ce^{(-2\bet+\frac{\bet^2}{2})\frac{m^2}{\log(m)}+O(\frac{m^2}{(\log(m))^2})}.
\end{eqnarray*}
$-2\bet+\frac{\bet^2}{2}$ is minimized when $\bet=2$ and we get,
\begin{equation}
 \P[M(r)\le e^{-m}] \le e^{-2\frac{m^2}{\log(m)}(1+o(1))}.  
\end{equation}
\end{proof}

\section{Overcrowding - The hyperbolic case}\label{sec:hyperbolic}

\subsection{The determinantal case}
We give a quick proof of Theorem~\ref{thm:hyperbolic} in the special
case $L=1$, as it is much easier and moreover we get matching upper
and lower bounds. The proof is similar to the case of the Ginibre
ensemble dealt with in Theorem~\ref{thm:ginibrecrowd} and is based on the
fact that the set of absolute values of the zeros of $\f_1$ is
distributed the same as a certain set of independent random
variables. The reason for this similarity between the two cases owes
to the fact that both of them are determinantal. The zero set of
$\f_1$ is a determinantal process with the 
Bergman kernel for the unit disk, namely 
\[ \Kdet_B(z,w)=\frac{1}{\pi}\frac{1}{(1-z{\bar w})^2}, \]
  as discovered by Peres and Vir{\'a}g~\cite{pervir}. 

\begin{proof}[Proof of Theorem~\ref{thm:hyperbolic} for $L=1$]
 By the result of Peres and Vir\'{a}g quoted in
 Theorem~\ref{thm:pervir}, the set of 
 absolute values of the zeros of $\f_1$ has the same distribution as
 the set $\{U_n^{1/2n}\}$ where $U_n$ are i.i.d. uniform$[0,1]$ random
 variables. Therefore, 
\begin{eqnarray*} 
\P\l[n(r)\ge m\r] &\ge& \prodd_{n=1}^{m} \P\l[U_n^{1/2n}<r \r] \\
 &=& \prodd_{n=1}^{m} r^{2n} \\
 &=& r^{m(m+1)}.
\end{eqnarray*}
To prove the inequality in the other direction, note that
\begin{eqnarray*}
\P\l[n(r)\ge m\r] &\le& \P\l[\summ_{n=1}^{m^2}\ind(U_n^{1/2n}<r)\ge m \r]+\summ_{n=m^2+1}^{\infty} \P\l[U_n^{1/2n}<r\r] \\
&\le& {m^2 \choose m}\prodd_{n=1}^m\P\l[U_n^{1/2n}<r \r] + \summ_{n>m^2}r^{2n}\\
&=& {m^2 \choose m}r^{m(m+1)}+ \frac{r^{2m^2+2}}{1-r^2}\\
&=& r^{m(m+1)}\l(1+ O\l(e^{m\log(m)}\r)\r).
\end{eqnarray*}
This completes the proof of the theorem for $L=1$. 
\end{proof}

\subsection{All values of $L$}


%



\begin{remark} Overall, the idea of proof is the same as in that of
  Theorem~\ref{thm:planar}. However we do not get matching upper and
  lower bounds in the present case, the reason being that in the
  hyperbolic analogue of Lemma~\ref{lem:maxmodulus}, the leading term
  in the exponent of the upper bound does depend on $r$, unlike in the
  planar case. (An examination of the proof of
  Theorem~\ref{thm:planar} reveals that we get a matching upper bound
  only because replacing $r$ by $\eps$  does not affect the leading
  term in the exponent in the upper bound in
  Lemma~\ref{lem:maxmodulus}). However we still expect that the lower
  bound in Theorem~\ref{thm:hyperbolic} is tight. (See remark after
  the proof).  
\end{remark}

\begin{proof}[Proof of Theorem~\ref{thm:hyperbolic}] 

\noindent{\bf Lower Bound }  As before we find a lower bound for the probability that the $m^{\mb{th}}$ term dominates the rest. Note that if $|z|=r$, 
\begin{equation}\label{eq:threeconditions}
\Mid\f_{L}(z)-{-L \choose m}^{1/2} a_mz^m \Mid \le \summ_{n=0}^{m-1} |a_n| {-L \choose n}^{1/2}r^n + \summ_{n=m+1}^{\infty} |a_n| {-L \choose n}^{1/2}r^n 
\end{equation}
Now suppose the following happen-
\begin{enumerate}
\item $|a_n| \le \sqrt{n}$ $\forall n\ge m+1$.
\item $|a_m| \ge (\alp+1)\sqrt{m}$ where $\alp$ will be chosen shortly.
\item $|a_n|{-L\choose n}^{1/2}r^n < \frac{1}{\sqrt{m}}{-L\choose m}^{1/2}r^m$ for every $0\le n\le m-1$.
\end{enumerate}
Then the right hand side of (\ref{eq:threeconditions}) is bounded by
\begin{eqnarray*}
\mb{RHS of } (\ref{eq:threeconditions}) &\le& \summ_{n=0}^{m-1}|a_n|{-L\choose n}^{1/2}r^n  + \summ_{n=m+1}^{\infty}|a_n|{-L\choose n}^{1/2}r^n \\
&\le& \summ_{n=0}^{m-1}\frac{1}{\sqrt{m}} {-L\choose m}^{1/2}r^m + \summ_{n=m+1}^{\infty}\sqrt{n}{-L\choose n}^{1/2}r^n \\
&\le& \sqrt{m} {-L\choose m}^{1/2}r^m + C \sqrt{m}{-L\choose m}^{1/2}r^m \hsp{1.5cm} \mbox{ for some }C \\
&=&  (C+1)\sqrt{m} {-L\choose m}^{1/2}r^m\\
&\le& |a_m| {-L\choose m}^{1/2}r^m
\end{eqnarray*}
if $\alp=C$. Thus if the above three events occur with $\alp=C$, then the $m^{\mb{th}}$ term dominates the sum of all the other terms on $\d D(0;r)$. Also these events have probabilities as follows. 
\begin{enumerate}
\item $\P[|a_n| \le \sqrt{n}$ $\forall n\ge m+1]\ge 1- \summ_{n=m+1}^{\infty} e^{-n}\ge 1-C'e^{-m}$. 
\item $\P[|a_m|\ge (\alp+1)\sqrt{m}] = e^{-(\alp+1)^2m}$.
\item The third event has probability as follows. Recall again that
  $\P\l[\xi<x\r]\ge \frac{x}{2}$ if $x<1$ and $\xi$ has exponential distribution with
  mean $1$. We apply this below with $x=\l(\frac{{-L\choose
  m}^{1/2}r^{m-n}}{\sqrt{m} {-L\choose n}^{1/2}}\r)^2$. This is
  clearly less than $1$. Thus 
\begin{eqnarray*}
\P\l[|a_n|\le \frac{ \binom{-L}{m}^{1/2}r^{m-n}}{\sqrt{m}\binom{-L}{n}^{1/2} } \hsp{2mm}\forall n\le m-1 \r] &=& \prodd_{n=0}^{m-1} \P\l[|a_n|\le \frac{ \binom{-L}{m}^{1/2}r^{m-n}}{\sqrt{m} \binom{-L}{n}^{1/2} } \r] \\
&\ge& \prodd_{n=0}^{m-1}\frac{\binom{-L}{m}r^{2m-2n}}{2 m \binom{-L}{n} } \\
&=& r^{m(m+1)} m^{-m} \prodd_{n=0}^{m-1} \frac{(m+1)\ldots (m+L-1)}{ (n+1)\ldots (n+L-1)} \\
&\ge& r^{m(m+1)} m^{-m} \prodd_{n=0}^{m-1}\frac{m^{L}}{(n+L)^{L}} \\
&\ge& r^{m(m+1)+O(m\log(m))}.
\end{eqnarray*}
\end{enumerate}
Since these three events are independent, we get the lower bound in the theorem. 

\noindent{\bf Upper Bound } The proof will proceed along the same lines as in
Theorem~\ref{thm:planar}. We need the following analogue of
Lemma~\ref{lem:maxmodulus}.  

\begin{lemma}\label{lem:hypmaxmod} Fix $r<1$.  Let $M(r)=\sup_{z\in D(0;r)} |\f_{L}(z)|$. Then
\[ \P[M(r)\le e^{-m}] \le e^{-\frac{m^2}{|\log(r)|}(1+o(1))}. \]
\end{lemma}
\begin{proof}
By Cauchy's theorem, for every $n\ge 0$,
\[ a_n{-L \choose n}^{1/2} = \frac{1}{2\pi i} \intt_{rT} \frac{\f(\zet)}{\zet^{n+1}} d\zet. \]
From this we get
\[ |a_n|^2 \le \frac{M(r)^2}{{-L \choose n} r^{2n}}. \]
Since ${-L \choose n}\ge \frac{n^{L-1}}{\Gam(L+1)}$, we obtain
\begin{eqnarray*}
\P[M(r)\le m] &\le& \prodd_n \P[|a_n|^2\le \frac{\Gam(L+1)e^{-2m}}{n^{L-1}r^{2n}}] \\
&\le& \prodd_{n=0}^\frac{m}{\log(1/r)}\frac{\Gam(L+1) e^{-2m}}{r^{2n}n^{L-1}} \\
&\le& e^{-\frac{2m^2}{\log(1/r)}+\l(\frac{m}{\log(1/r)}\r)^2\log(r)+O(m\log(m))}\\
&=& e^{-\frac{m^2}{\log(1/r)}+O(m\log(m))}. 
\end{eqnarray*}
\end{proof}

Coming back to the proof of the upper bound in the theorem, fix $R$ such that $r<R<1$. Then by Jensen's formula,
\begin{equation}\label{eq:jensenss}
n(r)\log\l(\frac{R}{r}\r) \le \intt_r^R \frac{n(u)}{u} du =
\intt_{R\T}\log|f(Re^{i\theta})| \frac{d\theta}{2\pi} - \intt_{r\T}
\log|f(re^{i\theta})| \frac{d\theta}{2\pi}.\end{equation}

Now consider the first summand in the right hand side of (\ref{eq:jensens}).
\begin{eqnarray*}
\P\l[\intt_{R\T}\log|\f(Re^{i\theta})|\frac{d\theta}{2\pi}
>\sqrt{m}\r] &\le& \P\l[\log M(R)\ge \sqrt{m}\r].
\end{eqnarray*}
Now suppose that $|a_n|<\lam^n$ $\forall n\ge m+1$ where
$1<\lam<1/R$. This has probability at least
$C_1e^{-\lam^{2m}/2}$. Then, 
\begin{eqnarray*}
M(R) &\le& \summ_{n-0}^{\infty} |a_n|{-L \choose n}^{1/2} R^n \\
&\le& \l(\summ_{n=0}^m |a_n|^2 \r)^{1/2}C_R + C_{R'}
\end{eqnarray*}
for some constants $C_R$ and $C_{R'}$.

Thus if $M(R)>e^{\sqrt{m}}$  then either $\summ_{n=0}^m |a_n|^2 >Ce^{2\sqrt{m}}$ or else $|a_n|>\lam^n$ for some  $n\ge m+1$. Thus
\begin{equation*}
\P\l[M(R)>\sqrt{m}\r] \le e^{-e^{c\sqrt{m}}}.
\end{equation*}
This proves that
\begin{equation*}
\P\l[\intt_{R\T}\log|\f(Re^{i\theta})|\frac{d\theta}{2\pi} >\sqrt{m}\r]\le e^{-e^{c\sqrt{m}}}.
\end{equation*}
Fix $\delta>0$ and $R$ close enough to $1$ such that
$\log(R)>-\delta$. Then with probability $\ge 1-e^{-e^{c\sqrt{m}}}$,
we obtain from (\ref{eq:jensenss}), 
\[ -\intt_{r\T} \log|\f(re^{i\theta})|\frac{d\theta}{2\pi} \ge m\l(\log\l(\frac{1}{r}\r) -\delta \r)-\sqrt{m}. \]

Now the calculations in the proof of Lemma~\ref{lem:ubdonintegral} show that
\[ \log M(\eps) \le B_{\eps} \intt_0^{2\pi} \log|\f(re^{i\theta})| P(re^{i\theta},w) \frac{d\theta}{2\pi}  + A_{\eps}\log_+ M(r). \]
Here $0<\eps<r$ is arbitrary and $A_{\eps},B_{\eps}$ are as defined in Lemma~\ref{lem:ubdonintegral}.  By the same computations as in that Lemma, we obtain,
we obtain the inequality
\begin{equation*} \P\l[\intt_0^{2\pi} \log |\f(re^{i\theta})|\frac{d\theta}{2\pi}\le
-m(|\log r|-\del)+\sqrt{m} \r] \le
 e^{-B_{\eps}^2\frac{m^2\log^2(r)(1-\delta)}{|\log (\eps)|}}+e^{-e^{cm}}.
\end{equation*}
Therefore, by (\ref{eq:jensenss})
\begin{equation*}
  \P\l[n(r)\ge m \r] \le e^{-\kappa m^2\log^2(r)(1+o(1))},
\end{equation*}
where $\kappa =\sup \l\{\frac{B_{\eps}^2}{|\log (\eps)|}: 0<\eps <r \r\}$.
 However it is clear that this cannot be made to match the lower bound
 by any choice of $\eps$.  
\end{proof}
 
\begin{remark} If we could prove   
\begin{equation*}\P\l[\intt_0^{2\pi} \log
  |\f(re^{i\theta})|\frac{d\theta}{2\pi}\le -x \r] \le
  e^{-\frac{x^2}{|\log(r)|}}, 
\end{equation*} 
that would have given us a matching upper bound. Now, one way for the
event $\intt_0^{2\pi} \log |\f(re^{i\theta})|\frac{d\theta}{2\pi}\le
-x$ to occur is to have $\log M(r)<-x$ which,  by
Lemma~\ref{lem:hypmaxmod} has probability at most
$e^{-x^2/\log(\frac{1}{r})}$. One way to proceed could be to show that
if the integral is smaller than $-x$, so is $\log M(s)$ for $s$
arbitrarily close to $r$ (with high probability). Alternately, if we
could bound the coefficients directly by the bound on the integral (as
in Lemma~\ref{lem:hypmaxmod}), that would also give us the desired
bound. For these reasons, and keeping in mind the case $L=1$, where we do have a
matching upper bound, we believe that the lower bound in
Theorem~\ref{thm:hyperbolic} is tight. 
\end{remark}


\chapter{Moderate  and very large deviations for zeros of the planar
GAF}\label{chap:moderate}

Inspired by the results obtained (using not entirely rigorous physical
arguments) by Jancovici, Lebowitz and 
Manificat~\cite{janlebmag} for Coulomb gases in the plane (eg.,
Ginibre ensemble), M.Sodin~\cite{sod2} conjectured the following. 

\begin{conjecture}[Sodin] Let $n(r)$ be the number of zeroes of the planar
GAF $\g$ in the disk $D(0,r)$. Then, as $r\tends \infty$ 
\begin{equation}\label{eq:conjecture}
\frac{\log \log \l(\frac{1}{\P[|n(r)-r^2|>r^{\alp}]}\r)}{\log r} \tends \l\{ \begin{array}{cc}
         2\alp -1, & \frac{1}{2}\le \alp \le 1; \\
         3\alp -2, & 1\le \alp \le 2; \\
         2\alp,    & 2\le \alp. \end{array} \r. 
\end{equation}
\end{conjecture}
The idea here is that the deviation probabilities undergo a
qualitative change in behaviour when the deviation under consideration
becomes comparable to the perimeter ($\alp=1$) or to the area
($\alp=2$) of the domain.

Sodin and Tsirelson~\cite{ST3} had already settled the case
$\alp=2$ by showing that for any $\del>0$, $\exists
c_1(\del),c_2(\del)$ such that 
\begin{equation*}
  e^{-c_1(\del)r^4} \le \P[|n(r)-r^2|>\del r^2] \le e^{-c_2(\del)r^4}.
\end{equation*}

Here we consider $\P[n(r)-r^2>r^{\alp}]$ and prove that a ``phase
transition'' in the exponent occurs at  $\alp=2$. More precisely we
prove that the conjecture holds for $\alp>2$ and show the lower bound
for $1<\alp<2$. 

\begin{theorem}\label{thm:verylarge} Fix $\alp>2$. Then 
\[ \P\l[n(r)\ge r^2+\gam r^{\alp} \r] = e^{-\l(\frac{\alp}{2} -
  1\r)\gam^2 r^{2\alp}\log r (1+o(1))}. \] 
\end{theorem}

 \begin{theorem}\label{thm:moderate} Fix $1<\alp < 2$. Then  for
   any $\gam>0$, 
\[ \P\l[n(r)\ge r^2+\gam r^{\alp} \r]\ge e^{-\gam^3 r^{
    3\alp-2}(1+o(1))}. \] 
\end{theorem}

We prove Theorem~\ref{thm:verylarge} in Section~\ref{sec:verylarge}
and Theorem~\ref{thm:moderate} in Section~\ref{sec:moderate}. 
Taken together these show that the asymptotics
of $\P\l[ n(r)\ge r^2+\gam r^{\alp}\r]$ does undergo a qualitative change at
$\alp=2$.
\begin{remark}  Nazarov, Sodin and Volberg have recently proved all
  the remaining 
  parts of the conjecture (personal communication).
\end{remark}

\section{Very large deviations for the planar GAF}\label{sec:verylarge}
In this section we prove Theorem~\ref{thm:verylarge}.
\begin{remark}In the case $\alp\ge 2$, one side of the estimate as
  asked for in 
  the conjecture (with $\log \log$ of the probability) follows
  trivially from the results in  
Sodin and Tsirelson~\cite{ST3}. They prove that for any $\del>0$,
there exists a constant $c(\del)$ such that 
\[ \P\l[ |n(r)-r^2| >\del r^2 \r] \le e^{-c(\del)r^4}. \]
When $\alp\ge 2$, clearly $n((1-\del)r^{\sqrt{\alp}})\ge n(r)$, whence
from the above result it follows that  
\begin{eqnarray*}
\P\l[n(r)\ge r^2+r^{\alp} \r] &\le& \P\l[n((1-\del)r^{\sqrt{\alp}})\ge r^{\alp} \r]\\
&\le& e^{-c(\del)r^{2\alp}}.
\end{eqnarray*}
This gives 
\begin{equation}\label{eq:mildbd}
\limsup_{r\tends \infty} \frac{\log \log
         \l(\frac{1}{\P[|n(r)-r^2|>r^{\alp}]}\r)}{\log r} \le  2\alp.
\end{equation}
\end{remark}

The obviously loose inequality $n((1-\del)r^{\sqrt{\alp}})\ge n(r)$
that we used, suggests that (\ref{eq:mildbd}) can be improved when $\alp>2$
to Theorem~\ref{thm:verylarge}. 

\begin{proof}[Proof of Theorem~\ref{thm:verylarge}]

\noindent{\bf Lower Bound }Let $m=r^2+\gam r^{\alp}$. Suppose the $m^{\mb{th}}$ term dominates the sum of all the  other terms on $\d D(0;r)$, i.e., suppose
\begin{equation}\label{eq:dominate2} 
\Mid\frac{a_m z^m }{\sqrt{m!}}\Mid \ge \Mid \summ_{n\neq m} \frac{a_n z^n}{\sqrt{n!}}\Mid \hspace{2cm} \mb{ whenever }|z|=r.
\end{equation}

Now we want to find a lower bound for the probability of the event in (\ref{eq:dominate2}). Note that the left side of (\ref{eq:dominate2}) is identically equal to $\frac{|a_m|r^m}{\sqrt{m!}}$. 

Now suppose the following happen-
\begin{enumerate}
\item $|a_n| \le n$ $\forall n\ge m+1$.
\item $|a_m| \ge m$.
\item $|a_n|\frac{r^{n}}{\sqrt{n!}} < \frac{\gam r^{\alp}}{m}\frac{r^m}{\sqrt{m!}}$ for every $0\le n\le m-1$.
\end{enumerate}
Then the right hand side of (\ref{eq:dominate2}) is bounded by
\begin{eqnarray*}
\mb{RHS of } (\ref{eq:dominate2}) &\le& \summ_{n=0}^{m-1}|a_n|\frac{r^{n}}{\sqrt{n!}}  + \summ_{n=m+1}^{\infty} \frac{|a_n|r^n}{\sqrt{n!}} \\
&\le& \summ_{n=0}^{m-1}\frac{\gam r^{\alp}}{m}\frac{r^{m}}{ \sqrt{m!}} + \summ_{n=m+1}^{\infty} \frac{nr^n}{\sqrt{n!}} \\
&\le& \frac{mr^{m}}{\sqrt{m!}}\l( \frac{\gam r^{\alp}}{ m} + o(1) \r) \\
&\le& |a_m|\frac{r^m}{m!}
\end{eqnarray*}

 Thus if the above three events occur, then the $m^{\mb{th}}$ term dominates the sum of all the other terms on $\d D(0;r)$. Also these events have probabilities as follows.
\begin{enumerate}
\item $\P[|a_n| \le n$ $\forall n\ge m+1]\ge 1- \summ_{n=m+1}^{\infty} e^{-n^2}\ge 1-C'e^{-m^2}=1-o(1)$. 
\item $\P[|a_m|\ge m] = e^{-m^2}=e^{-\gam^2 r^{2\alp}(1+o(1))}$.
\item The third event has probability as follows. Recall again that
  $\P\l[\xi<x\r]\ge \frac{x}{2}$ if $x<1$ and $\xi$ is exponential with
  mean $1$. We apply this below with $x=\frac{\gam^2 r^{2\alp}}{m^2}\frac{r^{2m-2n}n!}{
    m!}$. This is clearly less than $1$ if $n\ge
  r^2$. Therefore if $m$ is sufficiently large it is easy to see that
  for all $0\le n\le m-1$, the same is valid. Thus
\begin{eqnarray*}
\P\l[|a_n|\le \frac{\gam r^{\alp}}{m}\frac{r^{m-n}\sqrt{n!}}{\sqrt{m!}} \hsp{2mm}\forall n\le m-1 \r] &=& \prodd_{n=0}^{m-1} \P\l[|a_n|\le \frac{\gam r^{\alp}}{m}\frac{r^{m-n}\sqrt{n!}}{\sqrt{m!}} \r] \\
&\ge& \prodd_{n=0}^{m-1} \frac{\gam^2 r^{2\alp}}{m^2}\frac{r^{2m-2n}n!}{2 m!} \\
&=& r^{2\alp (m+1)+m(m+1)} 2^{-m}m^{-2m}e^{-\summ_{k=1}^m k\log k}  \\
&=& e^{m^2\log(r)-\frac{1}{2}m^2\log(m)+O(m^2)}\\
&=& e^{-(\frac{\alp}{ 2}-1)\gam^2 r^{2\alp} \log(r) +O(r^{2\alp})}
\end{eqnarray*}

\end{enumerate}

Since these three events are independent, we get 
\begin{equation}
\P\l[n(r)\ge r^2+r^{\alp} \r] \ge e^{-\l(\frac{\alp}{ 2}-1\r)\gam^2 r^{2\alp}\log r + O(r^{2\alp})}.
\end{equation}
\noindent {\bf Upper Bound} We omit the proof of the upper bound, as it follows 
the same lines as that of Theorem~\ref{thm:planar} and we have
already seen such arguments again in the proof of
Theorem~\ref{thm:hyperbolic} (In those two cases as well as the
present case, we are looking at very large deviations, and that is the
reason why the same tricks work). 

Moreover note that the lower bound along with (\ref{eq:mildbd}) proves the
statement in the conjecture. 
\end{proof}

\section{Moderate deviations for the planar GAF}\label{sec:moderate}
In this section we prove Theorem~\ref{thm:moderate}. 

\begin{proof}[Proof of Theorem~\ref{thm:moderate}] Write $m=r^2+\gam
  r^{\alp}$. As usual, we bound $\P\l[n(r)\ge m\r]$ from below by the
  probability of the event that the $m^{\mb{th}}$ term dominates the
  rest of the series.

Firstly, we need a couple of estimates. Consider $\frac{r^{2n}}{n!}$
as a function of $n$. This increases 
monotonically  up to $n=r^2$ and then decreases
monotonically. $m=r^2+\gam r^{\alp}$ is on the latter part. Write
$M=r^2-\gam r^{\alp}$.

Firstly,  observe that $(r^2-k)(r^2+k)<(r^2)^2$, for  $1\le k\le
\gam r^{\alp}$, whence $r^{2m-2M}>\prodd_{j=M+1}^{m-1} j$. This
implies that 
\begin{equation}\label{eq:relationMm}
\frac{r^M}{\sqrt{M!}}<\frac{r^m}{\sqrt{m!}}
\end{equation}

Secondly, note that for any $n=M-p$, 
\begin{eqnarray*}
  \frac{r^{2n}/n!}{r^{2M}/M!} &=&\prodd_{j=0}^{p-1} \frac{M-j}{r^2} \\
             &=& \prodd_{j=0}^{p-1} (1-\gam r^{\alp-2}-jr^{-2}) \\
             &\le& e^{-\summ_{j=0}^{p-1} (\gam r^{\alp-2}+jr^{-2})} \\
             &=& e^{-\gam pr^{\alp-2} -\frac{p(p+1)}{2}r^{-2}}.
\end{eqnarray*}
Now we set $p=Cr^{2-\alp}$ with $C$ so large that 
$ e^{-\gam C}\le \frac{1}{4}$.

Then also note that if $n<M-kp$, it follows that
\begin{equation}
  \label{eq:decayofcoeffs}
  \frac{r^{2n}/n!}{r^{2m}/m!}\le \frac{1}{4^k},
\end{equation}
where we used (\ref{eq:relationMm}) to replace $M$ by $m$.

Thirdly, if $n=m+p$ with $p\le r^2-\gam r^{\alp}$, then,
\begin{eqnarray*}
  \frac{r^{2n}/n!}{r^{2m}/m!} &=&\prodd_{j=1}^{p} \frac{r^2}{m+j} \\
             &=& \prodd_{j=1}^{p} (1+\gam r^{\alp-2}+jr^{-2})^{-1} \\
             &\le& e^{-\frac{1}{2}\summ_{j=1}^{p} (\gam
             r^{\alp-2}+jr^{-2})} \\ 
             &=& e^{-\frac{1}{2}(\gam p r^{\alp-2} +\frac{p(p+1)}{2}r^{-2})}.
\end{eqnarray*}
If $p=2Cr^{2-\alp}$, where $C$ was as chosen before, then for
$n>m+kp$, we get   
\begin{equation}
  \label{eq:decayofcoeffs2}
  \frac{r^{2n}/n!}{r^{2m}/m!}\le \frac{1}{4^k}.
\end{equation}
From now on $p=2Cr^{2-\alp}$ is fixed so that (\ref{eq:decayofcoeffs})
and (\ref{eq:decayofcoeffs2}) are satisfied.

Next we divide the coefficients other than the $m^{\mb{th}}$ one into groups: 
\begin{itemize}
\item $A_k=\{n:n\in (M-kp,M-(k-1)p]$ for $1\le k\le \lceil \frac{M}{p}
  \rceil \}$.
\item $D_k=\{n:n\in [m+(k-1)p,m+kp)$ for $1\le k\le \lceil \frac{M}{p}
  \rceil$\}.
\item $B=\{n:n\in [M+1,m-1]\}$.
\item $C=\{n:n\in [2r^2,\infty)\}$.
\end{itemize}

\begin{remark}
  As defined, there is an overlap between $D_{\lceil \frac{M}{p}
  \rceil}$ and $C$. This is inconsequential, but for definiteness, let
  us truncate the former interval at $r^2$ (just as $A_{\lceil \frac{M}{p}
  \rceil}$ is understood to be truncated at $0$).
\end{remark}

Now consider the following events.

\begin{enumerate}
\item $|a_n| \le \frac{2^k}{M}$ for $n\in A_k$ for $k\le \lceil
  \frac{M}{p}\rceil$ \}.
\item $|a_n|\le \frac{2^k}{M}$ for $n\in D_k$ for $k\le
  \lceil \frac{M}{p}\rceil \}$.
\item $\summ_{n\in B} |a_n|\frac{r^n}{\sqrt{n!}} \le 4 \frac{r^m}{\sqrt{m!}}$.
\item $|a_n|<n-2r^2$ for $n\in C$.
\item $|a_m| \ge 15$.
\end{enumerate}

Suppose all these events occur. Then

\begin{enumerate}
\item The event $|a_n| \le \frac{2^k}{M}$ for $n\in A_k$, $k\le \lceil
  \frac{M}{p}\rceil$ gives
\begin{eqnarray}\label{eq:A}
  \sup\{\Mid \summ_{n=0}^M \frac{a_nz^n}{\sqrt{n!}} \Mid:|z|=r \} &\le&
  \summ_{k=1}^{\lceil M/p \rceil} \summ_{n=M-kp+1}^{M-(k-1)p}
  |a_n|\frac{r^n}{\sqrt{n!}} \\
     &\le&  \summ_{k=1}^{\lceil M/p
  \rceil}\frac{1}{2^{k}}\frac{r^{m}}{\sqrt{m!}}\frac{2^kp}{M}
  \hsp{1cm} \mb{by }(\ref{eq:decayofcoeffs})\\
 &\le& \frac{r^{m}}{\sqrt{m!}} \summ_{k=1}^{\lceil M/p \rceil}\frac{p}{M} \\
   &\le& \frac{r^{m}}{\sqrt{m!}} \l(1+\frac{p}{M}\r). 
\end{eqnarray}

\item The event $|a_n|\le \frac{2^k}{M}$ for $n\in D_k$, $k\le
  \lceil \frac{M}{p}\rceil$ gives
\begin{eqnarray}\label{eq:D}
\sup\{\Mid \summ_{n=m+1}^{2r^2} \frac{a_nz^n}{\sqrt{n!}}\Mid:|z|=r \} &=& 
\summ_{k=1}^{\lceil M/p \rceil} \summ_{n=m+(k-1)p+1}^{M+kp}
  |a_n|\frac{r^n}{\sqrt{n!}} \\
 &\le&\summ_{k=1}^{\lceil M/p
  \rceil}\frac{1}{2^{k}}\frac{r^{m}}{\sqrt{m!}}\frac{2^kp}{M}
  \hsp{1cm} \mb{by }(\ref{eq:decayofcoeffs2}) \\
 &\le& \frac{r^m}{\sqrt{m!}} \l(1+\frac{p}{M}\r). 
\end{eqnarray}

\item The third event gives 
\begin{equation}\label{eq:B}
\summ_{n\in B} |a_n|\frac{r^n}{\sqrt{n!}} \le 4 \frac{r^m}{\sqrt{m!}},
\end{equation}
by assumption.

\item The event $|a_n|<n-2r^2$ for $n\in C$: Since $n>2r^2$,
\begin{eqnarray*}
\frac{r^n}{\sqrt{n!}} &=& \frac{r^m}{\sqrt{m!}}\prodd_{k=m+1}^{n} \frac{r}{\sqrt{k}} \\
&\le&  \frac{r^m}{\sqrt{m!}}\prodd_{k=2r^2+1}^{n} \frac{r}{\sqrt{k}} \\
&\le& \frac{r^m}{\sqrt{m!}}\l(\frac{1}{\sqrt{2}}\r)^{n-2r^2}.
\end{eqnarray*}
Therefore we get (using $|a_n|<n-2r^2$ $\forall n>2r^2$)
\begin{eqnarray}
\summ_{n\in C} |a_n|\frac{r^n}{\sqrt{n!}} &\le& \frac{r^m}{\sqrt{m!}} \summ_{n>2r^2} (n-2r^2) \l(\frac{1}{\sqrt{2}}\r)^{n-2r^2} \\
 &=& \frac{r^m}{\sqrt{m!}} \frac{\sqrt{2}}{(\sqrt{2}-1)^2}.
\end{eqnarray}
\end{enumerate}

Putting together the contributions from these four groups of terms,
and using $|a_m|>15$, we get (for large values of $r$)
\begin{eqnarray*}
\summ_{n\not=m}|a_n|\frac{r^n}{\sqrt{n!}} \le |a_m|\frac{r^m}{\sqrt{m!}}.
\end{eqnarray*}

Now we compute the probabilities of the events enumerated above.

\begin{enumerate}
\item The event $|a_n| \le \frac{2^k}{M}$ for $n\in A_k$ for $k\le \lceil
  \frac{M}{p}\rceil$.
Now for a fixed $k\le 3\log_2(r)$, we deduce
\begin{eqnarray*}
  \P\l[|a_n|\le \frac{2^k}{M} \mb{ for } n\in A_k \r] &\ge&
  \P\l[|a_n|\le \frac{1}{M}  \mb{ for } n\in A_k \r] \\
   &\ge& \l(  \frac{1}{2M^2} \r)^{p}.
\end{eqnarray*}
Therefore
\begin{eqnarray*}
  \P\l[|a_n|\le \frac{2^k}{M} \mb{ for } n\in A_k \mb{
  for every }k\le 3\log_2(r)\r] &\ge& \l(  \frac{1}{2M^2}
  \r)^{3p\log_2(r)}\\
&\ge& e^{-cr^{2-\alp}(\log(r))^2}.
\end{eqnarray*}
for some $c$.

Next we deal with $k> 3\log_2(r)$.

\begin{eqnarray*}
  \P\l[|a_n|\le \frac{2^k}{M} \mb{ for } n\in A_k \mb{
  for every }k>3\log_2(r)\r] &\ge& 1- \summ_{k>3\log_2(r)}p \P\l[|a|>
  \frac{2^k}{M}\r]  \\
 &=& 1- \summ_{k>3\log_2(r)}p e^{-2^{2k}M^{-2}}.
\end{eqnarray*}
Now the summation in the the last line has rapidly decaying terms and 
starts with $pe^{-2^{6\log_2(r)}M^{-2}}$ which is smaller than
$pe^{-r^2}$. Thus 
\begin{equation*}
  \P\l[|a_n|\le \frac{2^k}{M} \mb{ for } n\in A_k \mb{
  for every }k>3\log_2(r)\r] = 1-o(1).
\end{equation*}

Thus the event in question has probability at least $e^{-
  cr^{2-\alp}(\log(r))^2(1+o(1))}$.



\item The event $|a_n|\le \frac{2^k}{M}$ for $n\in D_k$ 
  for $k\le \lceil \frac{M}{p}\rceil$. Following exactly the same
  steps as above we can prove that
\begin{equation*}
\P\l[|a_n|\le \frac{2^k}{M} \mb{ for } n\in D_k \r]
  \ge e^{- cr^{2-\alp}(\log(r))^2(1+o(1))}.
\end{equation*}

\item The event $\summ_{n\in B} |a_n|\frac{r^n}{\sqrt{n!}} \le 4
  \frac{r^m}{\sqrt{m!}}$.  By Cauchy-Schwarz, 
\begin{equation*}
\l(\summ_{n\in B} |a_n|\frac{r^n}{\sqrt{n!}}\r)^2 \le \l( \summ_{n\in
  B} |a_n|^2 \r) \l(\summ_{n\in B} \frac{r^{2n}}{n!}
\r).
\end{equation*}
 $Y=\summ_{n\in B} |a_n|^2$ has $\Gam(|B|,1)$ distribution. Also 
$\summ_{n\in B} \frac{r^{2n}}{n!} \le e^{r^2}$, since the left hand
is part of the Taylor series of $e^{r^2}$.
Therefore the event in question has probability, 
\begin{eqnarray*}
\P[\mb{event in question}] &\ge& \P\l[Y< 16 \frac{r^{2m}}{m!}e^{-r^2}\r] \\
             &\ge& \varphi\l(8 \frac{r^{2m}}{m!}e^{-r^2}\r)8
             \frac{r^{2m}}{m!}e^{-r^2}, 
\end{eqnarray*} 
where $\varphi$ is the density of the $\Gam(|B|,1)$ distribution.  This last follows because $\varphi$ is increasing on $[0,|B|]$ and thus $\P[Y<x]\ge \varphi\l(\frac{x}{2}\r)\frac{x}{2}$. for $x<|B|$.
Continuing,

\begin{eqnarray}
\P[\mb{event in question}] &\ge&\frac{1}{(2\gam r^{\alp})!} e^{-8
  \frac{r^{2m}}{m!}e^{-r^2}}\l(8 \frac{r^{2m}}{m!}e^{-r^2}\r)^{2\gam
  r^{\alp}} \\ 
&\ge& C e^{2\gam r^{\alp} (m\log
  (r^2)-r^2 -m\log (m)+m)+ O(r^{\alp}\log(r)) }.
\label{eq:probofB} 
\end{eqnarray}

where we used Stirling's approximation.

The exponent needs simplification. Take the first and third terms in
the exponent. We have $-2\gam
mr^{\alp}\log\l(\frac{m}{r^2}\r)$. Recall that $m=r^2+\gam r^{\alp}$
and that $\alp<2$. Therefore by Taylor's expansion of $\log(1+\gam
r^{\alp-2})$ we get 
\begin{equation}\label{eq:stirl}
-2\gam mr^{\alp}\log \l(\frac{m}{r^2}\r) = \l\{  \begin{array}{ccc}
           -2\gam^2 r^{2\alp} &+ \gam^3 r^{3\alp-2} &-\frac{2}{3} \gam^4  r^{4\alp-4}+\ldots\\
           -2\gam^3 r^{3\alp-2} &+\gam^4 r^{4\alp-4} &-\frac{2}{3} \gam^5 r^{5\alp-6}+\ldots \end{array} \r.
\end{equation}
Now consider (\ref{eq:probofB}). Expand the fourth term in the 
exponential as $2\gam r^{2+\alp}+2\gam^2 r^{2\alp}$. We get the following
terms 
\begin{enumerate}
\item $r^{2+\alp}(-2\gam+2\gam)=0$, from the second and fourth terms
  (first piece of the fourth term) in the exponential in (\ref{eq:probofB}).
\item $r^{2\alp}(-2\gam^2 + 2\gam^2)=0$, from the sum of the first
  term in the expansion (\ref{eq:stirl}) and the second piece of the
  fourth term in the exponential  in (\ref{eq:probofB}).  
\item $r^{3\alp-2}(\gam^3-2\gam^3)=-\gam^3r^{3\alp-2}$, from the
  expansion (\ref{eq:stirl}). 
\item Other terms such as
  $r^{\alp}\log(m),r^{\alp}\log(r),r^{\alp},r^{4\alp-4},r^{5\alp-6}$
  etc. All these are of lower order than $r^{3\alp-2}$ when $1<\alp<2$.
\end{enumerate}
Hence,
\[ \P[\mb{event in quesion}] \ge e^{-\gam^3 r^{3\alp-2}(1+o(1))}. \]
 
\item The event $|a_n|<n-2m$ for $n\in C$. This is just an event for a sequence of i.i.d.  complex Gaussians. It has a fixed probability $p_0$ (say). 
\item The event $|a_m| \ge 15$ also has a constant probability (not
  depending on $r$, that is).

\end{enumerate}
This completes the estimation of probabilities. Among these five
events, the third one, namely $\summ_{n\in B} |a_n|\frac{r^n}{n!} \le
4 \frac{r^m}{\sqrt{m!}}$ has the least probability (Recall that $1 <
 \alp < 2$).  

Also these events are all independent, being dependent on disjoint
sets of coefficients. Thus $\P\l[n(r)\ge r^2+\gam r^{\alp} \r]\ge
e^{-\gam^3r^{3\alp-2}(1+o(1))}$.  

\end{proof} 
\newpage
\addcontentsline{toc}{chapter}{Bibliography}
\bibliography{refs}


\end{document}